\documentclass[11pt]{article}
\usepackage{amsmath,amsfonts}
\usepackage{xcolor}
\usepackage{amsthm}         
\usepackage{dsfont}			
\usepackage{hyperref}
\usepackage{multirow}
\usepackage{graphicx}
\usepackage{subcaption}
\usepackage{epsfig}
\usepackage{stmaryrd}
\usepackage{mathtools}
\usepackage{amsrefs}
\usepackage{soul}
\usepackage{float}
\def\cb{{\bf c}}
\def\tb{{\bf t}}
\def\Xb{{\boldsymbol X}}
\def\Yb{{\boldsymbol Y}}

\def\Ukd{\mathcal{U}_{k,d}}

\usepackage{colortbl}
\definecolor{regionone}{RGB}{255,232,228}
\definecolor{regiontwo}{RGB}{229,240,252}
\definecolor{regionthree}{RGB}{228,244,230}

\newlength{\resultwidth}
\newcommand{\tablezero}{\makebox[\resultwidth][c]{0}}

\newcommand*{\Ztu}{{}_{n}\mathring{\mathbb{Z}}_{u_k}}
\newcommand*{\thetab}{{\boldsymbol{\theta}}}
\newcommand*{\etab}{{\boldsymbol{\eta}}}

\newcommand*{\lb}{{\boldsymbol{\ell}}}
\newcommand*{\E}[1]{\operatorname{E}_{#1}}

\newcommand*{\Prob}[2]{\operatorname{P}_{#1} \left( {#2} \right)}

\newcommand*{\deq}{{\, \stackrel{\text{def}}{=} \; }}
\newtheorem{theorem}{Theorem}

\hypersetup{
	pdfborder = {0 0 0},
	urlbordercolor = {0 0 0},
	colorlinks = {false},
}

\def\underacc #1/#2{\mathchoice{\uacc\textstyle{#1}{#2}}{\uacc\textstyle{#1}{#2}}
	{\uacc\scriptstyle{#1}{#2}}{\uacc\scriptscriptstyle{#1}{#2}}}
\def\uacc#1#2#3{\mathop{#2{}}\limits_{#1#3{}}}

\newcommand*{\ind}[1]{\mathds{1} \left( {#1} \right)}

\begin{document}
	
	\title{Adaptive almost full recovery in sparse nonparametric models}
	\author{ Xiang Zhao$^a$,
		Marie Turcicova$^{b0}$, and Natalia A. Stepanova$^a$ \\
		\\
		\small {\textit{$^a$School of Mathematics and Statistics, Carleton University, 4302 Herzberg Laboratories,}}\\
		\small {\textit{1125 Colonel By Drive, Ottawa, ON, K1S 5B6, Canada} } \\ 
		\small{\textit{$^b$Institute of Computer Science, Czech Academy of Sciences, }} \\
		\small{\textit{Pod Vodárenskou věží 271/2, Prague, 182 00,  Czech Republic}}}
	\date{} 
	
	\footnotetext{Corresponding author. Email address: turcicova@cs.cas.cz}
	\maketitle

	\begin{abstract}
		We observe an unknown function of $d$ variables $f(\tb)$, $\tb\in[0,1]^d$, in the Gaussian white noise model of intensity $\varepsilon>0$.
	We assume that the function $f$ is regular and that it is a sum of $k$-variate functions, where $k$ varies from $1$ to $s$ ($1\leq s\leq d$).
	These functions are unknown to us and only a few of them are nonzero.
	In this article, we address the problem of identifying the nonzero function components of $f$ almost fully  in the case when $d=d_\varepsilon\to \infty$ as $\varepsilon\to 0$ and $s$ is  either fixed or $s=s_\varepsilon\to \infty$, $s=o(d)$ as $\varepsilon\to 0$. This may be viewed as a variable selection problem.
	We derive the conditions when almost full variable selection in the model at hand is possible and provide a selection procedure
	that achieves this type of selection. The procedure is adaptive to the level of sparsity
	described by the sparsity index $\beta\in(0,1)$.
	We also derive conditions that  make almost full variable selection in the model of our interest impossible.
	In view of these conditions, the proposed selector is seen to be asymptotically optimal.
	The theoretical findings are illustrated numerically.
		
		\bigskip
		\noindent \textit{Keywords:} Gaussian white noise, functional ANOVA model, sparsity, almost full selection, asymptotically minimax selector
		
		\noindent \textit{2020 MSC}: Primary: 62G08, Secondary: 62H12, 62G20
	\end{abstract}
	
\section{Introduction} \label{sec:problem_statement}
In this work, we address the problem of sparse signal recovery
in a nonparametric regression model in continuous time, also known as the Gaussian white noise model,
and augment the results on exact variable selection obtained in \cite{ST-2023} and \cite{ST-2024}.
Specifically, we assume that an unknown signal $f$ of $d$ variables is observed in the Gaussian white noise model
\begin{gather}\label{model0}
	dX_\varepsilon(\tb)=f(\tb)d\tb+\varepsilon d W(\tb),\quad \tb\in[0,1]^d,
\end{gather}
where $dW$ is a $d$-parameter Gaussian white noise and $\varepsilon>0$ is the noise intensity.
The signal 	$f$ belongs to a subspace of $L_2([0,1]^d)=L_2^d$ with an inner product  $(\cdot,\cdot)_2$ and a norm $\|\cdot\|_2$
that consists
of regular enough functions, and we assume high-dimensional settings when $d=d_\varepsilon\to \infty$ as $\varepsilon\to 0$.
Consider an operator $\mathbb{W}:L_2^d\to {\cal G}_0$ taking values in
the set ${\cal G}_0$ of centered Gaussian random variables
such that if $\xi_0=\mathbb{W}(g_1)$ and $\eta_0=\mathbb{W}(g_2)$, where $g_1,g_2\in L_2^d$,
then ${\rm cov}(\xi_0,\eta_0)=(g_1,g_2)_{2}$. The $d$-parameter Gaussian white noise $dW$ in model~(\ref{model0}) is defined through the operator $\mathbb{W}$ by
$$\mathbb{W}(g)=\int_{[0,1]^d} g(\tb) dW(\tb)\sim N(0,\|g\|_2^2),\quad g\in L_2^d.$$
In particular, if $\{g_\lb\}_{\lb\in{\cal L}}$, $\mathcal{L} \subseteq \mathbb{Z}^d$, is an  orthonormal basis of $L_2^d$, then $\mathbb{W}(g_\lb)\sim N(0,1)$ for $\lb\in{\cal L}$
and, for any finite set $\{g_{\lb}\}$ of the basis functions, the family $\{\mathbb{W}(g_{\lb})\}$ forms a multivariate standard normal vector. Thus,
the centered Gaussian measure on $L_2^d$ determined by $\mathbb{W}$ has a diagonal covariance operator (i.e., the identity operator).
Furthermore, let $\mathbb{X}_\varepsilon:L_2^d\to {\cal G}$ be an operator taking values in the set
${\cal G}$  of Gaussian random variables such that if $\xi=\mathbb{X}_{\varepsilon}(g_1)$ and $\eta=\mathbb{X}_{\varepsilon}(g_2)$, where $g_1,g_2\in L_2^d$,
then ${\rm E}(\xi)=(f,g_1)_2$, ${\rm E}(\eta)=(f,g_2)_2$, and ${\rm cov}(\xi,\eta)=\varepsilon^2(g_1,g_2)_2$.
By ``observing the trajectory (\ref{model0})'', we mean observing
a realization of the Gaussian field $X_\varepsilon(\tb)$, $\tb\in[0,1]^d$,  defined
through the operator $\mathbb{X}_\varepsilon$
by
\begin{gather*} 
	\mathbb{X}_\varepsilon(g)=\int_{[0,1]^d}g(\tb)dX_\varepsilon(\tb)\sim N\Big((f,g)_2,\varepsilon^2\|g\|_2^2  \Big),\quad  g\in L_2^d.
\end{gather*}
In terms of the operators $\mathbb{W}$ and $\mathbb{X}_\varepsilon$, the stochastic differential equation (\ref{model0}) can be expressed as
\begin{equation}\label{model1}
	\mathbb{X}_{\varepsilon}=f+\varepsilon \mathbb{W},
\end{equation}
and ``observing the trajectory (\ref{model1})'' means that we observe all normal $N\left((f,g)_2,\varepsilon^2\|g\|_2^2  \right)$ random variables
when $g$ runs through $L_2^d$. For any $f\in L_2^d$, the ``observation'' $\mathbb{X}_\varepsilon$ in model~(\ref{model1})
defines the Gaussian measure ${\rm P}_{\varepsilon,f}$ on the Hilbert space $L_2^d$ with mean function $f$ and covariance operator
$\varepsilon^2I$, where $I$ is the identity operator (for references, see \cite{GN-2016, IH.97, SK}).
In addition to regularity constraints, we assume that $f$ has a sparse structure and consider the problem of recovering almost fully the
sparsity pattern of $f$ from the ``observation'' $\mathbb{X}_{\varepsilon}$  by using the asymptotically minimax approach.

\subsection{Sparsity conditions}\label{SCond}
To avoid the curse of dimensionality  stemming from high-dimensional settings,
we assume that $f$ has a sparse structure.
The notion of sparsity employed in this work will be formalized by assuming a~sparse functional ANOVA expansion for $f$ proposed in \cite{ST-2024}. 
The problem of our interest is to  ``almost fully'' recover the sparsity pattern of $f$ when
$d=d_\varepsilon\to \infty$ as $\varepsilon\to 0$. A precise definition of almost full recovery will be given in Section \ref{PS}, 
after the necessary notation and framework are introduced.

	To justify the choice of sparse functional models studied in this article, let us first consider the ANOVA expansion of a function $f\in L_2^d$ 
(see, for example, \cite{ISu-2008} and \cite{Lin-2000}).
	To this end, for $1\leq k\leq d$, denote by $\Ukd$ the set of all subsets $u_k\subseteq\{1,\ldots,d\}$ of cardinality $k$, that is,
	\begin{gather*} \label{def:Ukd}
		\Ukd=\{u_k:u_k\subseteq\{1,\ldots,d\}, \#(u_k)=k\},
	\end{gather*}
	and note that $\#\left( \mathcal U_{k,d}\right)={d\choose k}.$
We put $u_k=\emptyset$ and  ${\cal U}_{k,d}=\emptyset$ for $k=0.$
		For ${u_k}\in\Ukd$, $1\leq k\leq d,$ consider the set
	\begin{eqnarray*}
				\mathring{\mathbb{Z}}_{u_k}&=&\{\lb=(l_1,\ldots,l_d)\in \mathbb{Z}^d: l_j=0\;\mbox{for}\;j\notin u_k \mbox{ and }  l_j\neq 0\;\mbox{for}\;j\in u_k \}.
	\end{eqnarray*}
	The set $\mathring{\mathbb{Z}}_{u}$ for $u\subseteq\{1,\ldots,d\}$ is defined analogously.
	We also put $\mathring{\mathbb{Z}}_\emptyset=\underbrace{(0,\ldots,0)}_d$, $\mathring{\mathbb{Z}}=\mathbb{Z}\setminus\{0\},$
		and note that
		$\mathbb{Z}^d=\left( \mathring{\mathbb{Z}}\cup\{0\} \right)^d=\bigcup_{u\subseteq\{1,\ldots,d\}}\mathring{\mathbb{Z}}_{u}=
	\bigcup_{0\leq k\leq d}\bigcup_{u_k\in \Ukd}\mathring{\mathbb{Z}}_{u_k}.$
	The space $L_2^d$ has a tensor product structure $L_2^d=L_2^1\bigotimes\ldots \bigotimes L_2^1$ and the 
	tensor product Fourier basis $\{\phi_{\lb}(\tb)\}_{\lb\in{\mathbb{Z}}^d}$, defined as follows:
\begin{gather*}
		\phi_{\lb}(\tb)=\prod_{j=1}^d \phi_{l_j}(t_{j}),\quad \lb=(l_1,\ldots,l_d)\in {\mathbb{Z}}^d,
		\\
		\phi_0(t)=1,\quad \phi_l(t)=\sqrt{2}\cos(2\pi l t),\quad \phi_{-l}(t)=\sqrt{2}\sin(2\pi l t),\quad l>0.
	\end{gather*}
Note that $L_2^1={\cal L}_0\bigoplus{\cal L}_1$, where ${\cal L}_0$ consists of constants and ${\cal L}_1$ is its orthogonal complement,
which implies that $L_2^d=({\cal L}_0\bigoplus{\cal L}_1)^d$.
	
Using the notation  $\tb_{u_k}={(t_j)}_{j\in u_k}$ for  $u_k\in \Ukd$, observe that
	$\phi_{\lb}(\tb)=\phi_{\lb}(\tb_{u_k})$ for $\lb\in\mathring{\mathbb{Z}}_{u_k}$.
	Denote also $\tb_{u}={(t_j)}_{j\in u}$ for $u\subseteq\{1,\ldots,d\}$. Then we can write
	$$\{\phi_{\lb}(\tb)\}_{\lb\in\mathbb{Z}^d}=\bigcup_{u\subseteq\{1,\ldots,d\}}\{\phi_{\lb}(\tb_{u})\}_{\lb\in	\mathring{\mathbb{Z}}_{u}}=
	\bigcup_{0\leq k\leq d}\bigcup_{u_k\in \Ukd} \{\phi_{\lb}(\tb_{u_k})\}_{\lb\in	\mathring{\mathbb{Z}}_{u_k}}.$$
The tensor product structure in
$L_2^d=({\cal L}_0\bigoplus{\cal L}_1)^d$  leads to the orthogonal expansion
	\begin{gather}\label{eqq1}
		L_2^d=\bigoplus_{u\subseteq\{1,\ldots,d\}} {\cal L}_u^d,
	\end{gather}
	where $ {\cal L}_u^d$ consists of functions $f_u(\tb_u)$ that depend on $\tb_u$, and the sum is taken over all $2^d$ subsets $u$ of the set $\{1,\ldots,d\}$.
	The Fourier basis of ${\cal L}_u^d$ is $\{\phi_\lb\}_{\lb\in\mathring{\mathbb{Z}}_u }$.
	Clearly, if $f\in L_2^d$, then 
	$$f(\tb)=\sum_{u\subseteq\{1,\ldots,d\}} f_u(\tb_u) ,\quad \int_0^1f_u(\tb_u)\,dt_j=0\quad \mbox{ a.e. for all}\;j\in u,\; u\subseteq\{1,\ldots,d\},$$
where $f_u={\rm const}$ for $u=\emptyset$.
	Decomposition (\ref{eqq1}) can be expressed in the form of  ``ANOVA expansion'' as follows:
	\begin{equation*} 
		L_2^d=\bigoplus_{k=0}^d {\cal L}_k^d,\quad {\cal L}_k^d=\bigoplus_{u_k\in \Ukd}{\cal L}_{u_k}^d,
	\end{equation*}
	i.e., the subspace ${\cal L}_k^d$ consists of sums of functions that depend on $k$ variables,
	yielding the representation for $f\in L_2^d$ 
	in the form 
	\begin{equation}\label{function}
		f(\tb)=f_\emptyset+\sum_{k=1}^d \sum_{u_k\in \Ukd}f_{u_k}(\tb_{u_k}),\quad \int_0^1f_{u_k}(\tb_{u_k})\,dt_j=0\;\;
			\mbox{ a.e. for all}\;j\in u_k, \;u_k\in \Ukd,\; 1\leq k\leq d.
	\end{equation}
		From (\ref{function}) and the fact that $f_\emptyset=\mathrm{const}$, we have $f_\emptyset=\int_{[0,1]^d} f(\tb)\,d\tb$.
	The tensor product space ANOVA modelling, which extends the additive methodology in a natural way, is capable of dealing with interactions of all orders in a flexible way.

We now consider two sparse decompositions of $f$ in model (\ref{model1}) of our interest.  
	First, extending the sparse additive model studied in \cite{BSt-2017}, we assume that, for a given $k$ ($1\leq k\leq d$),
	an unknown (regular) signal $f\in L_2^d$ satisfies
\begin{gather}\label{f}
	f(\tb)=\sum_{u_k\in\Ukd}\eta_{u_k} f_{u_k}(\tb_{u_k}),\quad \tb\in[0,1]^d,
\end{gather}
where  each component  $f_{u_k}$, $u_k\in\Ukd$, is such that
\begin{gather}\label{orthcon}
	\int_{0}^1 f_{u_k}(\tb_{u_k})\,dt_j=0\quad  \mbox{  for all}\; j\in u_k,
\end{gather}
and the $\eta_{u_k}$s are unknown but deterministic quantities taking values in $\{0,1\}$:
$\eta_{u_k}=0$ (or \textit{irrelevant} $\eta_{u_k}$) means that the component $f_{u_k}$ is inactive, whereas $\eta_{u_k}=1$ (or \textit{relevant} $\eta_{u_k}$) means that the component $f_{u_k}$ is active.
The number $\sum_{u_k\in\Ukd  }\eta_{u_k}$ of active components
is set to be small compared to the total number of components ${d\choose k}$, specifically
	\begin{gather}\label{sparsitycond}
		\sum_{u_k\in\Ukd}\eta_{u_k} \leq C {d\choose k}^{1-\beta},
	\end{gather}
where $\beta\in(0,1)$ is the \textit{sparsity index}, and $1<C<\infty$ is a constant that is assumed to be fixed
throughout the article.
 Further, introduce the sets $H^{\beta}_{k,d}$, $1\leq k\leq d$, as follows:
	\begin{equation*}
H^{\beta}_{k,d} = \left\{ \boldsymbol{\eta}_k=( \eta_{u_k})_{ u_k \in \Ukd}: \eta_{u_k} \in \{0,1\}\;\mbox{and condition}\; (\ref{sparsitycond})\;\mbox{holds}\right\}.
\end{equation*}
The orthogonality conditions in (\ref{orthcon}) imply that if  $u_k\neq v_k$ are subsets of $\Ukd$, then $f_{u_k}(\tb_{u_k})$ and $f_{v_k}(\tb_{v_k})$
are orthogonal (in $L_2^d$) to each other and to a constant, which guarantees uniqueness of representation~(\ref{f}).
The signal $f$ as in (\ref{f}) is \textit{sparse} because the majority of the components $f_{u_k}$  are inactive,
and at most $C {d\choose k}^{1-\beta}$ components are active,
where ${d\choose k}^{1-\beta}=o\left( {d\choose k}\right)$ as $d\to \infty$, $k$~is either fixed or $k=o(d)$.
In other words,  $f$~is the sum of a small number of $k$-variate functions. The values of $\beta$ that are close to one make the signal $f$ in (\ref{f})
\textit{highly sparse}, with a very few components $f_{u_k}$ on the right side of (\ref{f}) being active,
whereas the values of $\beta$ that are close to zero make it \textit{dense}.

A more general problem of sparse signal recovery, in which for some $s$ ($1\leq s\leq d$) an unknown (regular) signal $f\in L_2^d$
observed in the Gaussian white noise model (\ref{model1})  has the form, cf. (\ref{function}) and (\ref{f}),
\begin{equation}
	f(\tb)=f_{\emptyset}+\sum_{k=1}^{s} \sum_{u_k\in\Ukd}\eta_{u_k} f_{u_k}(\tb_{u_k}),\quad \tb\in [0,1]^d,\quad\boldsymbol{t}_{u_k} =(t_j)_{j \in u_k} \in [0,1]^k,\;1\leq k\leq s,
	\label{fun_sub_anova2}
\end{equation}
where $f_{\emptyset}={\rm const}$, condition (\ref{sparsitycond}) holds for all $1\leq k\leq s$,
and 
\begin{gather}\label{orthcon1}
	\int_0^1f_{u_k}(\tb_{u_k})\,dt_j=0\quad 
	\mbox{ for all}\;j\in u_k, \;u_k\in \Ukd,\; 1\leq k\leq s,
\end{gather}
will also be addressed.
If $d=d_\varepsilon\to \infty$ as $\varepsilon\to 0$ and
$s$ is either fixed or $s=s_\varepsilon\to \infty$, $s=o(d)$, as $\varepsilon\to 0$, then 
\begin{equation*} \label{sparsitycond_s}
	\sum_{k=1}^{s} \sum_{u_k \in \Ukd} \eta_{u_k}\leq  C \sum_{k=1}^s  {d\choose k}^{1-\beta} = C  {d\choose s}^{1-\beta}(1+o(1)), \quad \varepsilon\to 0, 
	\end{equation*}
	that is, at most $O \left( {d\choose s}^{1-\beta}\right )=o({d\choose s})$ orthogonal components $f_{u_k}$ of $f$ in (\ref{fun_sub_anova2}) are active
	and the remaining
	components are inactive. This implies that the function $f$  is \textit{sparse} and is composed of functions of a small number of variables, 
	with all these functions being orthogonal
to each other and to  a constant in $L_2^d$.
For use later on, we also define the sets
$\mathcal{H}_{s,d}^{\beta}$, $1\leq s\leq d$, as follows:
\begin{gather}\label{cH}
	{\cal H}^{\beta}_{s,d}=\left\{\etab=(\etab_1,\ldots,\etab_s): \etab_k\in { H}_{k,d}^{\beta},\,1\leq k\leq s\right \}.
\end{gather}
In Section \ref{Ext}, based on the results of Section \ref{MR},  we shall study
a more general problem of the recovery of the relevant (i.e. nonzero) components  of a collection of vectors $\etab=(\etab_1,\ldots,\etab_s)\in {\cal H}^{\beta}_{s,d}$.

\subsection{Regularity conditions}\label{sec:reg_cond}
Functional ANOVA expansions appear in many contexts (see, for example, \cite{Lin-2000} and \cite{Wahba-1995}).
This makes the recovery of the sparsity pattern of $f$ observed in the Gaussian white noise model 
a problem of great interest. 
In order to obtain a meaningful problem of sparse signal recovery in the model given by  (\ref{model1}) and (\ref{f})--(\ref{sparsitycond}),
we have to assume that the set of signals $f$ in model (\ref{model1}) is not ``too large".
In this article, we will be interested  in periodic Sobolev classes described by means of Fourier coefficients.
Such classes are quite common in the literature on nonparametric estimation, signal detection, and variable selection.

Using the notation introduced in Section \ref{SCond}, let $\theta_{\lb}(u_k)=(f_{u_k},\phi_{\lb})_{L_2^d}$ be the $\lb$th Fourier coefficient of $f_{u_k}$ for $\lb \in 	 \mathring{\mathbb{Z}}_{u_k}$, $u_k\in\Ukd$, $1\leq k\leq d$.
Then, for $u_k=\{j_1,\ldots,j_k\}\in \Ukd$, where $1 \leq j_1 < \ldots < j_k \leq d$, the $k$-variate component $f_{u_k}$ on the right-hand side of (\ref{f}) can be expressed as
\begin{equation*}
f_{u_k}(\tb_{u_k})=\sum_{\lb\in \mathring{\mathbb{Z}}_{u_k}}\theta_{\lb}(u_k) \phi_{\lb}(\tb_{u_k}),
\end{equation*}
and the entire function $f$ in decomposition (\ref{f}) takes the form
\begin{gather*}
f(\tb)= \sum_{u_k\in \Ukd}\eta_{u_k} \sum_{\lb\in \mathring{\mathbb{Z}}_{u_k}}\theta_{\lb}(u_k) \phi_{\lb}(\tb_{u_k}).
\end{gather*}
Note that only those Fourier coefficients of $f$ that correspond to the components $f_{u_k}$ in~(\ref{f}) are nonzero and that
$\|f_{u_k}\|^2_{2}=(f_{u_k},f_{u_k})_{L_2^d}=\sum_{\lb\in \mathring{\mathbb{Z}}_{u_k}} \theta_{\lb}^2 (u_k)$.

For $u_k\in\Ukd$ as above, 
we assume that $f_{u_k}$ belongs to the Sobolev class of $k$-variate functions with integer smoothness parameter $\sigma\geq 1$ for which
the semi-norm $\|\cdot\|_{\sigma,2}$ is defined by
\begin{gather}\label{semi-norm}
\|f_{u_k}\|_{\sigma,2}^2=\sum_{i_1=1}^k\ldots\sum_{i_{\sigma}=1}^k \left\|\frac{\partial^{\sigma} f_{u_k}  }{\partial t_{j_{i_1}}\ldots\partial t_{j_{i_{\sigma}}  }}   \right\|^2_2.
\end{gather}
Under the periodic constraint, we can define the semi-norm $\|\cdot\|_{\sigma,2}$ for the general case $\sigma>0$
in terms of the Fourier coefficients $\theta_{\lb}(u_k)$, $\lb \in \mathring{\mathbb{Z}}_{u_k}$.
For this, assume that $f_{u_k}(\tb_{u_k})$ admits 1-periodic $[\sigma]$-smooth extension in each argument to $\mathbb{R}^k$,
i.e., for all derivatives $f_{u_k}^{(n)}$ of integer order $0\leq n\leq [\sigma]$, where $f_{u_k}^{(0)}=f_{u_k}$, one has
\begin{gather*}
f_{u_k}^{(n)}(t_{j_1},\ldots,t_{j_{i-1}},0,t_{j_{i+1}},\ldots,t_{j_k})=f_{u_k}^{(n)}(t_{j_1},\ldots,t_{j_{i-1}},1,t_{j_{i+1}},\ldots,t_{j_k}),\quad 2\leq i\leq k-1,
\end{gather*}
with obvious extension for $i=1,k$. Then, the expression in (\ref{semi-norm}) corresponds to
\begin{gather}
\|f_{u_k }\|_{\sigma,2}^2=\sum_{\lb\in \mathring{\mathbb{Z}}_{u_k} } \theta_{\lb}^2(u_k) c_{\lb}^2,\quad \quad  c_{\lb}^2=\left(\sum_{j=1}^d (2\pi l_{j})^2  \right)^{\sigma}
=\left(\sum_{i=1}^k (2\pi l_{j_i})^2  \right)^{\sigma}.
\label{def:c_lb}
\end{gather}
Indeed, for $u_k= \{j_1,\ldots,j_k\}$, $1\leq j_1<\ldots<j_k\leq d$, and an integer number $\sigma\geq 1$, 
noting that $\phi^{\prime}_l(t)=-(2\pi l)\phi_{-l}(t)$ for $l\in \mathbb{Z}$, we can write
\begin{gather*}
	\frac{\partial^{\sigma} f_{u_k}(\tb_{u_k})  }{\partial t_{j_{i_1}}\ldots\partial t_{j_{i_{\sigma}}  }}  =
	\sum_{\lb\in \mathring{\mathbb{Z}}_{u_k}}\theta_{\lb}(u_k)\frac{\partial^{\sigma} \phi_{\lb}(\tb_{u_k})  }{\partial t_{j_{i_1}}\ldots\partial t_{j_{i_{\sigma}}  }}=
	\sum_{\lb\in \mathring{\mathbb{Z}}_{u_k}} (-2\pi)^{\sigma}l_{j_{i_1}}\ldots l_{j_{i_\sigma}}\theta_{\lb}(u_k)\phi_{\lb^*}(\tb_{u_k}) 
\end{gather*}
where $\lb^*=(l_1^*,\ldots, l_d^*)$ with $l_j^*=-l_j$ for $j\in \{j_{i_1},\ldots, j_{i_{\sigma}}\}$ and $l_j^*=l_j$ otherwise. 
From this and orthonormality of the basis $\{\phi_{\lb}\}_{\lb\in \mathbb{Z}^d}$, we have 
\begin{align*}
	\left\|\frac{\partial^{\sigma} f_{u_k}  }{\partial t_{j_{i_1}}\ldots\partial t_{j_{i_{\sigma}}  }}   \right\|^2_2 &=
	\int_{[0,1]^d}\left( \sum_{\lb\in \mathring{\mathbb{Z}}_{u_k}}(-2\pi)^{\sigma}l_{j_{i_1}}\ldots l_{j_{i_\sigma}}\theta_{\lb}(u_k)\phi_{\lb^*}(\tb_{u_k}) \right)^2 d\tb \\
	&= \sum_{\lb\in \mathring{\mathbb{Z}}_{u_k}}(2\pi)^{2\sigma} l^2_{j_{i_1}}\ldots l^2_{j_{i_\sigma}}\theta^2_{\lb}(u_k).
\end{align*}
Therefore, the right-hand side of (\ref{semi-norm}) coincides with the expression for $\|f_{u_k }\|_{\sigma,2}^2$ in (\ref{def:c_lb}), since
\begin{multline*}
	\sum_{i_1=1}^k\ldots\sum_{i_{\sigma}=1}^k \left\|\frac{\partial^{\sigma} f_{u_k}  }{\partial t_{j_{i_1}}\ldots\partial t_{j_{i_{\sigma}}  }}   \right\|^2_2=
	\sum_{\lb\in \mathring{\mathbb{Z}}_{u_k}}(2\pi)^{2\sigma}\left(\sum_{i_1=1}^k\ldots\sum_{i_\sigma=1}^k l^2_{j_{i_1}}\ldots l^2_{j_{i_\sigma}} \right)
	\theta_\lb^2(u_k)\\=
	\sum_{\lb\in \mathring{\mathbb{Z}}_{u_k}}(2\pi)^{2\sigma}\left(\sum_{i_1=1}^k l^2_{j_{i_1}}\ldots \sum_{i_\sigma=1}^k l^2_{j_{i_\sigma}}\right)\theta_\lb^2(u_k)=
	\sum_{\lb\in \mathring{\mathbb{Z}}_{u_k}}\left( \sum_{i=1}^k (2\pi l_{j_i})^2 \right)^{\sigma}\theta_\lb^2(u_k).
\end{multline*}

Finally, denote by ${\cal F}_{\cb_{u_k}}$ the Sobolev ball of radius 1 with coefficients ${\cb}_{u_k}=({c}_{\lb})_{\lb\in\mathring{\mathbb{Z}}_{u_k} }$, that is, cf.~(\ref{def:c_lb}),
\begin{gather*}
{\cal F}_{\cb_{u_k}}=\left \{  f_{u_k}(\tb_{u_k})= \sum_{\lb\in\mathring{\mathbb{Z}}_{u_k} } \theta_{\lb}(u_k)\phi_{\lb}(\tb_{u_k}), \; \tb_{u_k} \in [0,1]^k:
\sum_{\lb\in \mathring{\mathbb{Z}}_{u_k} } \theta_{\lb}^2(u_k) c_{\lb}^2 \leq 1 \right \},
\end{gather*}
and assume that every component $f_{u_k}$ of $f$ in (\ref{f}) belongs to this Sobolev ball, that is,
\begin{gather}\label{Fu}
f_{u_k} \in {\cal F}_{\cb_{u_k}},\quad u_k \in \Ukd.
\end{gather}
Thus, the model of our interest is specified by equations (\ref{model1}), (\ref{f})--(\ref{sparsitycond}), and (\ref{Fu}).
Clearly, the Sobolev balls ${\cal F}_{\cb_{u_k}}$ are isomorphic for all $u_k$
of cardinality $k$ ($1\leq k\leq d$).

\subsection{Problem statement}\label{PS}

The problem of recovering
the sparsity pattern of a multivariate signal observed in the Gaussian white noise
has been previously  studied in \cite{BSt-2017, CD-2012,ISt-2014, ST-2023, ST-2024}.
In our context, the problem is to identify the relevant components of a binary-valued vector $\etab_k\in H_{k,d}^{\beta}$  based on an ``observation'' $\mathbb{X}_\varepsilon$.
This problem will be named the \textit{variable selection problem}, and
an estimator $\boldsymbol{\hat{\eta}}_k=\boldsymbol{\hat{\eta}}_k(\mathbb{X}_\varepsilon)=(\hat{\eta}_{u_k})_{u_k\in \Ukd}$ of $\etab_k=(\eta_{u_k})_{u_k\in \Ukd}\in H_{k,d}^{\beta}$
taking values in $\{0,1\}^{{d\choose k}}$ will be referred to as a \textit{selector}.
In the literature on variable selection in high dimensions, it is common to distinguish between exact and almost full selectors.
A selector $\boldsymbol{\hat{\eta}}_k$ is called \textit{exact} if its maximum risk is algebraically small for all large enough $d$, and it is called
\textit{almost full} if its maximum risk is small compared to the number of relevant components of $\etab_k\in H_{k,d}^{\beta}$.
Such a division of selectors into two groups was proposed in  \cite{GJWY-2012}.
For the exact variable selection problem to be meaningful, the function components $f_{u_k}$ of $f$ in the model
given by (\ref{model1}), (\ref{f})--(\ref{sparsitycond}), and (\ref{Fu})
should be separated from zero.
If at least one of the components $f_{u_k}$ is  ``too small'' and exact selection is impossible, a procedure that provides selection almost fully is sought after.
Unlike exact selection, almost full selection can be achieved  under {milder} assumptions on the statistical model.
In this article, we are interested in establishing conditions for the possibility and impossibility of almost full selection
in the model at hand, and providing an almost full selector that works for all values of the sparsity index $\beta$.
Compared to the exact selection framework studied in  \cite{ST-2023} and \cite{ST-2024},
the construction of an \textit{adaptive} almost full selector that works for all values of $\beta$, which is generally unknown, is
a more challenging problem. The reason for this is that the (nonadaptive) 
almost full selector $\hat{\etab}_k(\beta)=(\hat{\eta}_{u_k}(\beta))_{u_k\in\Ukd}$, 
defined  below by relations (\ref{def:selector_for_beta_known})--(\ref{cond_epsilon}), acts as a hard thresholding procedure that
identifies each function component $f_{u_k}$, $u_k\in \Ukd$, as active if an appropriate statistic exceeds 
a certain threshold, with both the statistic and threshold depending on $\beta$. A natural step towards 
making the proposed selector adaptive would be to estimate
$\beta$ by its efficient estimator. However, this  does not easily lead (if at all) to the adaptive selector with required properties. For this reason,
in order to construct the desired  selection procedure, we employ Lepski's method of adaptation from \cite{Lepski1990}.
This is different from the exact variable selection framework, where adaptation of the exact selector to unknown $\beta$
can be achieved by utilizing Bonferroni's approach (for details, see Section  2 of \cite{ST-2023}). 
The latter approach does not work for the problem under study, since the dependence of 
the almost full selector  $\hat{\etab}_k(\beta)=(\hat{\eta}_{u_k}(\beta))_{u_k\in\Ukd}$  in  (\ref{def:selector_for_beta_known}) 
on the parameter $\beta$ is more complicated than that of the exact selector in Section 2 of \cite{ST-2023}:
the threshold of the exact selector 
is independent of $\beta$, 
whereas the threshold for $\hat{\etab}_k(\beta)=(\hat{\eta}_{u_k}(\beta))_{u_k\in\Ukd}$  in  (\ref{def:selector_for_beta_known}), 
which is set at a lower level, does depend on $\beta$.

To quantify
the performance  of a selector, we shall study the \textit{Hamming risk} of $\boldsymbol{\hat{\eta}}_k$ as an estimator of $\etab_k\in  H_{k,d}^{\beta}$
defined by
\begin{equation*} 
{\rm E}_{\varepsilon,f}|\boldsymbol{\hat{\eta}}_k-\etab_k| :={\rm E}_{\varepsilon,f}\left(\sum_{u_k\in\Ukd}|\hat{\eta}_{u_k}-\eta_{u_k}|\right),
\end{equation*}
where $\E{\varepsilon,f}$ is the expectation with respect to the probability measure $\operatorname{P}_{\varepsilon,f}$, which was introduced in the Introduction. The Hamming risk corresponds to the expected
number of components for which the selector $\boldsymbol{\hat{\eta}}_k$ is not in agreement with $\etab_k$.
We define an \textit{almost full selector} $\boldsymbol{\hat{\eta}}_k$ to be a selector whose maximum Hamming risk is small compared to
the number of active components (see relation (\ref{r1}) below for the precise definition).
The problem of identifying almost fully
the active components of $f$ when the sparsity index $\beta$ is \textit{known}
can be settled without much difficulty, whereas the solution to the same problem when $\beta$ is \textit{unknown}
requires further subtle arguments.

In this article, we  establish conditions when almost full recovery
of the sparsity pattern of $f$ in the model given by (\ref{model1}), (\ref{f})--(\ref{sparsitycond}), and (\ref{Fu}) 
is possible, as well as when it is impossible, and propose
a~procedure that achieves this type of recovery. In other words, we
construct an estimator $\boldsymbol{\hat{\eta}}_k=\boldsymbol{\hat{\eta}}_k(\mathbb{X}_\varepsilon)=(\hat{\eta}_{u_k})_{u_k\in \Ukd}$ of $\etab_k=(\eta_{u_k})_{u_k\in \Ukd}\in H_{k,d}^{\beta}$ that would tell us which components  $f_{u_k}$ of $f$ in (\ref{f}) are active. 
Then, we extend the obtained results to the sparse regression model with function $f$ decomposed as in (\ref{fun_sub_anova2}).
Extensions will be provided for both cases, when $s$ is fixed and when $s=s_\varepsilon\to \infty$, $s=o(d)$ as $\varepsilon\to 0.$

Identifying the active components of $f$ in the model given by (\ref{model1}), (\ref{f})--(\ref{sparsitycond}), and (\ref{Fu}) 
is feasible  when,
in addition to the regularity constraints in (\ref{Fu}), the components $f_{u_k}$ of $f$ are not ``too small'', i.e., separated from zero in an appropriate way.
Therefore, following \cite{ST-2023}, for a given $u_k\in \Ukd$, $1\leq k\leq d$, and $r >0$, we define the set
$$\mathring{\cal F}_{\cb_{u_k}}(r)=\{ f_{u_k} \in \mathcal{F}_{\cb_{u_k}}: \| f_{u_k}  \|_2\geq r\}, $$
and consider testing
\begin{gather}
\mathbb{H}_{0,u_k} :  f_{u_k}=0 \quad\mbox{vs.}\quad \mathbb{H}^{\varepsilon}_{1,u_k}:  f_{u_k}\in \mathring{\cal F}_{\cb_{u_k}}(r_{\varepsilon,k}), 
\label{hypotheses0}
\end{gather}
for some positive family $r_{\varepsilon,k}\to 0$ as $\varepsilon\to 0$.
The hypothesis testing problem (\ref{hypotheses0}), known in the literature as the signal detection problem,
has been studied in \cite{IS-2005}. In the present context of sparse signal recovery,
this is an auxiliary problem that enables us to obtain the conditions when almost full selection of active components of $f$ is possible and
when this type of selection is impossible. Additionally, we shall use the
asymptotically minimax test statistics from the above signal detection problem (see Theorem 2 of \cite{IS-2005}) to design an almost full selector.

For a positive family $r_{\varepsilon,k}$ as above, we introduce the class of sparse multivariate functions of our interest as follows:
\begin{multline*}
{\cal F}_{k,d}^{\beta,\sigma} (r_{\varepsilon,k})=\Bigg\{f: f(\tb)=\sum\nolimits_{u_k\in \Ukd} \eta_{u_k} f_{u_k}(\tb_{u_k}), f_{u_k}\;\mbox{satisfies}\;(\ref{orthcon}),\\
f_{u_k}\in \mathring{\cal F}_{\cb_{u_k}}(r_{\varepsilon,k}), u_k\in \Ukd, \etab_k=(\eta_{u_k})_{u_k\in \Ukd}\in H_{k,d}^{\beta}\Bigg\}.
\end{multline*}
The dependence of ${\cal F}_{k,d}^{\beta,\sigma}$ on the smoothness parameter $\sigma$ is hidden in the coefficients $\cb_{u_k}=(c_{\lb})_{\lb\in\mathring{\mathbb{Z}}_{u_k}}$
defining the set $\mathring{\cal F}_{\cb_{u_k}}(r_{\varepsilon,k})$.
In this work, we are interested in selecting the active components of $f\in {\cal F}_{k,d}^{\beta,\sigma}$ almost fully. Therefore, we first establish the sharp almost full selection boundary
that allows us to verify whether the active components of $f$ can be selected almost fully,
and then construct a selector $\boldsymbol{\hat{\eta}}_k=\boldsymbol{\hat{\eta}}_k(\mathbb{X}_\varepsilon)\in\{0,1\}^{d\choose k}$ (for known and unknown $\beta$) attaining this boundary with the following property:
for all $\beta\in(0,1)$ and  $\sigma>0$,
\begin{gather}\label{r1}
\limsup_{\varepsilon\to 0}\sup_{\etab_k\in H_{k,d}^{\beta}} \sup_{f\in  {\cal F}_{k,d}^{\beta,\sigma}} {d \choose k}^{\beta-1} {\rm E}_{\varepsilon,f}|\boldsymbol{\hat{\eta}}_k-\etab_k| =0.
\end{gather}
Relation (\ref{r1}) says that the maximum normalized Hamming risk of $\boldsymbol{\hat{\eta}}_k$ is small relative to the total number of active components of $f$ in
the model given by (\ref{model1}), (\ref{f})--(\ref{sparsitycond}), and (\ref{Fu}), and thus $\boldsymbol{\hat{\eta}}_k$ recovers $\etab_k$ almost fully.
Additionally, we show that for all
those values of $r_{\varepsilon,k}$ that fall below the \textit{almost full selection boundary},
one has
\begin{gather}\label{r2}
\liminf_{\varepsilon\to 0}\inf_{\boldsymbol{\tilde{\eta}}_k} \sup_{\etab_k\in H_{k,d}^{\beta}}  \sup_{f\in  {\cal F}_{k,d}^{\beta,\sigma}} {d \choose k}^{\beta-1} {\rm E}_{\varepsilon,f}|\boldsymbol{\tilde{\eta}}_k-\etab_k|>0,
\end{gather}
where the infimum is taken over all selectors $\boldsymbol{\tilde{\eta}}_k$ of $\etab_k\in H_{k,d}^{\beta}$ in the model at hand,
that is, almost full recovery of the sparsity pattern of  $f\in {\cal F}_{k,d}^{\beta,\sigma}$ in the model given by 
(\ref{model1}), (\ref{f})--(\ref{sparsitycond}), and (\ref{Fu}) 
is impossible.
A similar problem for the case of exact selection has been addressed and solved in~\cite{ST-2023}.

In the problem of recovering the sparsity patter of $f$, it is more convenient to deal with the space of Fourier coefficients of the functions of interest rather than the functions themselves. 
In the sequence space of Fourier coefficients, the unit ball ${\cal F}_{\cb_{u_k}}$ corresponds to the ellipsoid  $ {\Theta}_{\cb_{u_k}}$ 
defined by
\begin{eqnarray*}
	{\Theta}_{\cb_{u_k}}&=&\bigg\{\thetab_{u_k}=(\theta_{\lb})_{\lb\in  \mathring{\mathbb{Z}}_{u_k}}\in l_2({\mathbb{Z}}^d):
	\sum_{\lb \in  \mathring{\mathbb{Z}}_{u_k} }c_{\lb}^2\theta_{\lb}^2(u_k)\leq 1\bigg\},
	\end{eqnarray*}
and the set $\mathring{{\cal F}}_{\cb_{u_k}}(r_{\varepsilon,k})$
corresponds to the ellipsoid with a small neighbourhood of the point $\thetab_{u_k}={\bf 0}$ removed:
\begin{eqnarray*}
		\mathring{\Theta}_{\cb_{u_k}}(r_{\varepsilon,k})&=&\bigg\{\thetab_{u_k}=(\theta_{\lb})_{\lb\in  \mathring{\mathbb{Z}}_{u_k}}\in {\Theta}_{\cb_{u_k}}:
	\sum_{\lb \in  \mathring{\mathbb{Z}}_{u_k} }\theta_{\lb}^2(u_k)\geq r^2_{\varepsilon,k}\bigg\}.
\end{eqnarray*}
Then, in terms of Fourier coefficients, the  testing problem in (\ref{hypotheses0})
takes the form
\begin{gather}\label{test1}
	{\bf H}_{0,u_k}:\thetab_{u_k}={\bf 0}\quad \mbox{vs.}\quad {\bf H}_{1,u_k}^\varepsilon: \thetab_{u_k}\in \mathring{\Theta}_{\cb_{u_k}}(r_{\varepsilon,k}).
\end{gather}
Note that $\mathring{\Theta}_{\cb_{u_k}}(r_{\varepsilon,k})=\emptyset$ for $r_{\varepsilon,k}>1/c_{\varepsilon,0},$ where, by
using (\ref{def:c_lb}), $c_{\varepsilon,0}:=\inf_{\lb\in \mathring{\mathbb{Z}}_u}c_{\lb}=(2\pi)^{\sigma}k^{\sigma/2}$.
Therefore, in what follows, we will be interested only in the case when $r_{\varepsilon,k}\in (0,(2\pi)^{-\sigma}k^{-\sigma/2}).$

The initial model given by (\ref{model1}), (\ref{f})--(\ref{sparsitycond}), and (\ref{Fu}) 
can be equivalently represented  in terms of the Fourier coefficients of the function components $f_{u_k}$  as (see, for example, Section 1.2 of  \cite{ST-2023})
\begin{gather}\label{model2}
X_{\lb} =\eta_{u_k}\theta_{\lb}(u_k)+\varepsilon\xi_{\lb},\quad \lb \in  \mathring{\mathbb{Z}}_{u_k},\quad u_k\in \Ukd,
\end{gather}
where $X_\lb= \mathbb{X}_\varepsilon(\phi_\lb)$ is the $\lb$th empirical Fourier coefficients, $\etab_k=(\eta_{u_k})_{u_k\in \Ukd}\in H_{k,d}^{\beta}$, $\xi_{\lb}=\mathbb{W}(\phi_{\lb})$
are iid standard normal random variables for ${\lb \in  \mathring{\mathbb{Z}}_{u_k}}$ and $u_k\in \Ukd$,
and  $\thetab_{u_k}=(\theta_{\lb}(u_k),{\lb\in \mathring{\mathbb{Z}}_{u_k}})$ consists of the
Fourier coefficients
$\theta_\lb(u_k)=(f_{u_k},\phi_\lb)_{L_2^d}$ of $f_{u_k}$ and belongs to the ellipsoid $	{\Theta}_{\cb_{u_k}}$.
Model (\ref{model2}) is known in literature as the Gaussian sequence space model.
From a technical point of view, it is more convenient to deal with ellipsoids in sequence spaces rather than Sobolev balls in function spaces.
Although the set of $\theta_{\lb}(u_k)$s in this model involves the orthonormal basis $\{\phi_{\lb}\}_{\lb\in \mathbb{Z}^d}$ of $L_2^d$, 
the results on minimax risks do not depend on the choice of this  basis because 
the random variables $X_{\lb}$ in model (\ref{model2}), which generate a sufficient $\sigma$-algebra for $f\in {\cal F}_{k,d}^{\beta,\sigma}$, are independent
normal $N(\eta_{u_k}\theta_{\lb}(u_k),\varepsilon^2)$ random variables. Thus the distribution of
$\{\Xb_{\!{u_k}}\}_{u_k\in \Ukd}$, where  $\Xb_{\!{u_k}}=(X_{\lb})_{\lb\in \mathring{\mathbb{Z}}_{u_k}}$,
depends on the Fourier coefficients $\theta_{\lb}(u_k)$ of $f_{u_k}$s with respect to the basis $\{\phi_{\lb}\}_{\lb\in \mathring{\mathbb{Z}}_{u_k}}$ 
of the set ${\cal L}_{u_k}^d$, 
but not on the choice of the basis itself.

Now, we define the set
\begin{equation*} 
\mathring{\Theta}^\sigma_{k,d}(r_{\varepsilon,k})=\left\{\thetab_k: \thetab_k=(\thetab_{u_k})_{u_k\in \Ukd},\;\mbox{where}\;
\thetab_{u_k}=(\theta_{\lb}(u_k))_{\lb\in \mathring{\mathbb{Z}}_{u_k}} \in \mathring{\Theta}_{\cb_{u_k}}(r_{\varepsilon,k})  \right\}.  
\end{equation*}
Then, in terms of model (\ref{model2}),
relations (\ref{r1}) and (\ref{r2}) take the form
\begin{gather}\label{r11}
\limsup_{\varepsilon\to 0} \sup_{\etab_k\in  H_{k,d}^{\beta}}\sup_{\thetab_k\in  \mathring{\Theta}^\sigma_{k,d}(r_{\varepsilon,k})} {d \choose k}^{\beta-1}{\rm E}_{\thetab_k,\etab_k}|\boldsymbol{\hat{\eta}}_k-\etab_k| =0,
\end{gather}
and
\begin{gather}\label{r22}
\liminf_{\varepsilon\to 0}\inf_{\tilde{\etab}_k}\sup_{\etab_k\in  H_{k,d}^{\beta}} \sup_{\thetab_k\in \mathring{\Theta}^\sigma_{k,d}(r_{\varepsilon,k})} {d \choose k}^{\beta-1} {\rm E}_{\thetab_k,\etab_k}|\boldsymbol{\tilde{\eta}}_k-\etab_k|>0,
\end{gather}
where  $\boldsymbol{\hat{\eta}}_k$ and $\boldsymbol{\tilde{\eta}}_k$ are estimators of $\etab_k$ based on  $\{\Xb_{\!{u_k}}\}_{u_k\in \Ukd}$, $\Xb_{\!{u_k}}=(X_{\lb})_{\lb\in \mathring{\mathbb{Z}}_{u_k}}$,
and ${\rm E}_{\thetab_k,\etab_k}$ is the expectation with respect to the probability distribution $\operatorname{P}_{\thetab_k,\etab_k}$ of $\{\Xb_{\!{u_k}}\}_{u_k\in \Ukd}$.
An \textit{almost full selector} $\boldsymbol{\hat{\eta}}_k\in\{0,1\}^{d\choose k}$ of ${\etab}_k$ in model (\ref{model2}) 
will be defined as a selector satisfying (\ref{r11}), cf. relation (17) in \cite{ST-2023}.
The limiting relations (\ref{r11}) and (\ref{r22}) will be referred to as the \textit{upper bound} on the normalized  maximum Hamming risk of $\boldsymbol{\hat{\eta}}_k$
and the \textit{lower bound} on the minimax normalized Hamming risk, respectively.
When the  upper bound in (\ref{r11}) holds true, the maximum Hamming risk of the selector $\boldsymbol{\hat{\eta}}_k$ is small
compared to the number of relevant components of $\etab_k\in  H_{k,d}^{\beta}$, which is at most $C {d\choose k}^{1-\beta}$,
and thus  $\boldsymbol{\hat{\eta}}_k$ achieves almost full selection.
Also, when the lower bound in (\ref{r22}) holds true, the minimax Hamming risk is at least $O\left({d\choose k}^{1-\beta}\right)$,
and thus any variable selection procedure fails completely.

\subsection{Notation}\label{N}
The majority of limits in the article are taken as $\varepsilon\to 0$.
We use the notation $a_\varepsilon \asymp b_\varepsilon$ when
$0 < \liminf_{\varepsilon \to 0} {a_\varepsilon}/{b_\varepsilon} \leq \limsup_{\varepsilon\to 0} {a_\varepsilon}/{b_\varepsilon} < \infty$
and the notation $a_\varepsilon \sim b_\varepsilon$ when
$\lim_{\varepsilon \to 0} {a_\varepsilon}/{b_\varepsilon} =1.$
For positive families $a_\varepsilon$ and $b_\varepsilon$, we shall occasionally write $a_\varepsilon\ll b_\varepsilon$ when $a_\varepsilon=o(b_\varepsilon)$ as $\varepsilon\to 0$.
For a real number $x$,  $x_+ = \max (0,x)$, and the symbol $\lceil x\rceil$  denotes the smallest integer strictly larger than $x$.
We use the symbol $\log x$ for the natural (base $e$) logarithm of the number~$x$.

\section{Construction of an almost full selector}\label{Selector}
For $u_k\in \Ukd$ and   $r_{\varepsilon,k}\in (0,(2\pi)^{-\sigma}k^{-\sigma/2})$, return to the problem of testing  $ {\bf H}_{0,u_k}$ against
${\bf H}^{\varepsilon}_{1,u_k}$ specified by (\ref{test1}) 
and consider the quantity
\begin{equation}
a^2_{\varepsilon,u_k} (r_{\varepsilon,k})= \frac{1}{2 \varepsilon^4} \inf_{\thetab_{u_k} \in \mathring{\Theta}_{\cb_{u_k}}(r_{\varepsilon,k})} \sum_{\lb\in \mathring{\mathbb{Z}}_{u_k}} \left(\theta_\lb(u_k)\right) ^4,
\label{def:a2}
\end{equation}
which is known to control the minimax total error probability and determines a cut-off point of the
asymptotically minimax test procedure in the problem of testing  $ {\bf H}_{0,u_k}$ against
${\bf H}^{\varepsilon}_{1,u_k}$ (for details, see Theorem 2 of \cite{IS-2005}).
Additionally, the function $a^2_{\varepsilon,u_k} (r_{\varepsilon,k})$ turns out to play a key role in establishing
conditions under which the active components $f_{u_k}$ of $f$ in the model given by (\ref{model1}), (\ref{f})--(\ref{sparsitycond}), and (\ref{Fu}) can be selected \textit{almost fully}.

Observe that $a_{\varepsilon,u_k} (r_{\varepsilon,k})$ is a nondecreasing function of its argument that possesses a kind of ``continuity'' property. Namely, for any $\gamma >0$, there exist $\varepsilon^*>0$ and $\delta^* >0$ such that (see Section 3.2 of~\cite{I1993a}) 
\begin{equation} 
a_{\varepsilon,u_k} (r_{\varepsilon,k}) \leq a_{\varepsilon,u_k} ((1+\delta)r_{\varepsilon,k}) \leq (1+\gamma) a_{\varepsilon,u_k} (r_{\varepsilon,k}), \quad \forall\, \varepsilon \in (0,\varepsilon^*), \forall\, \delta \in (0,\delta^*).
\label{def:continuity_of_a}
\end{equation} 
These and other general facts of the minimax hypothesis testing theory can be found in
a series of review articles \cite{I1993a}--\cite{I1993c} and monograph \cite{NGF-2003}.
Suppressing for brevity the dependence on $u_k$, denote the minimizing sequence in (\ref{def:a2}) by $(\theta^*_\lb (r_{\varepsilon,k}))_{\lb \in \mathring{\mathbb{Z}}_{u_k}}$, that is,
\begin{equation}
a^2_{\varepsilon,u_k} (r_{\varepsilon,k})= \frac{1}{2 \varepsilon^4} \sum_{\lb\in \mathring{\mathbb{Z}}_{u_k}} \left( \theta^{*}_\lb (r_{\varepsilon,k})\right)^4,
\label{def:a2_with_theta_star}
\end{equation}
and let $r^*_{\varepsilon,k}=r^*_{\varepsilon,k}(\beta)>0$ be determined by, cf. condition (\ref{cond:inf}) in Theorem \ref{th:sim_UB_fixedK_knownBeta} below,
\begin{equation}
a_{\varepsilon,u_k} (r^*_{\varepsilon,k}) = \sqrt{2 \beta \log {d\choose k}}.
\label{def:r_star}
\end{equation}
Assume for a while that the sparsity index $\beta$ is known. For the purpose of constructing an almost full
selector $\boldsymbol{\hat{\eta}}_k$ satisfying (\ref{r11}), we consider weighted $\chi^2$-type statistics
\begin{equation} \label{def:S_u}
S_{u_k}(\beta) = \sum_{\lb \in \mathring{\mathbb{Z}}_{u_k}} \omega_\lb (r^*_{\varepsilon,k}) \left( \left( {X_\lb}/{\varepsilon}  \right)^2 -1 \right), \quad u_k \in \Ukd,
\end{equation}
where
\begin{equation*}
\omega_\lb (r_{\varepsilon,k}) = \frac{1}{2 \varepsilon^2} \frac{\left(\theta^*_\lb (r_{\varepsilon,k})\right)^2}{a_{\varepsilon,u_k}(r_{\varepsilon,k})}, \quad \lb \in \mathring{\mathbb{Z}}_{u_k}.
\label{def:omega}
\end{equation*}
Due to (\ref{def:a2_with_theta_star}),
\begin{equation*}
\sum_{\lb \in \mathring{\mathbb{Z}}_{u_k}} \omega_\lb^2 (r_{\varepsilon,k}) ={1}/{2} \quad \mbox{for all}\;\; r_{\varepsilon,k} > 0.
\end{equation*}
It is also known  that for all those
$\lb \in \mathring{\mathbb{Z}}_{u_k}$ for which the weights  $\omega_\lb (r_{\varepsilon,k})$
do not vanish, we have (see, for example, relation (45) of \cite{ST-2023})
\begin{gather}\label{WL}
	\omega_\lb (r_{\varepsilon,k})\asymp r_{\varepsilon,k}^{k/(2\sigma)},\quad \varepsilon\to 0.
\end{gather}

For every $u_k\in \Ukd$, the statistic $S_{u_k}(\beta)$ is a test statistic of the asymptotically minimax test procedure in the problem of testing  $ {\bf H}_{0,u_k}$ against ${\bf H}^{\varepsilon}_{1,u_k}$ (for details, see Theorem 2 of \cite{IS-2005}). Note that  $S_{u_k}(\beta)$ depends on the sparsity index $\beta$ through the weights $\omega_\lb (r^*_{\varepsilon,k}),$ $\lb\in \mathring{\mathbb{Z}}_{u_k}$.

It is known that for any fixed $k$ ($1\leq k\leq d$), as $\varepsilon\to 0$ (see relation (27) in Section 2 of \cite{ST-2023} for details) 
\begin{gather}\label{thetal}
[\theta^*_\lb (r_{\varepsilon,k})]^2\sim
\frac{r_{\varepsilon,k}^{2+{k}/{\sigma}  } 2^k \pi^{{k}/{2}}(k+2\sigma)\Gamma\left(1+{k}/{2}\right) }{2\sigma\left(1+4\sigma/k  \right)^{{k}/{(2\sigma)}}}
\left(1-\left(\sum\nolimits_{j=1}^d (2\pi l_j)^2  \right)^{\sigma}\frac{r^2_{\varepsilon,k} }{(1+4\sigma/k)}  \right)_+,
\end{gather}	
and hence every statistic $S_{u_k}(\beta)$ 
consists of $O((r_{\varepsilon,k}^*)^{-{k}/{\sigma}})$ nonzero terms,
since as $\varepsilon\to 0$
\begin{align*}
\#\{\lb \in \mathring{\mathbb{Z}}_{u_k}: \theta^*_\lb (r_{\varepsilon,k}^*)\neq 0\}&=
\#\left\{\lb \in \mathring{\mathbb{Z}}_{u_k}:\left(\sum_{j=1}^d l_j^2  \right)^{1/2}< \frac{(1+4\sigma/k)^{1/(2\sigma)}}{2\pi (r_{\varepsilon,k}^*)^{1/\sigma}}\right\} \\
&=O\left( (r_{\varepsilon,k}^*)^{-{k}/{\sigma}}\right).
\end{align*}
It is also known that for any fixed $k$, the sharp asymptotics of $a_{\varepsilon,u_k}(r_{\varepsilon,k})$, $u_k\in \Ukd$, as $\varepsilon\to 0$  are given by (see Theorem 4 of \cite{IS-2005})
\begin{gather}\label{aek}
a_{\varepsilon,u_k}(r_{\varepsilon,k})\sim C(\sigma,k)r_{\varepsilon,k}^{2+{k}/{(2\sigma)}}\varepsilon^{-2},\quad \quad 
C^2(\sigma,k)=\frac{\pi^k (1+{2\sigma}/{k})\Gamma(1+k/2)}{(1+4\sigma/k)^{1+{k}/{(2\sigma)}}\Gamma^k(3/2)}.
\end{gather}

In order to estimate a vector $\etab_k=(\eta_{u_k})_{u_k\in \Ukd}\in H_{k,d}^{\beta}$ in case of known $\beta$, we can use a selector $\boldsymbol{\hat{\eta}}_k(\beta)=(\hat{\eta}_{u_k}(\beta))_{u_k\in \Ukd}$, depending on the data $\{\Xb_{\!u_k}\}_{u_k\in \Ukd}$, where $\Xb_{\!u_k}=(X_{\lb})_{\lb\in \mathring{\mathbb{Z}}_{u_k}}$,
through the statistics $S_{u_k}(\beta)$, ${u_k\in \Ukd}$,
of the following form:
\begin{equation}
\boldsymbol{\hat{\eta}}_k(\beta)=(\hat{\eta}_{u_k}(\beta))_{u_k\in \Ukd}, \quad \quad \hat{\eta}_{u_k} (\beta)= \ind{ S_{u_k}(\beta) > \sqrt{(2\beta+\epsilon) \log {d\choose k} }},\quad u_k\in \Ukd,
\label{def:selector_for_beta_known}
\end{equation}
where $\epsilon=\epsilon_{\varepsilon,k,d} > 0$ is such that
\begin{gather}\label{cond_epsilon}
	\epsilon \to 0,\quad  \epsilon \log {d \choose k} \to \infty,\quad \epsilon\log {d\choose k}=o(\log \varepsilon^{-1}),\quad
	\mbox{as}\;\varepsilon\to 0.
\end{gather}
In other words, $\hat{\eta}_{u_k}(\beta)$ identifies the component $\eta_{u_k}$ as relevant if the value of $S_{u_k}(\beta)$
exceeds the threshold $\sqrt{(2\beta+\epsilon) \log {d\choose k}}$.
However, if $\beta$ is unknown, the selector $\boldsymbol{\hat{\eta}}_k(\beta)$ is not applicable anymore.

To construct a selector adapted to unknown $\beta$, we shall act similar to Lepski's method of adaptive estimation
(see \cite{Lepski1990}). To this end, 
we assume that $\beta\in [b,B]$ for some $0<b<B<1$, which is the price that is paid for \textit{adaptive} almost full recovery of the 
sparsity pattern, and consider a grid of equidistant points on $[b,1)$ defined by
\begin{gather} \label{def:grid_for_beta}
\beta_{k,1}=b,\quad\beta_{k,m} =\beta_{k,1}+ (m-1) \rho_k=\beta_{k,m-1}+\rho_k,\quad m=2,\ldots,M_k,
\end{gather}
where $M_k=\lceil (B-b)/\rho_k \rceil +1$ for some $\rho_k=\rho_{k,\varepsilon}>0$ such that
$ \rho_k\to 0,$ $\rho_k\log {d\choose k}\to 0$,  $\rho_k{d\choose k}\to \infty$, as $\varepsilon\to 0$,
or, in terms of the number of nodes $M_k$,
\begin{equation} \label{cond:Mk}
M_k \to \infty, \quad {\log {d \choose k}}/{M_k}\to 0,\quad {{d \choose k}}/M_k\to \infty, \quad \text{as } \varepsilon \to 0.
\end{equation}
The second relation in (\ref{cond:Mk}) implies that for all small enough~$\varepsilon$
\begin{equation} \label{eq:d_choose_k_to_rho_const}
{d \choose k}^{\rho_k} \leq {\rm const}.
\end{equation}
By definition,  $b=\beta_{k,1}< \ldots< \beta_{k,M_k}\in(B,B+\rho_k] \subset (0,1)$ and hence, for all small enough $\varepsilon$,
the grid points $\beta_{k,m}$ are all separated from 0 and 1.
Note also that for any $\beta\in [b,B]$ and all $1\leq k\leq d$ there exists an index $m_0=m_{0,k}$ ($1\leq m_0\leq M_k-1$) 
such that $\beta_{k,m_0}\leq \beta<\beta_{k,m_0+1}$.

Next, for $m=1,\ldots,M_k$, let $r^*_{\varepsilon,k,m} >0$ be determined by, cf. (\ref{def:r_star}),
\begin{equation}
a_{\varepsilon,u_k} (r^*_{\varepsilon,k,m}) =  \sqrt{2 \beta_{k,m} \log {d\choose k} }.
\label{def:r_star_ekm}
\end{equation}
It is known that $r^*_{\varepsilon,k}$ and $r^*_{\varepsilon,k,m}$ satisfy (see relation (46) of \cite{ST-2023})
\begin{gather}\label{rr}
r_{\varepsilon,k}^*\asymp r_{\varepsilon,k,m}^*\asymp\left(\varepsilon\log^{1/4} {d\choose k} \right)^{{4\sigma}/{(4\sigma+k)}},
\quad m=1,\ldots,M_k.
\end{gather}
For every $u_k \in \Ukd$ consider the statistics
\begin{equation}
S_{u_k}(\beta_{k,m}) = \sum_{\lb \in \mathring{\mathbb{Z}}_{u_k}} \omega_\lb (r^*_{\varepsilon,k,m}) \left( \left( {X_\lb}/{\varepsilon}  \right)^2 -1 \right), \quad m=1,\ldots,M_k,
\label{def:S_um}
\end{equation}
where, in view of (\ref{WL}) and (\ref{rr}), for all $m=1,\ldots,M_k$ 
\begin{equation} \label{omegamax_asymp}
	\max_{\lb \in \mathring{\mathbb{Z}}_{u_k}} \omega_\lb(r^*_{\varepsilon,k,m}) \asymp \left( \varepsilon \log^{1/4} {d\choose k}  \right)^{{2k}/{(4 \sigma +k)}}, \quad \varepsilon \to 0.
\end{equation}

Employing the statistics $S_{u_k}(\beta_{k,m})$ defined in (\ref{def:S_um}), we now define an adaptive selector $\boldsymbol{\hat{\eta}}_k$ by, cf.~(\ref{def:selector_for_beta_known}),
\begin{equation}
\boldsymbol{\hat{\eta}}_k=(\hat{\eta}_{u_k})_{u_k \in \Ukd }, \quad \quad \hat{\eta}_{u_k} =
\ind{ S_{u_k}(\beta_{k,\hat{m}_k}) > \sqrt{(2 {\beta}_{k,\hat{m}_k}+\epsilon) \log {d\choose k}} }, \quad u_k \in \Ukd,
\label{def:selector_for_beta_unknown}
\end{equation}
where $\epsilon=\epsilon_{\varepsilon,k,d}>0$ satisfies (\ref{cond_epsilon})
and the random index $\hat{m}_k\in\{1,\ldots,M_k\}$
is chosen by Lepski's method (see Section 2 of \cite{Lepski1990}) as follows:
\begin{equation} \label{Lepski_m}
\hat{m}_k = \max \left\{ 1 \leq m \leq M_k: |\boldsymbol{\hat{\eta}}_k (\beta_{k,m}) - \boldsymbol{\hat{\eta}}_k(\beta_{k,j})| \leq v_j \text{ for all $j < m$} \right\},
\end{equation}
and $\hat{m}_k = 1$ if the set above is empty. Here, the selector $\boldsymbol{\hat{\eta}}_k(\beta_{k,m}) =(\hat{\eta}_{u_k}(\beta_{k,m}))_{u_k\in \Ukd}$
consists of the components $\hat{\eta}_{u_k}(\beta_{k,m})=\ind{S_{u_k}(\beta_{k,m})>\sqrt{(2\beta_{k,m}+\epsilon)\log {d\choose k} }  }$
for  $m=1,\ldots,M_k$, and $ |\boldsymbol{\hat{\eta}}_k (\beta_{k,m}) - \boldsymbol{\hat{\eta}}_k(\beta_{k,j})|  = \sum_{u_k\in\Ukd}|\hat{\eta}_{u_k}(\beta_{k,m})-
\hat{\eta}_{u_k}(\beta_{k,j})|$ for $m,j=1,\ldots,M_k$. The quantities $v_j$ are taken to be
\begin{gather}\label{vj}
v_j = v_{j,k,d} = {d \choose k}^{1-\beta_{k,j}}/\tau_{k,d}, \quad j=1,\ldots,M_k,
\end{gather}
where, for every fixed $k$, $\tau_{k,d}=\tau_{\varepsilon,k,d}>0$  is such that 
\begin{gather}
	\tau_{k,d}  \to \infty, \quad  \tau_{k,d}=o\left({d \choose k}^{\epsilon/2}   \right), \quad \text{as } \varepsilon \to 0, \label{def:tau}
\end{gather}
When dealing with a more general decomposition of $f$ as in (\ref{fun_sub_anova2}), in which
$s=s_\varepsilon\to \infty$, $s=o(d)$ as $\varepsilon\to 0,$ instead of (\ref{def:tau}), we will require that 
\begin{gather}
	\tau_{k,d}  \to \infty, \quad  {d \choose k}^{\epsilon/8}\ll \tau_{k,d}\ll {d \choose k}^{\epsilon/2},  \quad \text{as } \varepsilon \to 0. \label{def:tau_s}
\end{gather}

The idea behind Lepski's procedure is as follows.
Each estimator $\boldsymbol{\hat{\eta}}_k(\beta_{k,j})$ is associated with an interval of a size proportional to its standard deviation. When, in the process of iteration, for 
some index $m$ ($m=1,\ldots,M_k$) an estimator $\boldsymbol{\hat{\eta}}_k(\beta_{k,m})$ is found to be far enough from the previous ones, this means that 
we have detected a bias change, indicating that we are probably close to the true value of~$\beta$.

Algorithmically, Lepski's procedure for choosing $\hat{m}_k$ in (\ref{Lepski_m}) works as follows.
We start by setting $\hat{m}_k=1$ and attempt to increase the value of $\hat{m}_k$
from 1 to 2. If $ |\boldsymbol{\hat{\eta}}_k (\beta_{k,2}) - \boldsymbol{\hat{\eta}}_k(\beta_{k,1})|\leq v_1$, we set
$\hat{m}_k=2$; otherwise, we keep $\hat{m}_k$ equal to 1. In case $\hat{m}_k$ is increased to 2,
we continue the process attempting to increase it further.
If  $ |\boldsymbol{\hat{\eta}}_k (\beta_{k,3}) - \boldsymbol{\hat{\eta}}_k(\beta_{k,1})|\leq v_1$
and  $ |\boldsymbol{\hat{\eta}}_k (\beta_{k,3}) - \boldsymbol{\hat{\eta}}_k(\beta_{k,2})|\leq v_2$,
we set $\hat{m}_k=3$; otherwise, we keep $\hat{m}_k=2$; and so on. Note that by construction
$v_1 > v_2 > \ldots > v_{M_k}$.

It should be understood that the nonadaptive selector in (\ref{def:selector_for_beta_known}) with $\epsilon$ as in (\ref{cond_epsilon})
and the adaptive selector (\ref{def:selector_for_beta_unknown})--(\ref{def:tau}) with $\epsilon$ as in (\ref{cond_epsilon}) both depend on $\sigma>0$,
and therefore we have a whole class of adaptive selectors indexed by $\sigma$.
Below, we show that, under certain model assumptions, both selectors
achieve almost full selection in model~(\ref{model2}).

\section{Main results}\label{MR}
We first state the conditions when almost full variable selection in model (\ref{model2}) is possible and show that the proposed nonadaptive selector
$\boldsymbol{\hat{\eta}}_k(\beta)$ and adaptive selector $\boldsymbol{\hat{\eta}}_k$  achieve this type of selection. Then, we demonstrate
that our selectors are the best possible in the asymptotically minimax sense.
In the statements of Theorems \ref{th:sim_UB_fixedK_knownBeta} to \ref{th:sim_LB_fixedK} below, $u_k$ is an arbitrary element of $\Ukd$ for $1\leq k\leq d$.

When the level of sparsity $\beta$ is known, we have the following result.

\begin{theorem} \label{th:sim_UB_fixedK_knownBeta}
Let $\beta \in (0,1)$, $\sigma >0$, and $k \in \{1,\ldots,d\}$ be fixed numbers, $d=d_\varepsilon \to \infty$ as $\varepsilon\to 0$, and
\begin{equation} \label{th:assump_epsilon}
	\log {d\choose k}=o \left( \varepsilon^{-{2k}/{(2\sigma + k)}} \right), \quad \text{as } \varepsilon\to 0.
\end{equation}
Assume that the family $r_{\varepsilon,k}>0$ satisfies
\begin{equation}
	\liminf_{\varepsilon \to 0} \frac{a_{\varepsilon,u_k}(r_{\varepsilon,k})}{\sqrt{2\log {d\choose k}}} > \sqrt{\beta}.
	\label{cond:inf}
\end{equation}
Then
\begin{equation*}
	\limsup_{\varepsilon\to 0} \sup_{\etab_k\in H_{k,d}^{\beta}}\sup_{\thetab_k\in  \mathring{\Theta}^\sigma_{k,d}(r_{\varepsilon,k})} {d \choose k}^{\beta-1}{\rm E}_{\thetab_k,\etab_k}|\boldsymbol{\hat{\eta}}_k(\beta)-\boldsymbol{\eta}_k| =0,
\end{equation*}
where $\boldsymbol{\hat{\eta}}_k(\beta)$ is the selector in {\rm (\ref{def:selector_for_beta_known})} with $\epsilon$ as in {\rm (\ref{cond_epsilon})}.
\end{theorem}

When the level of sparsity $\beta$ is unknown, we assume that $\beta\in[b,B]$ for some $0<b<B<1$, where $b$ and $B$ can be arbitrarily close to 0 and 1, respectively.
This is the price that we pay for \textit{adaptive recovery} of the sparsity pattern of $f$.
We claim that the adaptive selector $ \boldsymbol{\hat{\eta}}_k$ of $\etab_k$ in model (\ref{model2}) achieves almost full selection.

\begin{theorem} \label{th:sim_UB_fixedK}
Let $\beta \in [b,B]\subset(0, 1)$, $\sigma >0$, and $k \in \{1,\ldots,d\}$ be fixed numbers, $d=d_\varepsilon \to \infty$ as $\varepsilon\to 0$,
and let condition {\rm (\ref{th:assump_epsilon})} be satisfied.
Assume that the family $r_{\varepsilon,k}>0$ is such that condition {\rm (\ref{cond:inf})} holds true.
Then
\begin{equation*}
	\limsup_{\varepsilon\to 0} \sup_{\etab_k\in H_{k,d}^{\beta}}\sup_{\thetab_k\in  \mathring{\Theta}^\sigma_{k,d}(r_{\varepsilon,k})} {d \choose k}^{\beta-1}{\rm E}_{\thetab_k,\etab_k}|\boldsymbol{\hat{\eta}}_k-\boldsymbol{\eta}_k| =0,
\end{equation*}
where $ \boldsymbol{\hat{\eta}}_k$ is the selector in {\rm (\ref{def:selector_for_beta_unknown})--(\ref{def:tau})}
with $\epsilon$ as in {\rm (\ref{cond_epsilon})}. 
\end{theorem}

Condition (\ref{cond:inf}) imposed on $r_{\varepsilon,k}$ will be named the \textit{selectability condition}. It
ensures that the norms $\|f_{u_k}\|_2$ for $u_k\in \Ukd$ are not too small, and hence the active components  $f_{u_k}$ are {selectable} almost fully.
Theorems \ref{th:sim_UB_fixedK_knownBeta} and  \ref{th:sim_UB_fixedK} show that, under the {selectability condition}~(\ref{cond:inf}), the selection procedures based on $\boldsymbol{\hat{\eta}}_k(\beta)$ and $\boldsymbol{\hat{\eta}}_k$ reconstruct the nonzero elements of $\etab_k\in H_{k,d}^{\beta}$ in such a way that their Hamming errors are small compared to the total number of nonzero elements. In particular, the adaptive selector $\boldsymbol{\hat{\eta}}_k$ provides almost full recovery of $\etab_k$, uniformly  over the sets
$H_{k,d}^{\beta}$ and $ \mathring{\Theta}^\sigma_{k,d}(r_{\varepsilon,k})$, for all $\beta\in[b,B]$ and $\sigma>0$. 
Theorem~\ref{th:sim_UB_fixedK} extends Theorem 3 of \cite{BSt-2017} from $k=1$ to the case $1\leq k\leq d$.

The next theorem shows that if the family $r_{\varepsilon,k}>0$ falls below a certain level,
the minimax normalized  Hamming risk is strictly positive in the limit,  and thus almost full selection in model (\ref{model2}) is impossible.

\begin{theorem} \label{th:sim_LB_fixedK}
Let $\beta \in (0,1)$, $\sigma >0$, and $k \in \{1,\ldots,d\}$ be fixed numbers,  $d=d_\varepsilon \to \infty$ as $\varepsilon \to 0$, and let
condition {\rm (\ref{th:assump_epsilon})} be satisfied.
Assume that the family $r_{\varepsilon,k}>0$ is such that
\begin{equation}
	\limsup_{\varepsilon \to 0} \frac{a_{\varepsilon,{u_k}}(r_{\varepsilon,k})}{\sqrt{2\log {d\choose k}}} < \sqrt{\beta}.
	\label{cond:inf_LB}
\end{equation}
Then
\begin{equation*}
	\liminf_{\varepsilon\to 0}	\inf_{\tilde{\etab}_k} \sup_{\etab_k\in H_{k,d}^{\beta}}\sup_{\thetab_k\in \mathring{\Theta}^\sigma_{k,d}(r_{\varepsilon,k})} {d \choose k}^{\beta-1}{\rm E}_{\thetab_k,\etab_k}|\boldsymbol{\tilde{\eta}}_k-\boldsymbol{\eta}_k| > 0,
\end{equation*}
where the infimum is taken over all selectors $\boldsymbol{\tilde{\eta}}_k = (\tilde{\eta}_{u_k})_{ u_k \in \Ukd}$ 
of $\etab_k=({\eta}_{u_k})_{ u_k \in \Ukd}\in H_{k,d}^{\beta}$ in model~{\rm (\ref{model2})}.
\end{theorem}

When the level of sparsity $\beta$ is unknown, Theorems \ref{th:sim_UB_fixedK} and  \ref{th:sim_LB_fixedK}  ensure that the adaptive selector $\boldsymbol{\hat{\eta}}_k$
is the best possible among all selectors
in model (\ref{model2}) with respect to the (normalized) Hamming risk  in the \textit{asymptotically minimax} sense.
(As in publications \cite{BSt-2017}, \cite{ISt-2014}, \cite{ST-2023}, and \cite{ST-2024}, the optimality of a selection
procedure in this article is understood in the minimax hypothesis testing sense.)

Inequalities {\rm (\ref{cond:inf})} and  {\rm (\ref{cond:inf_LB})} describe the \textit{sharp almost full selection boundary}, which defines a precise demarcation between what is possible and impossible in the problem at hand.
The boundary is determined in terms of the function $a^2_{\varepsilon,u_k}(r_{\varepsilon,k})$ defined in {\rm (\ref{def:a2})}
whose sharp asymptotics for every fixed $k$  are given by {\rm (\ref{aek})}.
Theorems \ref{th:sim_UB_fixedK} and \ref{th:sim_LB_fixedK} augment  Theorems 3.1 and 3.2 of \cite{ST-2023},
where, under similar model assumptions, the \textit{sharp exact selection boundary}
described  by the  inequalities, cf. {\rm (\ref{cond:inf})} and  {\rm (\ref{cond:inf_LB})},  
\begin{equation*}
\liminf_{\varepsilon \to 0}  \frac{a_{\varepsilon,u_k}(r_{\varepsilon,k})}{\sqrt{2\log {d\choose k}}} > 1+\sqrt{1-\beta}
\quad \mbox{and}\quad  \limsup_{\varepsilon \to 0}  \frac{a_{\varepsilon,u_k}(r_{\varepsilon,k})}{\sqrt{2\log {d\choose k}}} < 1+\sqrt{1-\beta}
\end{equation*}
has been established. Parameterizing $a_{\varepsilon,{u_k}}(r_{\varepsilon,k})$ through
$$a_{\varepsilon,{u_k}}(r_{\varepsilon,k})=\sqrt{2\gamma \log{d\choose k}},\quad \gamma>0,$$
yields a ``phase diagram'' for the problem of recovering the sparsity pattern in model {\rm (\ref{model2})}:
(i) if $\gamma>(1+\sqrt{1-\beta})^2$, exact variable selection is possible;
(ii) if $\gamma<(1+\sqrt{1-\beta})^2$, exact variable selection is impossible;
(iii) if $\gamma>\beta$, almost full variable selection is possible;
(iv) if $\gamma<\beta$, neither exact nor almost full selection are possible.
That is, the parameter space $\{(\beta,\gamma)\in \mathbb{R}^2: (\beta,\gamma)\in(0,1)\times (0,4)\}$ can be divided into three regions.
If $(\beta,\gamma)$ are such that exact selection is possible, then we say that $(\beta,\gamma)$ falls within the region of \textit{exact selection}.
If $(\beta,\gamma)$ are such that almost full  selection is possible, then we say that $(\beta,\gamma)$ falls within the region of \textit{almost full selection}.
If $(\beta,\gamma)$ do not fall within the region of exact nor almost full selection, it is said that \textit{no selection is possible}.
The division of the parameter space  into three subregions when
$a_{\varepsilon,{u_k}}(r_{\varepsilon,k})=\sqrt{2\gamma \log{d\choose k}}$ is shown in Figure 1. Note that for $(\beta,\gamma)\in(0,1)\times (4,\infty)$ exact selection is
always possible.
In the context of variable selection, the phase diagram of this kind was obtained for the first time in \cite{GJWY-2012},
see also \cite{BSt-2017}, \cite{GS-2020}, and \cite{MSt-2023}.

\begin{figure}[h]
\centering
\includegraphics[width=12cm,height=7cm]{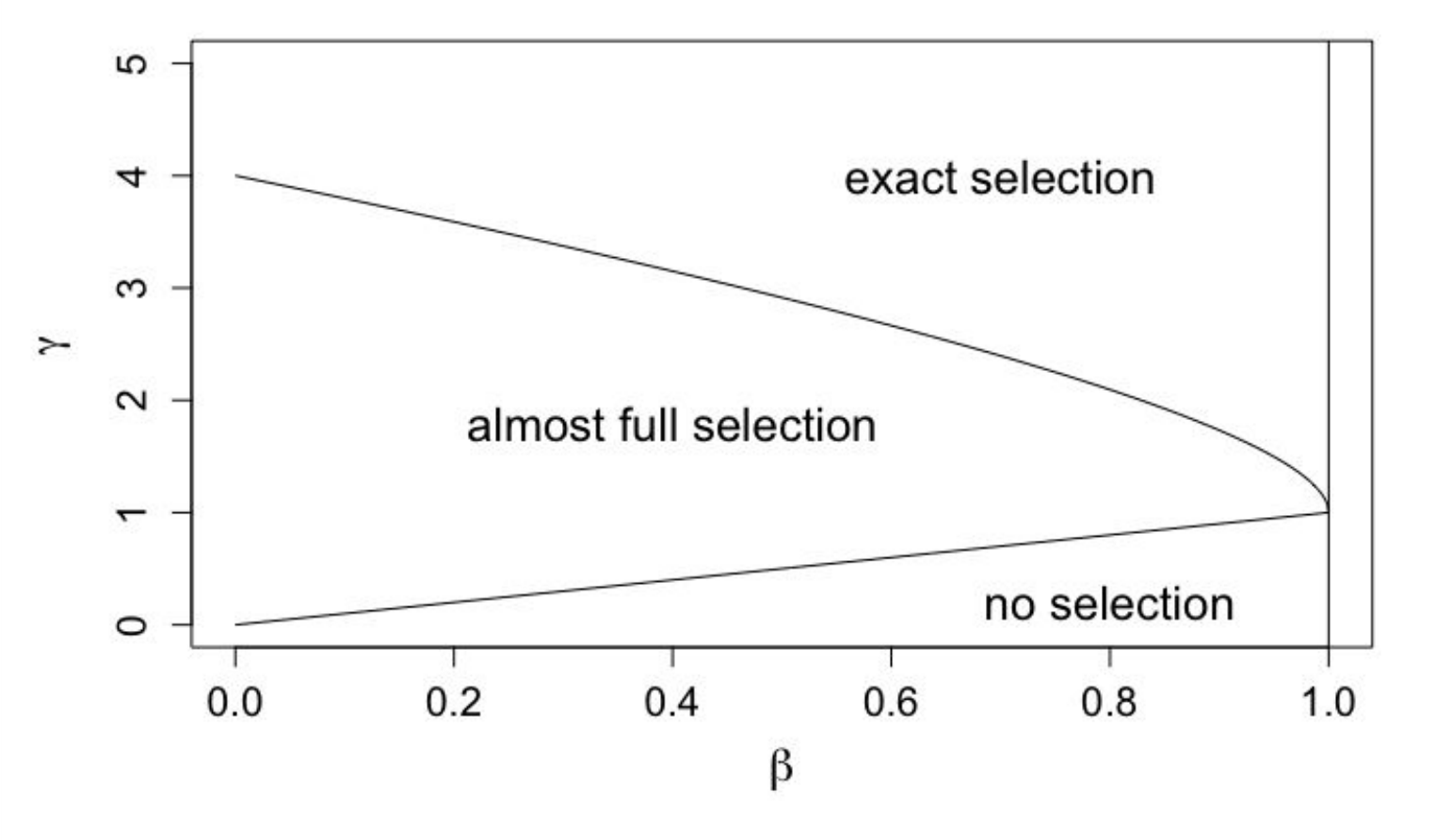}
\caption{Partition of the parameter space into the three regions of variables selection.}	
\label{fig1}
\end{figure}

\section{Simulation study}
In this section, we examine numerically the performance of the almost full selector $ \boldsymbol{\hat{\eta}}_k=(\hat{\eta}_{u_k})_{u_k\in {\cal U}_{k,d}}$
given by {\rm (\ref{def:selector_for_beta_unknown})--(\ref{def:tau})} with $\epsilon$ as in {\rm (\ref{cond_epsilon})}. To this end, we consider the sequence space model (\ref{model2}), in which we take $\varepsilon = 10^{-4}$, $\sigma=1$,  $d=50,100,200$, and $k=2,3$. 
The values for the number of active components $\sum\nolimits_{u_k \in \mathcal{U}_{k,d}}\eta_{u_k}$ are chosen to have $\beta$ ranging from nearly 0.4 to nearly 0.9,
which allows us to numerically study the adaptive selector’s performance across a range of sparsity levels. 
Table~\ref{tab1} contains the values of the sparsity index~$\beta$ and 
the number of active components $\sum\nolimits_{u_k \in \mathcal{U}_{k,d}}\eta_{u_k}$, as well as the  cardinality $\binom{d}{k}$ of the set $\mathcal{U}_{k,d}$ for different
combinations of $d$ and $k$. The $\beta$ values were computed from the equation $\sum\nolimits_{u_k \in \mathcal{U}_{k,d}}\eta_{u_k}={d\choose k}^{1-\beta}$
as $\beta = 1- \log\left({\sum\nolimits_{u_k \in \mathcal{U}_{k,d}} \eta_{u_k}}\right)/{\log {d \choose k}}$. 
\begin{table}[h]
	\centering
	\caption{The sparsity index $\beta$, the number of active components $\sum_{u_k \in \Ukd} \eta_{u_k}$, and the cardinality $\binom{d}{k}$ for all combinations of $d = 50,100,200$ and $k = 2,3$.}
	\label{tab1}
    \begin{tabular}{|c|ccc|ccc|}
    \hline
     & \multicolumn{3}{c|}{$k=2$} & \multicolumn{3}{c|}{$k=3$} \\ \cline{2-7} 
    $d$ & \multicolumn{1}{c|}{$\beta$} & \multicolumn{1}{c|}{$\sum \eta_{u_k}$} & $\binom{d}{2}$ & \multicolumn{1}{c|}{$\beta$} & \multicolumn{1}{c|}{$\sum \eta_{u_k}$} & $\binom{d}{3}$ \\ \hline
    50 & \multicolumn{1}{c|}{0.4005} & \multicolumn{1}{c|}{71} & 1225 & \multicolumn{1}{c|}{0.7021} & \multicolumn{1}{c|}{19} & 19600 \\
    100 & \multicolumn{1}{c|}{0.5006} & \multicolumn{1}{c|}{70} & 4950 & \multicolumn{1}{c|}{0.8001} & \multicolumn{1}{c|}{11} & 161700 \\
    200 & \multicolumn{1}{c|}{0.6008} & \multicolumn{1}{c|}{52} & 19900 & \multicolumn{1}{c|}{0.9016} & \multicolumn{1}{c|}{4} & 1313400 \\ \hline
    \end{tabular}
\end{table}

Consider the following functions defined on the interval $[0,1]$: 
\begin{align*}
g_1(t)&=t^2\left(2^{t-1}-(t-0.5)^2\right)e^t-0.5424,
& g_8(t)&=10(t^2-0.5)^5+0.068,\\
g_2(t)&=t^2\left(2^{t-1}-(t-1)^5\right)-0.2887,
& g_{9}(t)&=6t^2-6t+1,\\
g_3(t)&=0.1\left(15t^2 2^{t-1}\cos(15t)-0.5011\right),
& g_{10}(t)&=(2t-1)^3,\\
g_4(t)&=2t-1, & g_{11}(t)&=(2t-1)^5,\\
g_5(t)&=5(t-0.7)^3+0.29, & g_{12}(t)&=3t^2-1,\\
 g_6(t)&=2(t-0.4)^2-0.1867,
& g_{13}(t)&=3\sqrt{t}-2,\\
g_7(t)&=0.7(t^2-0.1)^3-0.0643,
& 
\end{align*}
and let the active components $f_{u_k}$ of $f$ defined on $[0,1]^k$ be as follows. For $k=2$, 
consider ${13\choose 2}=78$ pairs $\{i, j\}$, $1 \leq i<j \leq 13$,  ordered in the following way:
\begin{align*}
    \{1,2\},\{1,3\}, \ldots,\{1,13\},\{2,3\}, \ldots,\{2,13\}, \ldots,\{12,13\} .
\end{align*}
For $d=50, 100, 200$, we take the first $71$, $70$, $52$ pairs above to be the \textit{active set}, respectively. 
Then for $u_k = \{i, j\} \in \mathcal{U}_{2, d}$, we set $f_{u_k}\left(\mathbf{t}_{u_k}\right)= g_i\left(t_i\right) g_j\left(t_j\right)$ if 
$u_k$ belongs to the active set, 
and $f_{u_k}\left( \mathbf{t}_{u_k}\right)=0$ otherwise. For $k=3$, consider ${13\choose 3}=286$ triples $\{i, j, l\},$ $1 \leq i<j<l \leq 13$, ordered as follows:
\begin{align*}
    \{1,2,3\},\{1,2,4\}, \ldots,\{1,2,13\},\{1,3,4\}, \ldots,\{1,3,13\}, \ldots,\{11,12,13\}.
\end{align*}
For $d=50, 100, 200$, we take the first $19$, $11$, $4$ triples above to be the active set,  respectively. 
For $u_k=\{i, j, l\} \in \mathcal{U}_{3, d}$, we set  $ f_{u_k}\left(\mathbf{t}_{u_k}\right)=g_i\left(t_i\right) g_j\left(t_j\right) g_{l}\left(t_{l}\right)$ if $u_k$ belongs to the active set, 
and $f_{u_k}\left( \mathbf{t}_{u_k}\right)=0$ otherwise. For $k=2,3$ and $u_k\in {\cal U}_{k,d}$, the condition $\int_0^1 f_{u_k} (\tb_{u_k})\, dt_j =0$ 
holds true up to four decimal places for all $j\in u_k$.
Note that for $k=2$ the active set is large and hence the signal  $f=\sum_{u_k\in {\cal U}_{k,d}} \eta_{u_k}f_{u_k}$ is relatively dense, whereas for $k=3$ 
the active set is small and hence the signal is highly sparse.

For smooth Sobolev functions, the absolute values of their Fourier coefficients decay to zero at a~polynomial rate.
Therefore, although in theory $\lb=(l_1,\ldots,l_d) \in \mathring{\mathbb{Z}}_{u_k}$,
we shall restrict ourselves to  $\lb \in \Ztu \deq  \mathring{\mathbb{Z}}_{u_k} \cap [-n,n]^d,$
where for all values of $d$ under study we take $n=344$ for $k=2$ and  $n=127$ for $k=3.$
The chosen values of $n$ ensure that none of the nonzero coefficients $\theta^*_{\lb}(r^*_{\varepsilon,k,m})$, and hence none of the
nonzero weights $\omega_{\lb}(r^*_{\varepsilon,k,m})$, is missing in the evaluation of the statistics $S_{u_k}(\beta_{k,\hat{m}_k})$, $u_k\in {\cal U}_{k,d}$, in
the definition of $\hat{\eta}_{u_k}$ in (\ref{def:selector_for_beta_unknown}).
In this definition, we choose $ \epsilon = 2.4\log^{-1/2} \binom{d}{k}$  to satisfy condition (\ref{cond_epsilon}).
The random index  $\hat{m}_k$ is determined by Lepski's method in accordance with (\ref{Lepski_m})--(\ref{def:tau}).
The grid of equidistant points $\beta_{k,m}$, $m = 1, \ldots , M_k$, which is used to obtain  $\hat{m}_k$,
is chosen as in Section~\ref{Selector} with $b=0.001,$ $B=0.999$, and $M_k=25$. Thus, the model we are dealing with in this section  is as follows, cf. (\ref{model2}):
\begin{equation*}
	X_{\lb}=\eta_{u_k}\theta_\lb(u_k)+\varepsilon \xi_{\lb},\quad \lb \in \Ztu,\quad u_k\in{\cal U}_{k,d},
\end{equation*}
where  the component $\eta_{u_k} $ of $\etab_k=(\eta_{u_k})_{u_k\in{\cal U}_{k,d}}\in{H}^{\beta}_{k,d}$ equals 1
if $u_k$ belongs to the active set, 
and zero otherwise.
Note that, in all cases under study, $\log{d\choose k}$ is much smaller than $\varepsilon^{-2k/(2\sigma+k)}$, and thus condition (\ref{th:assump_epsilon}) of Theorem~\ref{th:sim_UB_fixedK} is satisfied.

For  $d=50, 100, 200$ and $k=2,3$, we run $J=100$ independent cycles of simulations and
estimate the normalized Hamming risk 
${d\choose k}^{\beta-1}{\rm E}_{\thetab_k,\etab_k}\left(\sum_{u_k\in{\cal U}_{k,d}}|\hat{\eta}_{u_k}-\eta_{u_k}|\right)$
by means of the quantity
$${\rm Err}(\boldsymbol{\hat{\eta}}_k):=\frac{1}{J}\sum_{j=1}^J{d\choose k}^{\beta-1}\sum_{u_k \in \mathcal{U}_{k,d}}|\hat{\eta}_{u_k}^{(j)}-\eta_{u_k}|,$$
where $\hat{\eta}_{u_k}^{(j)}$ is the value of $\hat{\eta}_{u_k}$ obtained in the $j$th repetition of the experiment. 
We shall use ${\rm Err}(\boldsymbol{\hat{\eta}}_k)$ 
to judge the quality of the adaptive selector $\boldsymbol{\hat{\eta}}_k$ under the settings described below,
and  numerically illustrate the impact of the \textit{sharp almost full selection boundary}, given by inequalities \eqref{cond:inf} and (\ref{cond:inf_LB}),
on the ability of $ \boldsymbol{\hat{\eta}}_k $ to select almost fully the active components of  $f=\sum_{u_k\in {\cal U}_{k,d}} \eta_{u_k}f_{u_k}$. 

To do this, we modify the signal under study by multiplying the active component $f_{u_k}$, with $u_k=\{1,2\}$ for $k=2$ and $u_k=\{1,2,3\}$ for $k=3$, by a positive constant $\alpha$, while keeping the other active components unchanged. First, take $k=2$ and consider the  active component $f_{\{1,2\}}=g_1 g_2$ of $f$,
whose $L_2$-norm after rescaling is $\left\|\alpha f_{\{1,2\}}\right\|_2 =\alpha \left\| {g_1} \right\|_{2} \left\| {g_2} \right\|_{2}$, and the condition
$\|f_{\{1,2\}}\|_2\geq r_{\varepsilon,k}$  implies $\left\|\alpha f_{\{1,2\}}\right\|_2\geq \alpha r_{\varepsilon,k}.$
In the worst case scenario, when $\left\|\alpha f_{\{1,2\}}\right\|_2=\alpha r_{\varepsilon,k},$ 
in view of the asymptotic expression in \eqref{aek}, we can write for $u_k=\{1,2\}$ as $\varepsilon\to 0$
\begin{align}\label{alphar}
	a_{\varepsilon, u_k}\left( \alpha r_{\varepsilon, k} \right) \sim C(\sigma, k) 
\left(\alpha \left\| {g_1} \right\|_{2} \left\| {g_2} \right\|_{2}\right)^{2+k /(2 \sigma)}  \varepsilon^{-2},
\end{align}
with constant $C(\sigma, k)$ as  in \eqref{aek}. The value $\alpha^*=\alpha^*_{\varepsilon,k,d}$ of $\alpha$ that makes the signal component $\alpha f_{\{1,2\}}$ with $\alpha<\alpha^*$ non-selectable  is
obtained by equating the right-hand side of (\ref{alphar}) to $\sqrt{2 \beta \log \binom{d}{k}}$ (see condition (\ref{cond:inf_LB})) and solving the resulting equation for $\alpha$. This gives
\begin{align*}
	\alpha^* = \left( \left\| {g_1} \right\|_{2} \left\| g_{2} \right\|_{2} \right)^{-1} \left( \frac{\varepsilon^2\sqrt{2 \beta \log \binom{d}{k}}  }{C(\sigma, k) } \right)^{{2\sigma}/{(4 \sigma+k)}}.
\end{align*}
For $\alpha<\alpha^*$, we have $a_{\varepsilon, u_k}( \alpha r_{\varepsilon, k} ) / \sqrt{2 \log \binom{d}{k}}  < \sqrt{ \beta}$, and hence
the norm $\|\alpha f_{\{1,2\}}\|_2$, 
which in this case falls below the almost full selection boundary,
is too small for $\alpha f_{\{1,2\}}$ to be selectable almost fully. 
On the other hand,  for $\alpha>$ $\alpha^*$, we have $a_{\varepsilon, u_k}( \alpha r_{\varepsilon, k} ) / \sqrt{2 \log \binom{d}{k}} > \sqrt{ \beta}$, and hence 
the norm $\|\alpha f_{\{1,2\}}\|_2$,
which now falls above the almost full selection boundary,  is
large enough for $\alpha f_{\{1,2\}}$ to be selectable almost fully.
Similar arguments apply to the active component $\alpha f_{u_k}$ with $u_k=\{1,2,3\}$ for $k=3$.

\begin{figure}[h]
	\centering
    \includegraphics[width=14cm]{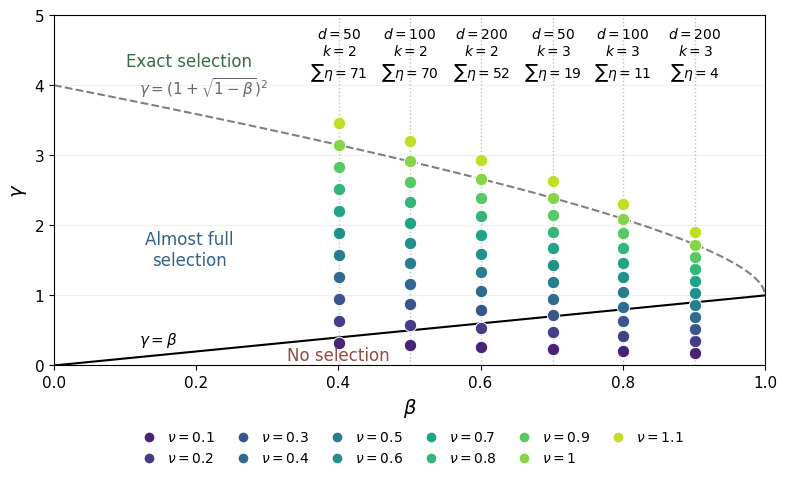}    
		\caption{Simulation grid of signal strengths superimposed on the $(\beta,\gamma)$ phase diagram in Figure \ref{fig1}. The solid and dashed curves are the almost full and exact selection boundaries, respectively.}
			\label{fig2}
\end{figure}

Now, recall the parametrization of $a_{\varepsilon,u_k}(r_{\varepsilon,k})$ through  
$a_{\varepsilon,u_k}(r_{\varepsilon,k})=\sqrt{2\gamma\log {d\choose k}}$, $\gamma>0$, introduced in Section \ref{MR}.
Then, graphically, in terms of the phase diagram depicted in Figure \ref{fig1}, 
varying $\alpha$ shifts the value of the signal strength (i.e. the $L_2$-norm of $\alpha f_{\{1,2\}}$ and $\alpha f_{\{1,2,3\}}$ for $k=2$ and 
 $k=3$, respectively) along the $\gamma$ axis, while $\beta$ is kept fixed (see Figure \ref{fig2}).
Aiming to examine the impact of the signal strength inside and near the almost full selection region, 
we shall use the exact selection boundary given by $\gamma=\left(1+\sqrt{1-\beta}  \right)^2$ as a reference for choosing $\alpha$. 
Specifically, for each value of $\beta$ under study, we set 
\begin{align*} 
	\alpha^*(\nu) = \left( \prod_{i=1}^k \| {g_i} \|_{2} \right)^{-1} \left( \frac{\varepsilon^2\sqrt{2 \gamma^*(\nu)  \log \binom{d}{k}}  }{C(\sigma, k) } \right)^{{2\sigma}/{(4 \sigma+k)}},\quad k=2,3,
\end{align*}
where $\gamma^*(\nu) =\nu\left(1+\sqrt{1-\beta}\right)^2$ for $\nu \in \{0.1,0.2, \ldots, 1.1\}$,
and obtain  a grid of signal strengths being below, on, and just above the exact selection boundary, as shown in Figure \ref{fig2}. 
In Figure \ref{fig2}, for each value of $\beta$, increasing $\nu$ from 0.1 to 1.1 shifts the norm $\|\alpha^*(\nu)f_{u_k}\|_2$, with $u_k=\{1,2\}$ for $k=2$ and $u_k=\{1,2,3\}$ for $k=3$, from the region of no selection, through the region of almost full selection, to the region of exact selection. 
For different combinations of $d$ and $k$, the values of  $\alpha^*$ and  $\alpha^*(\nu)$, as defined above, are shown in Table $\ref{tab2}$.

\begin{table}[h!]
\centering
\caption{Signal multipliers $\alpha^*(\nu)$ corresponding to
$\gamma^*(\nu)=\nu(1+\sqrt{1-\beta})^2$, with
$\varepsilon=10^{-4}$ and $\sigma=1$.
The $\alpha^*$ row corresponds to the almost full
selection boundary $\gamma=\beta$, while the row $\nu=1.0$
corresponds to the exact selection boundary $\gamma=\left(1+\sqrt{1-\beta}  \right)^2$.
Cells shaded in pink, blue, and green indicate that
$\alpha^*(\nu)<\alpha^*$,
$\alpha^*<\alpha^*(\nu)\leq\alpha^*(1)$,
and $\alpha^*(\nu)>\alpha^*(1)$, respectively.}
\label{tab2}
\begin{tabular}{|c|cccccc|}
	\hline
	& \multicolumn{3}{c|}{$k=2$}
	& \multicolumn{3}{c|}{$k=3$} \\ \cline{2-7}
	
	$d$ & 50 & 100 & \multicolumn{1}{c|}{200}
	& 50 & 100 & 200 \\
	
	$\beta$ & 0.4005 & 0.5006 & \multicolumn{1}{c|}{0.6008}
	& 0.7021 & 0.8001 & 0.9016 \\
	
	$\alpha^*$ & 0.0146 & 0.0156 & \multicolumn{1}{c|}{0.0165}
	& 0.0837 & 0.0876 & 0.0912 \\ \hline
	
	$\nu$ & \multicolumn{6}{c|}{$\alpha^*(\nu)$} \\ \hline
	
	1.1
	& \cellcolor{regionthree}{0.0209}
	& \cellcolor{regionthree}{0.0212}
	& \multicolumn{1}{c|}{\cellcolor{regionthree}{0.0214}}
	& \cellcolor{regionthree}{0.1010}
	& \cellcolor{regionthree}{0.1019}
	& \cellcolor{regionthree}{0.1014} \\ \cline{2-7}
	
	1.0
	& \cellcolor{regiontwo}{0.0205}
	& \cellcolor{regiontwo}{0.0209}
	& \multicolumn{1}{c|}{\cellcolor{regiontwo}{0.0211}}
	& \cellcolor{regiontwo}{0.0997}
	& \cellcolor{regiontwo}{0.1005}
	& \cellcolor{regiontwo}{0.1001} \\
	
	0.9
	& \cellcolor{regiontwo}{0.0202}
	& \cellcolor{regiontwo}{0.0205}
	& \multicolumn{1}{c|}{\cellcolor{regiontwo}{0.0207}}
	& \cellcolor{regiontwo}{0.0982}
	& \cellcolor{regiontwo}{0.0990}
	& \cellcolor{regiontwo}{0.0986} \\
	
	0.8
	& \cellcolor{regiontwo}{0.0198}
	& \cellcolor{regiontwo}{0.0201}
	& \multicolumn{1}{c|}{\cellcolor{regiontwo}{0.0203}}
	& \cellcolor{regiontwo}{0.0965}
	& \cellcolor{regiontwo}{0.0974}
	& \cellcolor{regiontwo}{0.0969} \\
	
	0.7
	& \cellcolor{regiontwo}{0.0193}
	& \cellcolor{regiontwo}{0.0197}
	& \multicolumn{1}{c|}{\cellcolor{regiontwo}{0.0199}}
	& \cellcolor{regiontwo}{0.0947}
	& \cellcolor{regiontwo}{0.0955}
	& \cellcolor{regiontwo}{0.0951} \\
	
	0.6
	& \cellcolor{regiontwo}{0.0188}
	& \cellcolor{regiontwo}{0.0192}
	& \multicolumn{1}{c|}{\cellcolor{regiontwo}{0.0194}}
	& \cellcolor{regiontwo}{0.0926}
	& \cellcolor{regiontwo}{0.0935}
	& \cellcolor{regiontwo}{0.0930} \\ \cline{7-7}
	
	0.5
	& \cellcolor{regiontwo}{0.0183}
	& \cellcolor{regiontwo}{0.0186}
	& \multicolumn{1}{c|}{\cellcolor{regiontwo}{0.0188}}
	& \cellcolor{regiontwo}{0.0903}
	& \multicolumn{1}{c|}{\cellcolor{regiontwo}{0.0911}}
	& \cellcolor{regionone}{0.0906} \\
	
	0.4
	& \cellcolor{regiontwo}{0.0176}
	& \cellcolor{regiontwo}{0.0179}
	& \multicolumn{1}{c|}{\cellcolor{regiontwo}{0.0181}}
	& \cellcolor{regiontwo}{0.0874}
	& \multicolumn{1}{c|}{\cellcolor{regiontwo}{0.0882}}
	& \cellcolor{regionone}{0.0878} \\ \cline{6-6}
	
	0.3
	& \cellcolor{regiontwo}{0.0168}
	& \cellcolor{regiontwo}{0.0171}
	& \multicolumn{1}{c|}{\cellcolor{regiontwo}{0.0173}}
	& \multicolumn{1}{c|}{\cellcolor{regiontwo}{0.0839}}
	& \cellcolor{regionone}{0.0847}
	& \cellcolor{regionone}{0.0843} \\ \cline{4-5}
	
	0.2
	& \cellcolor{regiontwo}{0.0157}
	& \multicolumn{1}{c|}{\cellcolor{regiontwo}{0.0160}}
	& \multicolumn{1}{c|}{\cellcolor{regionone}{0.0161}}
	& \cellcolor{regionone}{0.0792}
	& \cellcolor{regionone}{0.0799}
	& \cellcolor{regionone}{0.0795} \\ \cline{2-3}
	
	0.1
	& \cellcolor{regionone}{0.0140}
	& \cellcolor{regionone}{0.0142}
	& \multicolumn{1}{c|}{\cellcolor{regionone}{0.0144}}
	& \cellcolor{regionone}{0.0717}
	& \cellcolor{regionone}{0.0724}
	& \cellcolor{regionone}{0.0720} \\ \hline
\end{tabular}
\end{table}

\begin{figure}[h!]
	\centering
	\includegraphics[width=1\textwidth]{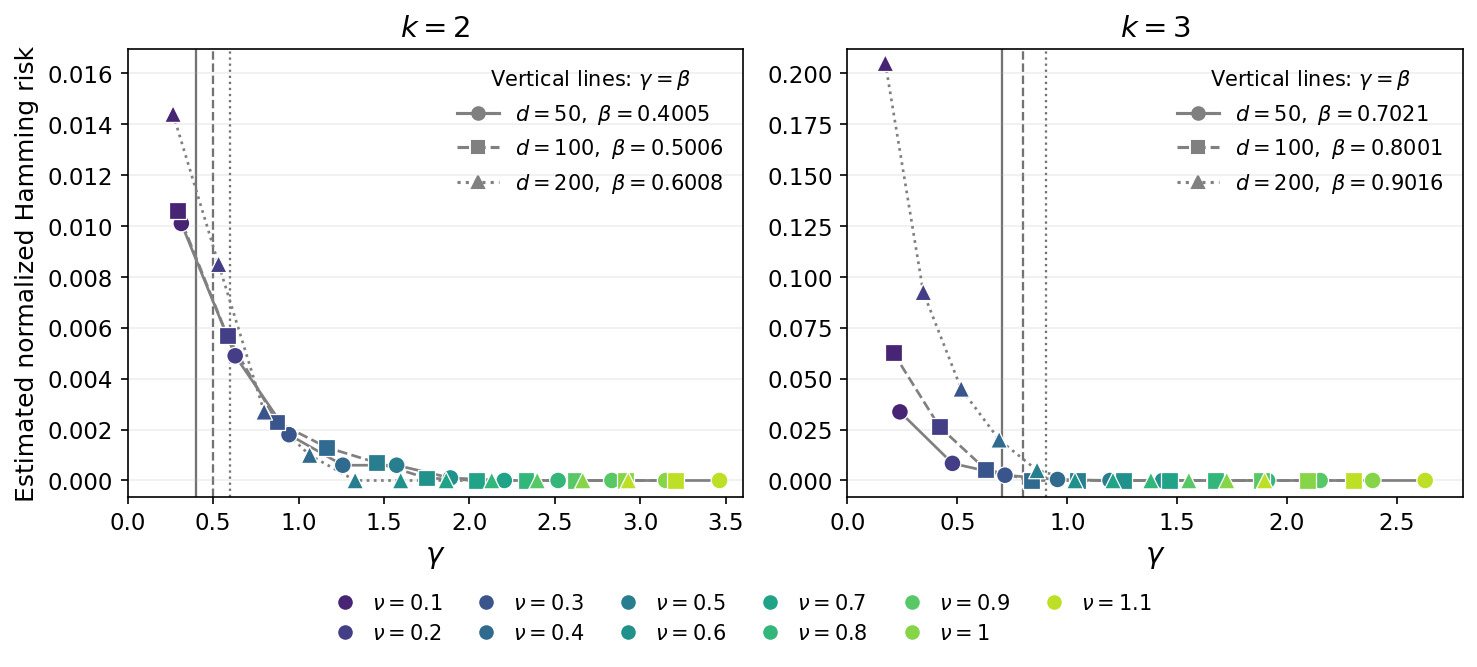}
		\caption{A graphical representation of the estimated normalized Hamming risk $\operatorname{Err}(\boldsymbol{\hat{\eta}}_k)$  with vertical lines at $\gamma=\beta$.}
			\label{fig3}
\end{figure}

The values of the estimated normalized Hamming risk $\operatorname{Err}(\boldsymbol{\hat{\eta}}_k)$ 
obtained from $J=100$ simulation cycles are reported in Table $\ref{tab3}$.  
The same values of $\operatorname{Err}(\boldsymbol{\hat{\eta}}_k)$ 
as a function of $\gamma$ are sketched in Figure~\ref{fig3},
in which the vertical lines indicate the almost full selection boundary $\gamma=\beta$.
The use of $\gamma$ allows us to make a connection between the plots in Figure  \ref{fig3} and the phase diagram in Figure \ref{fig1}. 
The values of $\operatorname{Err}(\boldsymbol{\hat{\eta}}_k)$ 
for different $\nu \in \{0.1,0.2, \ldots, 1.1\}$ are specified in Table \ref{tab3}, 
where the entries shaded in pink correspond to 
the no selection region when $\gamma^*(\nu)<\beta$.

\begin{table}[h!]
    \centering
    \setlength{\tabcolsep}{2.45pt}
    \caption{The estimated normalized Hamming risk $\operatorname{Err}(\boldsymbol{\hat{\eta}}_k)$ obtained from $J=100$ simulation cycles. 
    Cells shaded in pink, blue, and green  indicate that
    $\gamma^*(\nu)<\beta$, $\beta<\gamma^*(\nu)\leq(1+\sqrt{1-\beta})^2$, and $\gamma^*(\nu)>(1+\sqrt{1-\beta})^2$, respectively.}
    \label{tab3}
    \begin{tabular}{|c|c|c|ccccccccccc|}
    	\hline
    	& & & \multicolumn{11}{c|}{$\nu$} \\ \cline{4-14}
    	$k$ & $d$ & $\beta$
    	& 0.1 & 0.2 & 0.3 & 0.4 & 0.5 & 0.6
    	& 0.7 & 0.8 & 0.9 & 1.0 & 1.1 \\ \hline
    
    	\multirow{3}{*}{2}
    	& 50 & 0.4005
    	& \multicolumn{1}{c|}{\cellcolor{regionone}{0.0101}}
    	& \cellcolor{regiontwo}{0.0049}
    	& \cellcolor{regiontwo}{0.0018}
    	& \cellcolor{regiontwo}{0.0006}
    	& \cellcolor{regiontwo}{0.0006}
    	& \cellcolor{regiontwo}{0.0001}
    	& \cellcolor{regiontwo}{\tablezero}
    	& \cellcolor{regiontwo}{\tablezero}
    	& \cellcolor{regiontwo}{\tablezero}
    	& \multicolumn{1}{c|}{\cellcolor{regiontwo}{\tablezero}}
    	& \cellcolor{regionthree}{\tablezero} \\
    
    	& 100 & 0.5006
    	& \multicolumn{1}{c|}{\cellcolor{regionone}{0.0106}}
    	& \cellcolor{regiontwo}{0.0057}
    	& \cellcolor{regiontwo}{0.0023}
    	& \cellcolor{regiontwo}{0.0013}
    	& \cellcolor{regiontwo}{0.0007}
    	& \cellcolor{regiontwo}{0.0001}
    	& \cellcolor{regiontwo}{0}
    	& \cellcolor{regiontwo}{0}
    	& \cellcolor{regiontwo}{0}
    	& \multicolumn{1}{c|}{\cellcolor{regiontwo}{0}}
    	& \cellcolor{regionthree}{0} \\ \cline{5-5}
    
    	& 200 & 0.6008
    	& \cellcolor{regionone}{0.0144}
    	& \multicolumn{1}{c|}{\cellcolor{regionone}{0.0085}}
    	& \cellcolor{regiontwo}{0.0027}
    	& \cellcolor{regiontwo}{0.0010}
    	& \cellcolor{regiontwo}{0}
    	& \cellcolor{regiontwo}{0}
    	& \cellcolor{regiontwo}{0}
    	& \cellcolor{regiontwo}{0}
    	& \cellcolor{regiontwo}{0}
    	& \multicolumn{1}{c|}{\cellcolor{regiontwo}{0}}
    	& \cellcolor{regionthree}{0} \\ \hline
    
    	\multirow{3}{*}{3}
    	& 50 & 0.7021
    	& \cellcolor{regionone}{0.0337}
    	& \multicolumn{1}{c|}{\cellcolor{regionone}{0.0084}}
    	& \cellcolor{regiontwo}{0.0026}
    	& \cellcolor{regiontwo}{0.0005}
    	& \cellcolor{regiontwo}{0}
    	& \cellcolor{regiontwo}{0}
    	& \cellcolor{regiontwo}{0}
    	& \cellcolor{regiontwo}{0}
    	& \cellcolor{regiontwo}{0}
    	& \multicolumn{1}{c|}{\cellcolor{regiontwo}{0}}
    	& \cellcolor{regionthree}{0} \\ \cline{6-6}
    
    	& 100 & 0.8001
    	& \cellcolor{regionone}{0.0627}
    	& \cellcolor{regionone}{0.0264}
    	& \multicolumn{1}{c|}{\cellcolor{regionone}{0.0055}}
    	& \cellcolor{regiontwo}{0}
    	& \cellcolor{regiontwo}{0}
    	& \cellcolor{regiontwo}{0}
    	& \cellcolor{regiontwo}{0}
    	& \cellcolor{regiontwo}{0}
    	& \cellcolor{regiontwo}{0}
    	& \multicolumn{1}{c|}{\cellcolor{regiontwo}{0}}
    	& \cellcolor{regionthree}{0} \\ \cline{7-8}
    
    	& 200 & 0.9016
    	& \cellcolor{regionone}{0.2050}
    	& \cellcolor{regionone}{0.0925}
    	& \cellcolor{regionone}{0.0450}
    	& \cellcolor{regionone}{0.0200}
    	& \multicolumn{1}{c|}{\cellcolor{regionone}{0.0050}}
    	& \cellcolor{regiontwo}{0}
    	& \cellcolor{regiontwo}{0}
    	& \cellcolor{regiontwo}{0}
    	& \cellcolor{regiontwo}{0}
    	& \multicolumn{1}{c|}{\cellcolor{regiontwo}{0}}
    	& \cellcolor{regionthree}{0} \\ \hline
    \end{tabular}
\end{table}

\newpage 
Figure $\ref{fig3}$ shows that when the signal strength of $\alpha^*(\nu)f_{u_k}$, with $u_k=\{1,2\}$ for $k=2$ and $u_k=\{1,2,3\}$ for $k=3$, lies above the almost full selection boundary, the estimated normalized Hamming risk $\operatorname{Err}(\boldsymbol{\hat{\eta}}_k)$ is close to zero, indicating that $\boldsymbol{\hat{\eta}}_k$ achieves almost full selection. 
Conversely, when the signal strength at hand falls below the almost full selection boundary, variable selection becomes more difficult, resulting in 
larger values of $\operatorname{Err}(\boldsymbol{\hat{\eta}}_k)$. The highest observed level of 
$\operatorname{Err}(\boldsymbol{\hat{\eta}}_k)$ is obtained when the signal strength falls substantially below the almost full selection boundary.

It is seen from Table \ref{tab3} and Figure \ref{fig3} that, in general, the stronger the signal, the smaller the estimated risk ${\rm Err}(\boldsymbol{\hat{\eta}}_k)$.
It is also seen that almost full selection gets harder as $\beta$ gets larger, or as the model gets sparser. This is consistent with the exact selection numerical results presented in Section 4 of \cite{ST-2023}.
Note, however, that the exact selector proposed in \cite{ST-2023}
never detects a signal if there is none (see Section 4 of \cite{ST-2023}). At the same time, the adaptive almost full selector
$ \boldsymbol{\hat{\eta}}_k=(\hat{\eta}_{u_k})_{u_k\in {\cal U}_{k,d}}$
given by {\rm (\ref{def:selector_for_beta_unknown})--(\ref{def:tau})} with $\epsilon$ as in {\rm (\ref{cond_epsilon})} 
may produce a small number of false positives. 
This is, however, not surprising since, by definition, the almost full selector is less precise than the exact selector whose Hamming risk is algebraically small for small $\varepsilon$ (see Theorem 3.1 of \cite{ST-2023}).
Overall, the numerical results of this section are in agreement with the analytical 
findings of Theorems \ref{th:sim_UB_fixedK} and \ref{th:sim_LB_fixedK} in Section \ref{MR}.

\section{Extensions}\label{Ext}

In this section, we  shall extend Theorems  \ref{th:sim_UB_fixedK} and  \ref{th:sim_LB_fixedK} to a more general and interesting
case when the $d$-variate function $f$ has a more complex sparse structure and is decomposed as in  (\ref{fun_sub_anova2})
rather than as in  (\ref{f}), and the model under study is given by relations (\ref{model1}), (\ref{fun_sub_anova2})--(\ref{orthcon1}), and
relations (\ref{sparsitycond}) and  (\ref{Fu}) for all $1\leq k\leq s$.
The  sequence space model that is equivalent to this model 
is, cf. (\ref{model2}),
\begin{eqnarray*}
	X_{\lb}&=&(f_\emptyset,\phi_{\lb})_2+\varepsilon\xi_{\lb}=f_\emptyset+\varepsilon\xi_{\lb},\quad \lb\in\mathring{\mathbb{Z}}_\emptyset,\nonumber\\
	X_{\lb}&=&(f,\phi_{\lb})_2+\varepsilon \xi_{\lb}=(f_\emptyset,\phi_{\lb})_2+\sum_{v_k\in \Ukd}\eta_{v_k}(f_{v_k},\phi_{\lb})_2+
	\varepsilon \xi_{\lb}\label{model3_s} \\ &=&\eta_{u_k}\theta_{\lb}(u_k)+\varepsilon \xi_{\lb},\quad 
	\lb\in  \mathring{\mathbb{Z}}_{u_k},\quad u_k\in \Ukd,\quad 1\leq k\leq s,\nonumber 
\end{eqnarray*}
where  $X_\lb=\mathbb{X}_\varepsilon(\phi_{\lb})$ is the $\lb$th empirical Fourier coefficient, 
$\etab_{k}=(\eta_{u_k})_{u_k \in \Ukd}\in H_{k,d}^{\beta}$,
$\thetab_{k}=(\theta_{\lb}(u_k))_{\lb\in  \mathring{\mathbb{Z}}_{u_k}}$ $\in {\Theta}_{\cb_{u_k}}$,
and $\xi_{\lb}=\mathbb{W}(\phi_{\lb})$ are iid standard normal random variables for all $\lb\in  \mathring{\mathbb{Z}}_{u_k}$, $u_k\in \Ukd,$ $0\leq k\leq s$. 
Since making inference about $f_\emptyset={\rm const}$ (either estimating $f_\emptyset$ or testing the hypothesis $H_0: f_\emptyset=0$) based on the observation  $X_{\lb}=f_\emptyset+\varepsilon\xi_{\lb}$, $\lb\in \mathring{\mathbb{Z}}_\emptyset$, 
is, in general, a simple problem and is not of interest in this article, the sequence space model under study takes the form 
\begin{eqnarray}\label{model2_s}
		X_{\lb}	=\eta_{u_k}\theta_{\lb}(u_k)+\varepsilon \xi_{\lb},\quad 
	\lb\in  \mathring{\mathbb{Z}}_{u_k},\quad u_k\in \Ukd,\quad 1\leq k\leq s,
\end{eqnarray}
where $\etab_{k}=(\eta_{u_k})_{u_k\in \Ukd}$, $\thetab_{k}=(\theta_{\lb}(u_k))_{\lb\in  \mathring{\mathbb{Z}}_{u_k}}$,
and $\{\xi_{\lb}\}_{\lb\in  \mathring{\mathbb{Z}}_{u_k}}$ are as  above.

For ${u_k}\in \Ukd$, $1\leq k\leq s$, denote as before $\Xb_{\!u_k}=(X_{\lb})_{\lb\in \mathring{\mathbb{Z}}_{u_k}}$
and let $\boldsymbol{\hat{\eta}}_k=(\hat{\eta}_{u_k})_{u_k\in \Ukd}$, where 
$\hat{\eta}_{u_k}=\hat{\eta}_{u_k}(\Xb_{\!u_k})\in \{0,1\}$, be an estimator of 
$\boldsymbol{{\eta}}_k=({\eta}_{u_k})_{u_k\in \Ukd}\in H_{k,d}^{\beta}$.
An  estimator $\boldsymbol{\hat{\eta}}=(\boldsymbol{\hat{\eta}}_1,\ldots,\boldsymbol{\hat{\eta}}_s)$ for
$\boldsymbol{{\eta}}=(\boldsymbol{{\eta}}_1,\ldots,\boldsymbol{{\eta}}_s)\in {\cal H}^{\beta}_{s,d}$ will be named an \textit{aggregate selector}.
Let ${\rm E}_{\thetab,\boldsymbol{\eta}}$ denote the expectation with respect to the joint distribution of 
$\Xb_{\!u_k}=(X_{\lb})_{\lb\in \mathring{\mathbb{Z}}_{u_k}}$, $u_k\in \Ukd$, $1\leq k\leq s$, in model (\ref{model2_s});
the expectation ${\operatorname{E}}_{\thetab_k,\etab_k} $ is defined accordingly. 
When dealing with the problem of selecting relevant $\eta_{u_k}$ in model (\ref{model2_s}) almost fully, 
we consider the \textit{maximum normalized Hamming risk}
of an aggregate selector $\boldsymbol{\hat{\eta}}=(\boldsymbol{\hat{\eta}}_1,\ldots,\boldsymbol{\hat{\eta}}_s)$ 
given by
\begin{eqnarray}
	\mathcal{R}^{\varepsilon}_{s,d}(\boldsymbol{\hat{\eta}})&:=&\sup\limits_{\boldsymbol{\eta} \in \mathcal{H}^{\beta}_{s,d}} \sup\limits_{\thetab\in  \Theta_{s,d}^{\sigma}(r_{\varepsilon})} \E{\thetab,\etab} \left\{\sum_{k=1}^s{d\choose k}^{\beta-1}  |\boldsymbol{\hat{\eta}}_{k} - {\etab}_{k}|\right\} \nonumber \\
	&=&\sup\limits_{\boldsymbol{\eta} \in \mathcal{H}_{s,d}^{\beta}} \sup\limits_{\thetab\in  \Theta_{s,d}^{\sigma}(r_{\varepsilon})}
	\sum_{k=1}^s{d\choose k}^{\beta-1}{\operatorname{E}}_{\thetab_k,\etab_k} |\boldsymbol{\hat{\eta}}_k - {\etab_k}|, \label{def:Ham_risk_s}
\end{eqnarray}
where  $r_\varepsilon=\{r_{\varepsilon,k},$ $1\leq k\leq s\}$, $r_{\varepsilon,k}>0$, 
is a positive family of collections,  ${\cal H}_{s,d}^{\beta}$ is the set defined in (\ref{cH}), and
\begin{gather*}
	\Theta_{s,d}^{\sigma} (r_{\varepsilon})=\left \{\thetab=(\thetab_1,\ldots,\thetab_s): \thetab_k
	\in
	\mathring{\Theta}^{\sigma}_{k,d}(r_{\varepsilon,k}),\,1\leq k\leq s\right \}.
\end{gather*}
(For brevity, the dependence of the maximum normalized Hamming risk 
${\cal R}^{\varepsilon}_{s,d}(\boldsymbol{\hat{\eta}})$ on $\sigma$ and $\beta$ is suppressed.)
An \textit{almost full aggregate  selector} $\boldsymbol{\hat{\eta}}$ is defined for all $\beta\in (0,1)$ and $\sigma>0$, such that
$$\limsup_{\varepsilon\to 0}{\cal R}_{s,d}^{\varepsilon}(\boldsymbol{\hat{\eta}})=0.$$
If such a selector exists, we say that almost full selection in model (\ref{model2_s}) is possible.
The almost full selection of the relevant $\eta_{u_k}$ in the sequence space model (\ref{model2_s}) is possible when all elements of the family of collections 
$r_\varepsilon=\{r_{\varepsilon,k},$ $1\leq k\leq s\}$ defining the set $ \Theta^{\sigma}_{s,d} (r_{\varepsilon})$ are not ``too small'', in which case, 
the hypotheses ${\bf H}_{0,u_k}$ and ${\bf H}_{1,u_k}^{\varepsilon}$ in 
(\ref{test1}) separate asymptotically for $u_k\in\Ukd$, $1\leq k\leq s$. In this section, we derive the \textit{sharp almost full selection boundary}, stated in terms of
$r_\varepsilon=\{r_{\varepsilon,k},$ $1\leq k\leq s\}$, which defines when almost full selection of the relevant index variables $\eta_{u_k}$ in model (\ref{model2_s})
is possible and when this type of selection is impossible. The \textit{sharp exact selection boundary} associated with model (\ref{model2_s})
was obtained in Section 3 of \cite{ST-2024}.

In order to state Theorems \ref{th4} and \ref{th5}, take $\epsilon=\epsilon_{\varepsilon,s,d}>0$ to be such that
\begin{gather}\label{eps}
	\epsilon\to 0,\quad\epsilon\log d\to \infty,\quad \epsilon\log{d\choose s}=o(\log \varepsilon^{-1}),\quad \mbox{as}\quad \varepsilon\to 0,
\end{gather}
and consider the adaptive aggregate selector $\boldsymbol{\hat{\eta}}=(\boldsymbol{\hat{\eta}}_1,\ldots,\boldsymbol{\hat{\eta}}_s)$, where $\boldsymbol{\hat{\eta}}_k$, $1\leq k\leq s$, is
given by (\ref{def:selector_for_beta_unknown})--(\ref{def:tau}) with $\epsilon$ as in (\ref{eps})  instead of (\ref{cond_epsilon}).
We claim that the selector $ \boldsymbol{\hat{\eta}}=(\boldsymbol{\hat{\eta}}_1,\ldots,\boldsymbol{\hat{\eta}}_s)$
achieves almost full selection, and thus extends Theorem \ref{th:sim_UB_fixedK} to a more general sparse model.
The precise statement is as follows.

\begin{theorem} \label{th4}
Let $\beta \in [b,B]\subset(0,1)$, $\sigma >0$, and $s \in \{1,\ldots,d\}$ be fixed numbers,
and let $d=d_\varepsilon \to \infty$ and $\varepsilon^2(\log d)^{1+2\sigma}=o(1)$ as $\varepsilon\to 0$.
Assume that the family of collections $r_\varepsilon=\{r_{\varepsilon,k},1\leq k\leq s\}$, $r_{\varepsilon,k}>0$, is such that
\begin{gather}\label{nomer0}
	\liminf_{\varepsilon \to 0} \min_{1\leq k\leq s}\frac{a_{\varepsilon,{u_k}}(r_{\varepsilon,k})}{\sqrt{2\log {d\choose k}}} > \sqrt{\beta}.
\end{gather}
Then
\begin{equation*}
	\limsup_{\varepsilon\to 0} \mathcal{R}^{\varepsilon}_{s,d}(\boldsymbol{\hat{\eta}})=\limsup_{\varepsilon\to 0} \sup\limits_{\boldsymbol{\eta} \in \mathcal{H}_{s,d}^{\beta}} \sup\limits_{\thetab\in  \Theta_{s,d}^{\sigma}(r_{\varepsilon})}  {\rm E}_{\thetab,\etab} \left\{\sum_{k=1}^s {d \choose k}^{\beta-1}|\boldsymbol{\hat{\eta}}_k-\boldsymbol{\eta}_k| \right\}=0,
\end{equation*}
where $ \boldsymbol{\hat{\eta}}=(\boldsymbol{\hat{\eta}}_1,\ldots,\boldsymbol{\hat{\eta}}_s)$ is the aggregate selector with components defined in {\rm (\ref{def:selector_for_beta_unknown})--(\ref{def:tau})} and with $\epsilon$ as in {\rm (\ref{eps})}.
\end{theorem}

The analogue of Theorem \ref{th:sim_LB_fixedK} for the  general model at hand is as follows.

\begin{theorem} \label{th5}
Let $\beta \in (0,1)$, $\sigma >0$, and $s \in \{1,\ldots,d\}$ be fixed numbers,
and let $d=d_\varepsilon \to \infty$ and $\varepsilon^2(\log d)^{1+2\sigma}=o(1)$ as $\varepsilon\to 0$.
Assume that the family of collections $r_\varepsilon=\{r_{\varepsilon,k},1\leq k\leq s\}$, $r_{\varepsilon,k}>0$, is such that
\begin{gather}\label{nomer1}
	\limsup_{\varepsilon \to 0} \min_{1\leq k\leq s}\frac{a_{\varepsilon,{u_k}}(r_{\varepsilon,k})}{\sqrt{2\log {d\choose k}}} < \sqrt{\beta}.
\end{gather}
Then
\begin{equation*}
	\liminf_{\varepsilon\to 0}\inf_{\boldsymbol{\tilde{\eta}}} \mathcal{R}^{\varepsilon}_{s,d}(\boldsymbol{\tilde{\eta}}) =	\liminf_{\varepsilon\to 0}\inf_{\boldsymbol{\tilde{\eta}}} \sup\limits_{\boldsymbol{\eta} \in \mathcal{H}_{s,d}^{\beta}}
	\sup\limits_{\thetab\in  \Theta_{s,d}^{\sigma}(r_{\varepsilon})}  {\rm E}_{\thetab,\etab}
	\left\{\sum_{k=1}^s {d \choose k}^{\beta-1}|\boldsymbol{\tilde{\eta}}_k-\boldsymbol{\eta}_k| \right\}>0,
\end{equation*}
where the infimum is taken over all aggregate selectors $\boldsymbol{\tilde{\eta}}=(\boldsymbol{\tilde{\eta}}_1,\ldots,\boldsymbol{\tilde{\eta}}_s)$ of 
$\etab=(\etab_1,\ldots,\etab_s)\in \mathcal{H}_{s,d}^{\beta}$ in model {\rm (\ref{model2_s})}.
\end{theorem}

Theorems \ref{th4} and \ref{th5} imply that, in the problem of identifying nonzero $\eta_{u_k}$s in model (\ref{model2_s}) almost fully, the aggregate selector
$ \boldsymbol{\hat{\eta}}=(\boldsymbol{\hat{\eta}}_1,\ldots,\boldsymbol{\hat{\eta}}_s)$ of $\etab=(\etab_1,\ldots,\etab_s)\in  \mathcal{H}_{s,d}^{\beta}$ given by (\ref{def:selector_for_beta_unknown})--\rm{(\ref{def:tau})} and (\ref{eps}) is \textit{optimal} in the asymptotically minimax sense.

We now state the analogues of Theorems \ref{th4} and \ref{th5} for the case when $s=s_\varepsilon\to \infty$, $s=o(d)$ as $\varepsilon\to 0$.
To do this, we need to slightly modify  the selector $ \boldsymbol{\hat{\eta}}=(\boldsymbol{\hat{\eta}}_1,\ldots,\boldsymbol{\hat{\eta}}_s)$ given by {\rm (\ref{def:selector_for_beta_unknown})--(\ref{def:tau})} and {\rm (\ref{eps})}. 
Specifically, in the definition of  $ \boldsymbol{\hat{\eta}}$, we need to replace 
condition (\ref{def:tau})  by  condition (\ref{def:tau_s}), and 
condition {\rm (\ref{eps})} by
the condition
\begin{gather}\label{eps1}
	\epsilon\to 0,\quad \epsilon\log (d/s)\to \infty,
	\quad \epsilon\log{d\choose s}=o(\log \varepsilon^{-1}),\quad\mbox{as}\;\; \varepsilon\to 0.
\end{gather}
The following results hold true.

\begin{theorem} \label{th6}
Let $\beta \in [b,B]\subset(0,1)$ and $\sigma >0$ be fixed numbers,
and let $d=d_\varepsilon \to \infty$ and $s=s_\varepsilon\to \infty$ be such that
$s=o(d)$, $s=o(\log \varepsilon^{-1})$, $\log \log d=o(s)$, and $\varepsilon^2(\log d)^{1+2\sigma}=o(1)$,
as $\varepsilon\to 0$.
Assume that the family of collections $r_\varepsilon=\{r_{\varepsilon,k},1\leq k\leq s\}$, $r_{\varepsilon,k}>0$, is as in Theorem {\rm \ref{th4}}.
Then
\begin{equation*}
	\limsup_{\varepsilon\to 0} \mathcal{R}^{\varepsilon}_{s,d}(\boldsymbol{\hat{\eta}})=\limsup_{\varepsilon\to 0} \sup\limits_{\boldsymbol{\eta} \in \mathcal{H}_{s,d}^{\beta}} \sup\limits_{\thetab\in  \Theta_{s,d}^{\sigma}(r_{\varepsilon})}  {\rm E}_{\thetab,\etab} \left\{\sum_{k=1}^s {d \choose k}^{\beta-1}|\boldsymbol{\hat{\eta}}_k-\boldsymbol{\eta}_k| \right\}=0,
\end{equation*}
where $ \boldsymbol{\hat{\eta}}=(\boldsymbol{\hat{\eta}}_1,\ldots,\boldsymbol{\hat{\eta}}_s)$ is  the aggregate selector with components defined 
in {\rm (\ref{def:selector_for_beta_unknown})--(\ref{vj}) and (\ref{def:tau_s})} and with $\epsilon$ as in {\rm (\ref{eps1})}.
\end{theorem}

\begin{theorem} \label{th7}
Let $\beta \in (0,1)$ and $\sigma >0$ be fixed numbers,
and let $d=d_\varepsilon \to \infty$ and $s=s_\varepsilon\to \infty$ be such that
$s=o(d)$, $s=o(\log \varepsilon^{-1})$, $\log \log d=o(s)$, and $\varepsilon^2(\log d)^{1+2\sigma}=o(1)$,
as $\varepsilon\to 0$.
Assume that the family of collections $r_\varepsilon=\{r_{\varepsilon,k},1\leq k\leq s\}$, $r_{\varepsilon,k}>0$, is as in Theorem {\rm \ref{th5}}.
Then
\begin{equation*}
	\liminf_{\varepsilon\to 0}\inf_{\boldsymbol{\tilde{\eta}}} \mathcal{R}^{\varepsilon}_{s,d}(\boldsymbol{\tilde{\eta}}) =	\liminf_{\varepsilon\to 0}\inf_{\boldsymbol{\tilde{\eta}}} \sup\limits_{\boldsymbol{\eta} \in \mathcal{H}_{s,d}^{\beta}}
	\sup\limits_{\thetab\in  \Theta_{s,d}^{\sigma}(r_{\varepsilon})}  {\rm E}_{\thetab,\etab}
	\left\{\sum_{k=1}^s {d \choose k}^{\beta-1}|\boldsymbol{\tilde{\eta}}_k-\boldsymbol{\eta}_k| \right\}>0,
\end{equation*}
where the infimum is taken over all aggregate selectors $\boldsymbol{\tilde{\eta}}=(\boldsymbol{\tilde{\eta}}_1,\ldots,\boldsymbol{\tilde{\eta}}_s)$ of 
$\etab=(\etab_1,\ldots,\etab_s)\in \mathcal{H}_{s,d}^{\beta}$ in model {\rm (\ref{model2_s})}.
\end{theorem}

Inequalities {\rm (\ref{nomer0})} and  {\rm (\ref{nomer1})} describe the \textit{sharp almost full selection boundary}
for the general sequence space model (\ref{model2_s}).
This boundary is determined in terms of the function $a_{\varepsilon,u_k}(r_{\varepsilon,k})$ defined in {\rm (\ref{def:a2})}
whose sharp asymptotics for fixed $k$ are given by (\ref{aek}) and, when
$k=k_\varepsilon\to \infty$, cf. formula (39) in \cite{IS-2005},
\begin{gather*} 
a_{\varepsilon,u_k}(r_{\varepsilon,k})\sim \left({2\pi k}/{e}\right)^{{k}/{4}}e^{-1}(\pi k)^{1/4}r_{\varepsilon,k}^{2+{k}/{(2\sigma)}}\varepsilon^{-2}.
\end{gather*}
Theorems \ref{th4} to \ref{th7} augment Theorems 1 to 4 of \cite{ST-2024}, which were established in the context of exact
variable selection in the Gaussian sequence space model (\ref{model2_s}). In the latter theorems, the \textit{sharp exact selection boundary} is found to be given by
the inequalities, cf.
(\ref{nomer0}) and (\ref{nomer1}), 
\begin{equation*}\label{testing5}
\liminf_{\varepsilon \to 0} \min_{1\leq k\leq s} \frac{a_{\varepsilon,u_k}(r_{\varepsilon,k})}{\sqrt{2\log {d\choose k}}} > 1+\sqrt{1-\beta}
\quad \mbox{and}\quad  \limsup_{\varepsilon \to 0} \min_{1\leq k\leq s} \frac{a_{\varepsilon,u_k}(r_{\varepsilon,k})}{\sqrt{2\log {d\choose k}}} < 1+\sqrt{1-\beta}
\end{equation*}
for both cases (i) when $s$ is fixed and (ii) when $s=s_\varepsilon\to \infty$, $s=o(d)$ as $\varepsilon\to 0$.

\section{Proofs of Theorems}
Before proving Theorems \ref{th:sim_UB_fixedK_knownBeta} to \ref{th7}, we recall some 
properties of the statistics $S_{u_k}(\beta)$ defined in (\ref{def:S_u}).
These properties can be found, for example, in Section 7.1 of \cite{DSA-2012} and Section 5 of \cite{ST-2023}.

Consider testing ${\bf H}_{0,u_k}$  versus $ {\bf H}^{\varepsilon}_{1,u_k}$ as in (\ref{test1}),
and let the expectation and variance under ${\bf H}_{0,u_k}$ be denoted by ${\rm E}_{\bf 0}$ and ${\rm var}_{\bf 0}$, and  under $ {\bf H}^{\varepsilon}_{1,u_k}$
by ${\rm E}_{\thetab_{u_k}}$ and ${\rm var}_{\thetab_{u_k}}$. Clearly, ${\rm E}_{\bf 0}(S_{u_k}(\beta))=0$ and ${\rm var}_{\bf 0} (S_{u_k}(\beta))=1$ for all $\beta\in(0,1)$.
Next, for all $\thetab_{u_k}\in \mathring{\Theta}_{\cb_{u_k}}(r_{\varepsilon,k})$ and  $\beta\in(0,1)$, 
\begin{align} 
{\rm E}_{\thetab_{u_k}}(S_{u_k}(\beta))&=\sum_{\lb\in\mathring{\mathbb{Z}}_{u_k}}\omega_\lb (r^*_{\varepsilon,k}) (\theta_\lb(u_k)/\varepsilon)^2, \label{Es_knownBeta} \\
{\rm var}_{\thetab_{u_k}} (S_{u_k}(\beta))& = \sum_{\lb \in \mathring{\mathbb{Z}}_{u_k}} \omega^2_\lb (r_{\varepsilon,k}^*) \left(2+4 (\theta_{\lb}(u_k)/\varepsilon)^2\right) = 
1+4\sum_{\lb \in \mathring{\mathbb{Z}}_{u_k}} \omega^2_\lb (r_{\varepsilon,k}^*)(\theta_\lb(u_k)/\varepsilon)^2\nonumber \\
& = 
1 + O \left( \max_{\lb\in  \mathring{\mathbb{Z}}_{u_k}}\omega_\lb(r_{\varepsilon,k}^*)\E{\thetab_{u_k}}( S_{u_k}(\beta)) \right), \quad \varepsilon \to 0.
\label{VarS_knownBeta}
\end{align}
Moreover,  if  $T_k =T_{k,\varepsilon}\rightarrow \infty$ as $\varepsilon\to 0$ is such that
\begin{gather}
T_k \max_{\lb \in \mathring{\mathbb{Z}}_{u_k}} \omega_\lb (r^*_{\varepsilon,k})=o(1), \quad \varepsilon\to 0,
\label{eq:exp_bound_condT}
\end{gather}
then 
for all $\beta\in(0,1)$ 
as $\varepsilon \to 0$ 
\begin{gather}\label{eq:exp_bound}
\operatorname{P}_{\bf 0}\left({ S_{u_k}(\beta) > T_k}\right) \leq \exp \left( - \dfrac{T_k^2}{2}(1+o(1)) \right).
\end{gather}
If, in addition to (\ref{eq:exp_bound_condT}), it holds for $\thetab_k \in \mathring{\Theta}_{c_k}(r_{\varepsilon,k})$ 
and some $\beta\in(0,1)$ that 
\begin{gather*}
{\rm E}_{\thetab_{u_k}} (S_{u_k}(\beta)) \max_{\lb \in \mathring{\mathbb{Z}}_{u_k}} \omega_\lb (r^*_{\varepsilon,k}) =o(1), \quad \varepsilon\to 0,
\end{gather*}
then, for these $\thetab_k$ and all  $\beta\in(0,1)$, 
as $\varepsilon \to 0$  
\begin{gather}\label{eq:exp_bound2}
\operatorname{P}_{\thetab_{u_k}}\left({ S_{u_k}(\beta)- {\rm E}_{\thetab_{u_k}} (S_{u_k}(\beta))\leq - T_k}\right) \leq \exp \left( - \dfrac{T_k^2}{2}(1+o(1)) \right).
\end{gather}
The proof of (\ref{eq:exp_bound}) and (\ref{eq:exp_bound2}) is similar to that of Proposition 7.1 in \cite{DSA-2012} and therefore is omitted.

\medskip
\textit{Proof of Theorem \ref{th:sim_UB_fixedK_knownBeta}.}
For a given~$k$ ($1\leq k\leq d$), the maximum normalized Hamming risk of $\boldsymbol{\hat{\eta}}_k(\beta)$ satisfies
\begin{align*}
& R^{\varepsilon}_{k,d}(\boldsymbol{\hat{\eta}}_k(\beta))  :=	\sup_{\etab_k\in H^{\beta}_{k,d}}\sup_{\thetab_k\in  \mathring{\Theta}^\sigma_{k,d}(r_{\varepsilon,k})}
{d \choose k}^{\beta-1}{\rm E}_{\thetab_k,\etab_k} |\boldsymbol{\hat{\eta}}_k(\beta)-\boldsymbol{\eta}_k| \nonumber \\
&= \sup_{\etab_k\in H^{\beta}_{k,d}}\sup_{\thetab_k\in  \mathring{\Theta}^\sigma_{k,d}(r_{\varepsilon,k})} {d \choose k}^{\beta-1}  \Bigg( \sum_{\substack{u_k : \eta_{u_k}=0}} {\rm E}_{\boldsymbol{0}}  (\hat{\eta}_{u_k}(\beta)) + \sum_{\substack{u_k : \eta_{u_k}=1}} {\rm E}_{\thetab_{u_k}}  (1-\hat{\eta}_{u_k}(\beta))  \Bigg) 
\end{align*}
\begin{align*} 
\begin{split}
	& \leq \sup_{\etab_k\in H^{\beta}_{k,d}} {d \choose k}^{\beta-1} \sum_{\substack{u_k : \eta_{u_k}=0}}  \Prob{\boldsymbol{0}}{S_{u_k}(\beta) > \sqrt{(2\beta+\epsilon) \log {d \choose k}}} + \\
	& \quad + \sup_{\etab_k\in H^{\beta}_{k,d}}\sup_{\thetab_k\in  \mathring{\Theta}^\sigma_{k,d}(r_{\varepsilon,k})} 
	{d \choose k}^{\beta-1} \sum_{\substack{u_k : \eta_{u_k}=1}} \Prob{\thetab_{u_k}}{S_{u_k}(\beta) \leq \sqrt{(2\beta+\epsilon) \log {d \choose k}}}.
\end{split}
\end{align*}
From this, using the sparsity condition (\ref{sparsitycond})  and the fact that $\#(\Ukd)={d\choose k}$, 
we have for all small enough $\varepsilon$ 
\begin{equation} 
\begin{split}\label{eq:I1_I2_knownBeta}
	R^{\varepsilon}_{k,d}(\boldsymbol{\hat{\eta}}_k(\beta))  &\leq {d \choose k}^{\beta} \Prob{\boldsymbol{0}}{S_{u_k}(\beta) > \sqrt{(2\beta+\epsilon) \log {d \choose k}}} +  \\
	& \qquad +  C \sup_{\thetab_{u_k} \in \mathring{\Theta}_{c_{u_k}}(r_{\varepsilon,k}) }    \Prob{\thetab_{u_k}}{S_{u_k}(\beta) \leq \sqrt{(2\beta+\epsilon) \log {d \choose k}}}  \\
	&=: I^{(1)}_{\varepsilon,k} +  I^{(2)}_{\varepsilon,k},
\end{split}
\end{equation}
where $u_k$ on the right-hand side of (\ref{eq:I1_I2_knownBeta}) is  any (selected)  element of $\Ukd$.
Due to (\ref{WL}), (\ref{rr}), and (\ref{th:assump_epsilon}), we have for $T_k= \sqrt{(2\beta+\epsilon) \log {d \choose k}}$ that
\begin{equation*} 
T_k \max_{\lb \in \mathring{\mathbb{Z}}_{u_k}} \omega_\lb (r^*_{\varepsilon,k})\asymp \log^{1/2}{d \choose k}  \left( \varepsilon \log^{1/4} {d\choose k} \right)^{{2k}/{(4 \sigma +k)}}= o(1), \quad \varepsilon\to 0. 
\end{equation*}
Therefore, condition (\ref{eq:exp_bound_condT}) is fulfilled and the application of 
(\ref{eq:exp_bound}), together with (\ref{cond_epsilon}), gives
\begin{gather}\label{eq:bound_I1}
I_{\varepsilon,k}^{(1)} \leq {d \choose k}^{\beta} \exp \left\{- \left( \beta +{\epsilon}/{2}\right)  \log {d \choose k} (1 + o(1))\right\} 
= {O} \left({d \choose k}^{- \epsilon/2} \right) = o(1),\quad \varepsilon\to 0.
\end{gather}

Next, consider the term $I_{\varepsilon,k}^{(2)}$ on the right-hand side of (\ref{eq:I1_I2_knownBeta}).
In view of condition (\ref{cond:inf}), there exists a~constant $\Delta_{0,k} \in (0,1)$ such that for all small enough $\varepsilon$ 
\begin{equation*} 
\frac{a_{\varepsilon,u_k}(r_{\varepsilon,k})}{\sqrt{\log {d\choose k}}} > \sqrt{2\beta} (1+\Delta_{0,k}).
\end{equation*}
Therefore, by the definition of $r_{\varepsilon,k}^*$ in (\ref{def:r_star}) and using the ``continuity'' property of $a_{\varepsilon,u_k}$ as in~(\ref{def:continuity_of_a}),
for all small enough $\varepsilon$ and some (small) constants $\Delta_{1,k} \in (0,1)$ and $\Delta_{2,k} \in (0,1)$, we can write
\begin{equation}\label{bound_for_a_by_a}
a_{\varepsilon,u_k} \big((1+\Delta_{1,k}) r^*_{\varepsilon,k}\big) \leq (1+\Delta_{2,k}) a_{\varepsilon,u_k}(r^*_{\varepsilon,k}) = (1+\Delta_{2,k}) \sqrt{2 \beta \log {d\choose k}} \leq a_{\varepsilon,u_k}(r_{\varepsilon,k}).
\end{equation}  
From this, by the monotonicity of $a_{\varepsilon,u_k}$,
\begin{equation} \label{r_as_B_r_star_knownBeta}
r_{\varepsilon,k} \geq (1+\Delta_{1,k}) r^*_{\varepsilon,k} =: B_k r^*_{\varepsilon,k}, \quad B_k >1.
\end{equation}
Since for any $\mathbb{B}_k\geq 1,$ $\varepsilon>0$, $r_{\varepsilon,k}>0$, $1\leq k\leq d$, one has, cf. relation (52) in \cite{ST-2023},
\begin{equation}\label{lemma:ball}
{\varepsilon^{-2}} \inf_{\thetab_{u_k} \in \mathring{\Theta}_{\cb_{u_k}}(\mathbb{B}_k r_{\varepsilon,k})} \sum_{\lb \in \mathring{\mathbb{Z}}_{u_k}} \omega_\lb (r_{\varepsilon,k}) (\theta_\lb(u_k))^2 \geq \mathbb{B}_k^2 a_{\varepsilon,u_k}(r_{\varepsilon,k}),
\end{equation}  
it now follows that for $\beta\in(0,1)$ and all small enough $\varepsilon$
\begin{align} 	\label{bound_infES_knownBeta}
\inf_{\thetab_{u_k} \in \mathring{\Theta}_{\cb_{u_k}}(r_{\varepsilon,k})} & {\rm E}_{\thetab_{u_k}}(S_{u_k}(\beta)) \geq\varepsilon^{-2}
\inf_{\thetab_{u_k} \in \mathring{\Theta}_{\cb_{u_k}}(B_kr_{\varepsilon,k}^*)} \sum_{\lb \in \mathring{\mathbb{Z}}_{u_k}} \omega_\lb(r_{\varepsilon,k}^*) 
(\theta_\lb(u_k))^2 \nonumber  \\ 
&\geq B_k^2 a_{\varepsilon,k} (r_{\varepsilon,k}^*) 
= B_k^2 \sqrt{2\beta \log {d \choose k}} > \sqrt{(2 \beta+\epsilon)\log {d\choose k}},
\end{align}
and, in particular, 
\begin{equation*} 
\sqrt{(2 \beta+\epsilon)\log {d\choose k}} -	\inf_{\thetab_{u_k} \in \mathring{\Theta}_{\cb_{u_k}}(r_{\varepsilon,k})}  {\rm E}_{\thetab_{u_k}}(S_{u_k}(\beta)) \to -\infty, \quad \varepsilon \to 0.
\end{equation*}
Therefore, uniformly in $\thetab_{u_k} \in \mathring{\Theta}_{\cb_{u_k}}(r_{\varepsilon,k}),$ $u_k\in \Ukd$, $1\leq k\leq d$, for $\beta\in(0,1)$
\begin{equation} \label{sqrt_minus_ES_to_infty_knownBeta}
\sqrt{(2 \beta+\epsilon)\log {d\choose k}} -	 {\rm E}_{\thetab_{u_k}}(S_{u_k}(\beta)) \to -\infty, \quad \varepsilon \to 0.
\end{equation}
Moreover, in view of (\ref{bound_infES_knownBeta}), for $\beta\in(0,1)$ and all small enough $\varepsilon$,
\begin{equation} \label{infES_minus_sqrt_LB_knownBeta}
\inf_{\thetab_{u_k} \in \mathring{\Theta}_{\cb_{u_k}}(r_{\varepsilon,k})} {\rm E}_{\thetab_{u_k}}(S_{u_k}(\beta)) - \sqrt{(2 \beta+\epsilon)\log {d\choose k}} \geq \sqrt{2\beta \log {d \choose k}} \left(B_k^2 -1 + o(1) \right).
\end{equation}

Now, for the element $u_k$ of $\Ukd$ selected in (\ref{eq:I1_I2_knownBeta}) and some $\beta\in(0,1)$, consider the subsets $\mathring{\Theta}^{(p)}_{\cb_{u_k}}(r_{\varepsilon,k})$, $p=1,2$, of
$\mathring{\Theta}_{\cb_{u_k}}(r_{\varepsilon,k})$ defined by
\begin{equation}
\begin{split}
	\mathring{\Theta}^{(1)}_{\cb_{u_k}}(r_{\varepsilon,k})&= \left\{\thetab_{u_k}\in \mathring{\Theta}_{\cb_{u_k}}(r_{\varepsilon,k}):
	\limsup_{\varepsilon\to 0} {\rm E}_{\thetab_{u_k}}(S_{u_k}(\beta))  \max_{\lb \in \mathring{\mathbb{Z}}_{u_k}} \omega_\lb (r^*_{\varepsilon,k})=0 \right\},\\
	\mathring{\Theta}^{(2)}_{\cb_{u_k}}(r_{\varepsilon,k})&=\left\{\thetab_{u_k}\in \mathring{\Theta}_{\cb_{u_k}}(r_{\varepsilon,k}):a\leq \liminf_{\varepsilon\to 0}
	{\rm E}_{\thetab_{u_k}}(S_{u_k}(\beta)) \max_{\lb \in \mathring{\mathbb{Z}}_{u_k}} \omega_\lb (r^*_{\varepsilon,k})\leq\right. \\
	\leq & \left. \limsup_{\varepsilon\to 0} {\rm E}_{\thetab_{u_k}}(S_{u_k}(\beta))  \max_{\lb \in \mathring{\mathbb{Z}}_{u_k}} \omega_\lb (r^*_{\varepsilon,k})
	\leq A\; \mbox{for some}\; 0<a\leq A\leq\infty\right\}.
\end{split}
\label{def:Theta123_knownBeta}
\end{equation}
Note that $\mathring{\Theta}_{\cb_{u_k}}(r_{\varepsilon,k})\subseteq \mathring{\Theta}^{(1)}_{\cb_{u_k}}(r_{\varepsilon,k})\cup  \mathring{\Theta}^{(2)}_{\cb_{u_k}}(r_{\varepsilon,k})$ and, by means of (\ref{omegamax_asymp}) and (\ref{th:assump_epsilon}), 
for all $\thetab_{u_k}\in\mathring{\Theta}^{(2)}_{\cb_{u_k}}(r_{\varepsilon,k})$ and $\beta\in(0,1)$ one has 
\begin{equation} \label{eq:log_dk_is_oE}
{\log^{1/2} {d\choose k}}=o\left( {\rm E}_{\thetab_{u_k}}(S_{u_k}(\beta))\right), \quad \text{as}\; \varepsilon\to 0. 	
\end{equation}
Due to (\ref{WL}), (\ref{rr}), and 
(\ref{Es_knownBeta}), the rate at which $ {\rm E}_{\thetab_{u_k}}(S_{u_k}(\beta))$  varies, as $\varepsilon$ tends to zero, is independent of $\beta$.
Hence, the choice of $\beta\in (0,1)$ in the definition of the subsets in (\ref{def:Theta123_knownBeta}) is irrelevant. 
Also, by means of (\ref{WL}) and (\ref{rr}), for all  $m=1,\ldots, M_k$, 
$\max_{\lb \in \mathring{\mathbb{Z}}_{u_k}} \omega_\lb (r^*_{\varepsilon,k})\asymp\max_{\lb \in \mathring{\mathbb{Z}}_{u_k}} \omega_\lb (r^*_{\varepsilon,k,m})$,
and therefore replacing $r^*_{\varepsilon,k}$ by 
$r^*_{\varepsilon,k,m}$ in (\ref{def:Theta123_knownBeta}) 
does not change the sets $\mathring{\Theta}^{(p)}_{\cb_{u_k}}(r_{\varepsilon,k})$, $p=1,2$.
This observation will be used in the proof of Theorem \ref{th6}.

Thanks to (\ref{sqrt_minus_ES_to_infty_knownBeta}) and (\ref{infES_minus_sqrt_LB_knownBeta}), the application of Chebyshev's inequality gives as $\varepsilon \to 0$
\begin{align}	\label{eq:I2_v1}
{I_{\varepsilon,k}^{(2)}}  &\leq  C \sum_{p=1}^2\sup_{\thetab_{u_k} \in \mathring{\Theta}^{(p)}_{c_{u_k}}(r_{\varepsilon,k}) }   \operatorname{P}_{\thetab_{u_k}}\Bigg(S_{u_k}(\beta)-
{\rm E}_{\thetab_{u_k}}(S_{u_k}(\beta)) \leq \nonumber\\ 
&\hspace{50mm}   \leq  \sqrt{(2 \beta + \epsilon) \log{d \choose k}}-{\rm E}_{\thetab_{u_k}}\left(S_{u_k}(\beta)\right)\Bigg) \nonumber  \\
\begin{split}
&\leq C	\sup_{\thetab_{u_k}\in \mathring{\Theta}^{(1)}_{\cb_{u_k}}(r_{\varepsilon,k})}
\operatorname{P}_{\thetab_{u_k}}\Bigg( S_{u_k}(\beta)- {\rm E}_{\thetab_{u_k}}(S_{u_k}(\beta))
\leq -\sqrt{2 \beta  \log{d \choose k}}\left(B_k^2-1+o(1)\right) \Bigg)	\\
&	\qquad \qquad +	C  \sup_{\thetab_{u_k}\in \mathring{\Theta}^{(2)}_{\cb_k}(r_{\varepsilon,k})      }\frac{ {\rm var}_{\thetab_{u_k}}(S_{u_k}(\beta))}{\left( {\rm E}_{\thetab_{u_k}}(S_{u_k}(\beta))  - \sqrt{(2\beta+\epsilon) \log {d\choose k}} \right)^2 }.
\end{split}
\end{align}  
From this and inequality (\ref{eq:exp_bound2}), taking into account 
(\ref{WL}), (\ref{rr}), (\ref{VarS_knownBeta}), and (\ref{eq:log_dk_is_oE}), 
for all small enough $\varepsilon$ and some positive constants $C_1$ and $C_2$, the same for all $1\leq k\leq d$,  
we can write
\begin{align}	
{I_{\varepsilon,k}^{(2)}}  
&\leq C \exp\left(- \beta\log {d\choose k}\left(B_k^2-1+o(1)\right)^2(1+o(1)) \right)\nonumber\\ &
\quad \quad + C
\sup_{\thetab_{u_k}\in  \mathring{\Theta}^{(2)}_{\cb_{u_k}}(r_{\varepsilon,k})      }\frac{  1+O \left( \max_{\lb\in  \mathring{\mathbb{Z}}_{u_k}}\omega_\lb(r_{\varepsilon,k}^*){\rm E}_{\thetab_{u_k}}(S_{u_k}(\beta)) \right) }{\left( {\rm E}_{\thetab_{u_k}}(S_{u_k}(\beta))  - \sqrt{(2\beta+\epsilon) \log {d\choose k}} \right)^2 }\nonumber\\
& \leq  C \exp\left(- \frac{\beta}{2}\left(B_k^2-1\right)^2 \log {d\choose k}\right)+
\frac{C_1\max_{\lb\in  \mathring{\mathbb{Z}}_{u_k}}\omega_\lb(r_{\varepsilon,k}^*)}{	\inf_{\thetab_{u_k} \in \mathring{\Theta}^{(2)}_{\cb_{u_k}}(r_{\varepsilon,k})} {\rm E}_{\thetab_{u_k}}(S_{u_k}(\beta)) }\nonumber	\\ 
&\leq C {d\choose k}^{-\beta(B_k^2-1)^2/2}+ C_2\log^{-1/2}  {d\choose k}\left( \varepsilon \log^{1/4} {d\choose k} \right)^{{2k}/{(4\sigma+k)}} =o(1),
\quad \varepsilon\to 0,
\label{qb_Chebychev_knownBeta}
\end{align}  
where the last equality follows from condition (\ref{th:assump_epsilon}) and the fact that for all $1\leq k\leq d$
\begin{gather*}\label{Bminimum}
\beta(B_k^2-1)^2 \geq \beta \liminf_{d\to \infty}\min_{1\leq k\leq d}(B_k^2-1)^2>0.
\end{gather*}  
Substituting (\ref{eq:bound_I1}) and (\ref{qb_Chebychev_knownBeta}) into (\ref{eq:I1_I2_knownBeta}) gives
\begin{equation*} 
R^{\varepsilon}_{k,d}(\boldsymbol{\hat{\eta}}_k(\beta)) \leq I^{(1)}_{\varepsilon,k} +  I^{(2)}_{\varepsilon,k} = o(1), \quad \varepsilon \to 0.
\end{equation*}
The proof of Theorem \ref{th:sim_UB_fixedK_knownBeta} is complete. 
\qedsymbol

\bigskip

\textit{Proof of Theorem \ref{th:sim_UB_fixedK}.} The proof goes partially along the lines of that of Theorem 3 of \cite{BSt-2017}, as well as that of Theorem \ref{th:sim_UB_fixedK_knownBeta} above, with a random point $\beta_{k,\hat{m}_k}$  in place of $\beta$ ($b\leq \beta\leq B$), selected from the grid of points $\beta_{k,m}$, $1\leq m\leq M_k$, as in (\ref{def:grid_for_beta})
by Lepski's method.
For a~given $k$, let index $m_0=m_{0,k}$ ($1\leq m_0\leq M_{k}-1$)  be such that
\begin{equation*} \label{interval_around_beta}
\beta_{k,m_0} \leq  \beta <\beta_{k,m_0+1}.
\end{equation*}
Then, using the law of total probability for expectations, we can write
\begin{align}
R^{\varepsilon}_{k,d}(\boldsymbol{\hat{\eta}}_k) &:=
\sup_{\etab_k\in H^{\beta}_{k,d}}\sup_{\thetab_k\in  \mathring{\Theta}^\sigma_{k,d}(r_{\varepsilon,k})} {d \choose k}^{\beta-1}{\rm E}_{\thetab_k,\etab_k}|\boldsymbol{\hat{\eta}}_k-\boldsymbol{\eta}_k|\nonumber \\
& \leq \sup_{\etab_k\in H^{\beta}_{k,d}}\sup_{\thetab_k\in  \mathring{\Theta}^\sigma_{k,d}(r_{\varepsilon,k})} {d \choose k}^{\beta-1}{\rm E}_{\thetab_k,\etab_k}\Big(| \boldsymbol{\hat{\eta}}_k-\boldsymbol{\eta}_k | \, \big| \, \hat{m}_k \geq m_0 \Big)
\operatorname{P}_{\thetab_k,\etab_k}\left(\hat{m}_k \geq m_0 \right) \nonumber \\
& \quad \quad  +
\sup_{\etab_k\in H^{\beta}_{k,d}}\sup_{\thetab_k\in \mathring{\Theta}^\sigma_{k,d}(r_{\varepsilon,k})} {d \choose k}^{\beta-1}{\rm E}_{\thetab_k,\etab_k}\Big( |\boldsymbol{\hat{\eta}}_k-\boldsymbol{\eta}_k| \, \big| \, \hat{m}_k < m_0 \Big)
\operatorname{P}_{\thetab_k,\etab_k}\left(\hat{m}_k < m_0 \right)\nonumber \\
& =: Q_{\varepsilon,k}^{(1)} + Q_{\varepsilon,k}^{(2)}. 	 \label{def_Q1_Q2}
\end{align}

To derive a good upper bound for the term $ Q_{\varepsilon,k}^{(1)}$, 
we note that, by the triangle inequality and (\ref{Lepski_m}), when $\hat{m}_k\geq m_0$, for all $\etab_k\in H^{\beta}_{k,d}$ and all $\thetab_k\in   \mathring\Theta^{\sigma}_{k,d}(r_{\varepsilon,k})$,
\begin{align*}
|\boldsymbol{\hat{\eta}}_k - \boldsymbol{\eta}_k| \leq  |\boldsymbol{\hat{\eta}}_k - \boldsymbol{\hat{\eta}}_k (\beta_{k,m_0})|  +
|\boldsymbol{\hat{\eta}}_k(\beta_{k,m_0}) - \boldsymbol{\eta}_k| \leq v_{m_0} +  |\boldsymbol{\hat{\eta}}_k (\beta_{k,m_0}) - \boldsymbol{\eta}_k|,
\end{align*}
where, in view of (\ref{eq:d_choose_k_to_rho_const}) and (\ref{vj}),
$	{d \choose k}^{\beta-1} v_{m_0} = {\tau_{k,d}}^{-1}{{d \choose k}^{\beta-\beta_{k,m_0}}} \leq {\tau_{k,d}}^{-1}{{d \choose k}^{\rho_k}} = O(\tau^{-1}_{k,d}) = o(1),$ as $ \varepsilon \to 0.$
From this, we obtain
\begin{align*}
Q_{\varepsilon,k}^{(1)} &= \sup_{\boldsymbol{\eta}_k \in {H}^{\beta}_{k,d}} \sup_{\thetab_k \in \mathring{\Theta}^\sigma_{k,d}(r_{\varepsilon,k}) }{d \choose k}^{\beta-1}  \E{\thetab_k,\etab_k} \Big( |\boldsymbol{\hat{\eta}}_k - \etab_k| \, \big| \, \hat{m}_k \geq m_0 \Big) \operatorname{P}_{\thetab_k,\etab_k} \left( \hat{m}_k \geq m_0 \right) \nonumber \\
&\leq  \sup_{\boldsymbol{\eta}_k \in {H}^{\beta}_{k,d}} \sup_{\thetab_k \in \mathring{\Theta}^\sigma_{k,d}(r_{\varepsilon,k}) }{d \choose k}^{\beta-1}  \E{\thetab_k,\etab_k}  | \boldsymbol{\hat{\eta}}_k (\beta_{k,m_0}) - \boldsymbol{\eta}_k |  + {\tau_{k,d}}^{-1}{{d \choose k}^{\rho_k}} \nonumber \\
&=  \sup_{\boldsymbol{\eta}_k \in {H}^{\beta}_{k,d}} \sup_{\thetab_k \in \mathring{\Theta}^\sigma_{k,d}(r_{\varepsilon,k}) }{d \choose k}^{\beta-1} \sum_{u_k \in \Ukd}  \E{\thetab_{u_k},\eta_{u_k}}  | \hat{\eta}_{u_k} (\beta_{k,m_0}) - \eta_{u_k} |  + {\tau_{k,d}}^{-1}{{d \choose k}^{\rho_k}},\nonumber
\end{align*}
where  $\operatorname{E}_{\thetab_{u_k},\eta_{u_k}} $ is the expectation which corresponds 
to the measure $\operatorname{P}_{\thetab_{u_k},\eta_{u_k}}$ induced by the observation $\Xb_{\!u_k} = (X_\lb)_{\lb \in \mathring{\mathbb{Z}}_{u_k}}$.
  Next, using the sparsity condition (\ref{sparsitycond}) and the definition of $H^{\beta}_{k,d}$, we obtain
\begin{align*}
Q_{\varepsilon,k}^{(1)} 
&\leq  \sup_{\boldsymbol{\eta}_k \in {H}^{\beta}_{k,d}} {d \choose k}^{\beta-1}
\left\{ \sum_{u_k: \eta_{u_k}=0} \operatorname{P}_{\bf 0}\left({S_{u_k}(\beta_{k,m_0}) > \sqrt{(2 \beta_{k,m_0} + \epsilon) \log{d \choose k}}}\right) \right. \nonumber \\
&	 \quad \left. + \sum_{u_k: \eta_{u_k}=1}\sup_{\thetab_{u_k} \in \mathring{\Theta}_{c_{u_k}}(r_{\varepsilon,k}) } \operatorname{P}_{\thetab_{u_k}}\left(S_{u_k}(\beta_{k,m_0}) \leq  \sqrt{(2 \beta_{k,m_0} + \epsilon) \log{d \choose k}}\right)  \right\} + {\tau_{k,d}}^{-1}{{d \choose k}^{\rho_k}} \nonumber \\
&	\leq  {d \choose k}^{\beta} \operatorname{P}_{\bf 0}\left(S_{u_k}(\beta_{k,m_0}) > \sqrt{(2 \beta_{k,m_0} + \epsilon) \log{d \choose k}}\right) \nonumber \\
& \quad + C  \sup_{\thetab_{u_k} \in \mathring{\Theta}_{c_{u_k}}(r_{\varepsilon,k}) }   \operatorname{P}_{\thetab_{u_k}}\Bigg(S_{u_k}(\beta_{k,m_0}) \leq  \sqrt{(2 \beta_{k,m_0} + \epsilon) \log{d \choose k}}\Bigg)  + {\tau_{k,d}}^{-1}{{d \choose k}^{\rho_k}} 
\end{align*}
\begin{align}
 =:  \, q_{\varepsilon,k}^{(1)}  +  q_{\varepsilon,k}^{(2)} + {\tau_{k,d}}^{-1}{{d \choose k}^{\rho_k}}.
\label{eq:Q1_part1}
\end{align}

A good upper bound on the term $q_{\varepsilon,k}^{(1)}$ 
is obtained by using
the arguments leading to relation~(\ref{eq:bound_I1})
in  the proof of Theorem \ref{th:sim_UB_fixedK_knownBeta}, with $\beta_{k,m_0}$ in place of $\beta$.  Specifically, the application of  (\ref{eq:exp_bound}) with $T_k = \sqrt{(2 \beta_{k,m_0} + \epsilon) \log{d \choose k}}$,
together with condition~(\ref{cond_epsilon}) and the inequality  ${d\choose k}^{\beta-\beta_{k,m_0}}\leq  {d\choose k}^{\rho_k}\leq {\rm const}$, which holds true for all small enough $\varepsilon$, gives 
\begin{align} \label{eq:bound_qa}
q_{\varepsilon,k}^{(1)} &\leq {d \choose k}^{\beta} \exp \left\{- \left( \beta_{k,m_0} +{\epsilon}/{2}\right)  \log {d \choose k} (1 + o(1))\right\} = O \left({d \choose k}^{\beta-\beta_{k,m_0} - \epsilon/2} \right)  \nonumber \\
& = O \left({d \choose k}^{\rho_k - \epsilon/2} \right) = {O} \left({d \choose k}^{- \epsilon/2} \right) = o(1),\quad \varepsilon\to 0. 
\end{align}

The term $q_{\varepsilon,k}^{(2)}$ on the right-hand side of (\ref{eq:Q1_part1}) has the same formal appearance as the term $I_{\varepsilon,k}^{(2)}$ in the proof of Theorem~\ref{th:sim_UB_fixedK_knownBeta}, with $\beta_{k,m_0}$ in place of $\beta$. 
By the definition of $r_{\varepsilon,k,m}^*$ in (\ref{def:r_star_ekm}) and the choice of index $m_0$, we have 
\begin{equation*} \label{bound_for_a_r_star}
a_{\varepsilon,u_k} (r^*_{\varepsilon,k,m_0}) = \sqrt{2 \beta_{k,m_0} \log {d\choose k}} \leq \sqrt{2 \beta \log {d\choose k}}.
\end{equation*}  
Hence, relations (\ref{bound_for_a_by_a})--(\ref{lemma:ball}), originally stated for $r^*_{\varepsilon,k}$, also hold true for $r^*_{\varepsilon,k,m_0}$ (with possibly different  constants $\Delta_{1,k}$, $\Delta_{2,k}$ and $B_k$).
Therefore,  relation (\ref{bound_infES_knownBeta}) continues to hold with $\beta_{k,m_0}$
in place of $\beta$, yielding, uniformly in $\thetab_{u_k} \in \mathring{\Theta}_{\cb_{u_k}}(r_{\varepsilon,k})$ for $u_k\in \Ukd$,
\begin{equation} \label{sqrt_minus_ES_to_infty}
\sqrt{(2 \beta_{k,m_0}+\epsilon)\log {d\choose k}} -	 {\rm E}_{\thetab_{u_k}}(S_{u_k}(\beta_{k,m_0})) \to -\infty, \quad \varepsilon \to 0,
\end{equation}
and
for all small enough $\varepsilon$
\begin{equation} \label{infES_minus_sqrt_LB}
\inf_{\thetab_{u_k} \in \mathring{\Theta}_{\cb_{u_k}}(r_{\varepsilon,k})} {\rm E}_{\thetab_{u_k}}(S_{u_k}(\beta_{k,m_0})) - \sqrt{(2 \beta_{k,m_0}+\epsilon)\log {d\choose k}} \geq \sqrt{2\beta_{k,m_0} \log {d \choose k}} \left(B_k^2 -1 + o(1) \right).
\end{equation}

Using the same arguments as in the derivation of (\ref{eq:I2_v1})--(\ref{qb_Chebychev_knownBeta}), 
but applying
(\ref{sqrt_minus_ES_to_infty})--(\ref{infES_minus_sqrt_LB}) instead of 
(\ref{sqrt_minus_ES_to_infty_knownBeta})--(\ref{infES_minus_sqrt_LB_knownBeta}), we obtain for all small enough $\varepsilon$ and some positive constant $C_3$, the same for all $1\leq k\leq d$,  that as $\varepsilon\to 0$
\begin{align}	
{q_{\varepsilon,k}^{(2)}}  
&\leq C {d\choose k}^{-\beta_{k,m_0}(B_k^2-1)^2/2}+ C_3\log^{-1/2}  {d\choose k}\left( \varepsilon \log^{1/4} {d\choose k} \right)^{{2k}/{(4\sigma+k)}} =o(1),
\label{qb_Chebychev}
\end{align}  
where the last equality follows from condition (\ref{th:assump_epsilon}) and the fact that for all $1\leq k\leq d$ 
\begin{gather}\label{Bmin}
\beta_{k,m_0}(B_k^2-1)^2\geq b(B_k^2-1)^2 \geq b\liminf_{d\to \infty}\min_{1\leq k\leq d}(B_k^2-1)^2=:c_0>0.
\end{gather}  
Substituting (\ref{eq:bound_qa}) and (\ref{qb_Chebychev}) into (\ref{eq:Q1_part1}),
and using ${\tau_{k,d}}^{-1}{{d \choose k}^{\rho_k}}=o(1)$, gives
\begin{equation} \label{eq:Q1_bound}
Q_{\varepsilon,k}^{(1)} \leq q_{\varepsilon,k}^{(1)}  + q_{\varepsilon,k}^{(2)} + {\tau_{k,d}}^{-1}{{d \choose k}^{\rho_k}} = o(1), \quad \varepsilon \to 0.
\end{equation}

We now estimate the term $Q_{\varepsilon,k}^{(2)}$ on the right-hand side of (\ref{def_Q1_Q2}).
Noting that $ |\boldsymbol{\hat{\eta}}_k-\boldsymbol{\eta}_k| = \sum_{u_k \in \Ukd} |\hat{\eta}_{u_k} - \eta_{u_k}| \leq {d\choose k}$, we can write
\begin{align}
Q_{\varepsilon,k}^{(2)} &= \sup_{\boldsymbol{\eta}_k \in {H}^{\beta}_{k,d}} \sup_{\thetab_k \in \mathring\Theta_{k,d}^{\sigma}(r_{\varepsilon,k}) }{d \choose k}^{\beta-1}  \E{\thetab_k,\etab_k} \Big( |\boldsymbol{\hat{\eta}}_k-\boldsymbol{\eta}_k  | \, \big| \, \hat{m}_k < m_0 \Big)\operatorname{P}_{\thetab_k,\etab_k} \left( \hat{m}_k < m_0 \right) \nonumber \\
& \leq \sup_{\boldsymbol{\eta}_k \in {H}^{\beta}_{k,d}} \sup_{\thetab_k \in \mathring\Theta_{k,d}^{\sigma}(r_{\varepsilon,k}) }{d \choose k}^{\beta}   \operatorname{P}_{\thetab_k,\etab_k} \left( \hat{m}_k < m_0 \right),
\label{Q2_part1}
\end{align}
where, in view of (\ref{Lepski_m}), for all $\etab_k \in {H}^{\beta}_{k,d}$ and all $\thetab_k \in \mathring\Theta_{k,d}^{\sigma}(r_{\varepsilon,k})$
\begin{eqnarray}
\operatorname{P}_{\thetab_k,\etab_k} \left( \hat{m}_k < m_0 \right) &=& \sum_{j=1}^{m_0-1} \operatorname{P}_{\thetab_k,\etab_k}
\left( \hat{m}_k =j \right) \nonumber \\
&\leq &\sum_{j=1}^{m_0-1} \operatorname{P}_{\thetab_k,\etab_k} \big( \exists\, i \in \{1,\ldots,j\}: |\boldsymbol{\hat{\eta}}_k (\beta_{k,j+1}) - \boldsymbol{\hat{\eta}}_k(\beta_{k,i})|>v_i \big) \nonumber \\
&\leq &\sum_{j=1}^{m_0-1} \sum_{i=1}^{j} \operatorname{P}_{\thetab_k,\etab_k} \left( \sum_{u_k \in \Ukd} |\hat{\eta}_{u_k} (\beta_{k,j+1}) - \hat{\eta}_{u_k}(\beta_{k,i})|>v_i \right).
\label{P_mk_smaller_m0_v1}
\end{eqnarray}
Next, consider the independent events $A_{u_k}(\beta),$ $ u_k \in \Ukd,$ defined by
\begin{equation*}
A_{u_k}(\beta) = \left\{ S_{u_k}(\beta) \leq \sqrt{(2\beta + \epsilon) \log {d\choose k}} \right\}, \quad u_k \in \Ukd, \quad   \beta \in(0,1),
\end{equation*}
and denote by $\overline{A_{u_k}(\beta)}$ the complement of $A_{u_k}(\beta)$. Observe that for all 
$1\leq i\leq j\leq m_0-1$, 
the equality $|\hat{\eta}_{u_k} (\beta_{k,j+1}) - \hat{\eta}_{u_k}(\beta_{k,i})|=1$ holds true
if and only if either $\overline{A_{u_k}(\beta_{k,j+1})} \cap A_{u_k}(\beta_{k,i})$
or $A_{u_k}(\beta_{k,j+1}) \cap \overline{A_{u_k}(\beta_{k,i})}$ occurs.
Therefore, the inequality in (\ref{P_mk_smaller_m0_v1}) takes the form
\begin{equation}
\begin{split}
P_{\thetab_k,\etab_k} \left( \hat{m}_k < m_0 \right)  \leq \sum_{j=1}^{m_0-1} \sum_{i=1}^{j} \operatorname{P}_{\thetab_k,\etab_k} \Bigg( \sum_{u_k \in \Ukd}
\left[ \ind{ \overline{A_{u_k}(\beta_{k,j+1})} \cap A_{u_k}(\beta_{k,i}) }\right.\\ +\left. \ind{A_{u_k}(\beta_{k,j+1}) \cap \overline{A_{u_k}(\beta_{k,i})}}  \right]>v_i \Bigg).
\label{P_mk_smaller_m0_v2}
\end{split}
\end{equation}
  Define the random variables $W_{u_k} = W_{u_k}(\beta_{k,j+1},\beta_{k,i}), u_k \in \Ukd,  i=1,\ldots,j, \, j=1,\ldots,m_0-1$, by
\begin{equation*} 
\begin{split}
W_{u_k} &= \ind{ \overline{A_{u_k}(\beta_{k,j+1})} \cap A_{u_k}(\beta_{k,i}) } + \ind{A_{u_k}(\beta_{k,j+1}) \cap \overline{A_{u_k}(\beta_{k,i})}}
\\ & \quad - \bigg(  \operatorname{P}_{\thetab_{u_k},\eta_{u_k}}\left( \overline{A_{u_k}(\beta_{k,j+1})} \cap A_{u_k}(\beta_{k,i})\right)
+ \operatorname{P}_{\thetab_{u_k},\eta_{u_k}}\left(A_{u_k}(\beta_{k,j+1}) \cap \overline{A_{u_k}(\beta_{k,i})}\right)\bigg),
\end{split}
\end{equation*}
and note that $\E{\thetab_{u_k},\eta_{u_k}}(W_{u_k}) = 0$  and $|W_{u_k}| \leq {2}$. 
In the \hyperref[appn]{Appendix}, we shall prove the following result
regarding the normalizing terms of the random variables  $W_{u_k}$.
\medskip

{\textbf{Lemma 1.}} \textit{For all $\etab_k=(\eta_{u_k})_{u_k\in {\cal U}_{k,d}} \in {H}^{\beta}_{k,d}$ and $\thetab_k=(\thetab_{u_k})_{u_k\in{\cal U}_{k,d}} \in \mathring\Theta_{k,d}^{\sigma}(r_{\varepsilon,k})$, all $1\leq i\leq j\leq m_0-1$, as $\varepsilon\to 0$}
\begin{equation*} 
\sum_{u_k \in \Ukd} \left( \operatorname{P}_{\thetab_{u_k},\eta_{u_k}}\left( \overline{A_{u_k}(\beta_{k,j+1})} \cap A_{u_k}(\beta_{k,i})\right)
+ \operatorname{P}_{\thetab_{u_k},\eta_{u_k}}\left(A_{u_k}(\beta_{k,j+1}) \cap \overline{A_{u_k}(\beta_{k,i})}\right) \right) = o(v_i).
\end{equation*}

\medskip
From (\ref{P_mk_smaller_m0_v2}) and Lemma 1, for all $\etab_k=(\eta_{u_k})_{u_k\in {\cal U}_{k,d}} \in {H}^{\beta}_{k,d}$ 
and $\thetab_k=(\theta_{u_k})_{u_k\in{\cal U}_{k,d}} \in \mathring\Theta_{k,d}^{\sigma}(r_{\varepsilon,k})$,
\begin{align}
&P_{\thetab_k,\etab_k} \left( \hat{m}_k < m_0 \right) \leq \sum_{j=1}^{m_0-1} \sum_{i=1}^{j} P_{\thetab_k,\etab_k} \Bigg( \sum_{u_k \in \Ukd} W_{u_k} >  \nonumber \\  
& >v_i  - \sum_{u_k \in \Ukd} \bigg[ \operatorname{P}_{\thetab_{u_k},\eta_{u_k}}\left(\overline{A_{u_k}(\beta_{k,j+1})} \cap A_{u_k}(\beta_{k,i}) \right) + \operatorname{P}_{\thetab_{u_k},\eta_{u_k}}\left({A_{u_k}(\beta_{k,j+1})} \cap \overline{A_{u_k}(\beta_{k,i})} \right)  \bigg] \Bigg) \nonumber \\
&= \sum_{j=1}^{m_0-1} \sum_{i=1}^{j} \Prob{\thetab_{k},\etab_{k}}{\sum_{u_k \in \Ukd} W_{u_k} > v_i (1+o(1))},\quad \varepsilon\to 0.
\label{P_mk_smaller_m0_v3}
\end{align}
To bound from above the probability on the right-hand side of (\ref{P_mk_smaller_m0_v3}),
we shall need the following fact.

\medskip

\textbf{Fact 1. Bernstein's first inequality.}
\textit{For $k=1,\ldots, d$, let $Y_{u_k}, u_k\in \Ukd$, be independent random variables such that {\rm (a)} $\E{} (Y_{u_k}) = 0$ for $u_k \in \Ukd$, and
{\rm (b)} for some $H>0$ and all $l \geq 2$, $\left | \E{}\left(Y^l_{u_k}\right) \right | \leq \frac{\E{}(Y_{u_k}^2)}{2}H^{l-2} l! < \infty$, $u_k\in \Ukd$.
If  conditions {\rm (a)} and {\rm (b)} hold, then, using the notation $\mathbb{S}_k=\sum_{u_k \in \Ukd} Y_{u_k}$ and $\mathbb{D}^2_{k} = \sum_{u_k \in \Ukd} \E{} (Y_{u_k}^2)$,
\begin{equation*}
\max \left\{ \Prob{}{\mathbb{S}_k \geq t}, \Prob{}{\mathbb{S}_k\leq -t} \right\} \leq  \begin{cases}
	e^{ -{t^2}/{(4\mathbb{D}_{k}^2)}}, & \text{if $0 \leq t < {\mathbb{D}_k^2}/{H}$,}\\
	e^{ -{t}/{(4H)}}, & \text{if $t \geq {\mathbb{D}^2_{k}}/{H}$}.
\end{cases}
\end{equation*}
If $|Y_{u_k}|\leq L$ almost surely for some positive constant $L$ for all $u_k \in \Ukd$,
then the above condition {\rm(b)} holds with $H=L/3$.}\\
(This version of Bernstein's inequality can be found  on pp. 164--165 of \cite{BSN}; it is also stated as Theorem~2.8 in \cite{Petrov}.)
\medskip

In order to apply Fact 1 to the probability on the right-hand side of (\ref{P_mk_smaller_m0_v3}), we first observe that, by Lemma 1, 
for all $\etab_k=(\eta_{u_k})_{u_k\in {\cal U}_{k,d}} \in {H}^{\beta}_{k,d}$ and $\thetab_k=(\thetab_{u_k})_{u_k\in{\cal U}_{k,d}} \in \mathring\Theta_{k,d}^{\sigma}(r_{\varepsilon,k})$, all $1\leq i\leq j\leq m_0-1$,
as $\varepsilon\to 0$
\begin{align*}
&\sum_{u_k \in \Ukd} {\rm E}_{\thetab_{u_k},\eta_{u_k}} (W_{u_k}^2)= \\
& \quad = \sum_{u_k \in \Ukd}
\left[ \operatorname{P}_{\thetab_{u_k},\eta_{u_k}}\left(\overline{A_{u_k}(\beta_{k,j+1})} \cap A_{u_k}(\beta_{k,i})\right)
+ \operatorname{P}_{\thetab_{u_k},\eta_{u_k}}\left({A_{u_k}(\beta_{k,j+1})} \cap \overline{A_{u_k}(\beta_{k,i})}\right)- \nonumber \right.\\
& \hspace{15mm} \left.	 - \left(\operatorname{P}_{\thetab_{u_k},\eta_{u_k}}\left(\overline{A_{u_k}(\beta_{k,j+1})} \cap A_{u_k}(\beta_{k,i})\right)
+ \operatorname{P}_{\thetab_{u_k},\eta_{u_k}}\left({A_{u_k}(\beta_{k,j+1})} \cap \overline{A_{u_k}(\beta_{k,i})}\right) \right)^2 \right] \nonumber \\
& \quad = \Bigg( \sum_{u_k \in \Ukd} \bigg[ \operatorname{P}_{\thetab_{u_k},\eta_{u_k}}\left(\overline{A_{u_k}(\beta_{k,j+1})} \cap A_{u_k}(\beta_{k,i})\right)+ \nonumber \\
& \hspace{50mm}+ \operatorname{P}_{\thetab_{u_k},\eta_{u_k}}\left({A_{u_k}(\beta_{k,j+1})} \cap \overline{A_{u_k}(\beta_{k,i})}\right) \bigg] \Bigg) (1 + o(1))  = o(v_i).
\end{align*}
Next, since ${\rm E}_{\thetab_{u_k},\eta_{u_k}}(W_{u_k}) = 0$ and $|W_{u_k}| \leq {2}$, we can apply Fact 1 with $H={2/3}$.
From~(\ref{Q2_part1}),~(\ref{P_mk_smaller_m0_v3}), and Fact 1 for $t > \mathbb{D}_{k}^2/H$, we obtain
\begin{align}
Q_{\varepsilon,k}^{(2)} & \leq \sup_{\boldsymbol{\eta}_k \in {H}_{k,d}^{\beta}} \sup_{\thetab_k \in \mathring{\Theta}^\sigma_{k,d}(r_{\varepsilon,k}) }{d \choose k}^{\beta}   \operatorname{P}_{\thetab_k,\etab_k} \left( \hat{m}_k < m_0 \right) \nonumber 
\end{align}
\begin{align}
& \leq \sup_{\boldsymbol{\eta}_k \in {H}_{k,d}^{\beta}} \sup_{\thetab_k \in \mathring{\Theta}^\sigma_{k,d}(r_{\varepsilon,k})}{d \choose k}^{\beta} \sum_{j=1}^{m_0-1} \sum_{i=1}^{j} \operatorname{P}_{\thetab_k,\etab_k}\left(\sum_{u_k \in \Ukd} W_{u_k} > v_i (1+o(1))\right)\nonumber \\ 
& \leq {d \choose k}^{\beta} \sum_{j=1}^{m_0-1} \sum_{i=1}^{j} \exp \Big( -({3v_i}/8) (1+o(1)) \Big).
\label{eq:Q2_v1}
\end{align}	
From this, by means of conditions (\ref{cond:Mk}), (\ref{vj}), and (\ref{def:tau}), and recalling that $v_1 > v_2 > \ldots > v_{M_k}$, we have as $\varepsilon \to 0$ 
\begin{eqnarray}
Q_{\varepsilon,k}^{(2)} &\leq&  {d\choose k}^\beta M_k^2 \exp \left( - \frac{1}{{4}} v_{m_0-1} \right)  
=   \exp \left( -\frac{1}{{4}\tau_{k,d}}{d\choose k}^{1-\beta_{k,m_0 - 1}}(1+o(1))\right) 
\nonumber \\ &\leq& \exp\left(-\frac{1}{{8}}{d\choose k}^{1-\beta_{k,m_0 - 1}-\epsilon/2}   \right)
= o(1),
\label{Q2_part2}
\end{eqnarray}
where the last equality follows from the fact that $1-\beta_{k,m_0 - 1}-\epsilon/2\geq 1-B-\epsilon/2= 1-B(1+o(1))>0$ for all small enough $\varepsilon$.
Finally, the substitution of (\ref{eq:Q1_bound}) and (\ref{Q2_part2}) into (\ref{def_Q1_Q2}) yields
$$R^{\varepsilon}_{k,d}(\boldsymbol{\hat{\eta}}_k) = o(1), \quad \varepsilon \to 0,$$
which completes the proof of Theorem \ref{th:sim_UB_fixedK}.
\qedsymbol
\medskip

\textit{Proof of Theorem \ref{th:sim_LB_fixedK}.}
The proof partially goes along the lines of that of Theorem 3.2 in~\cite{ST-2023}. In particular,
we can restrict ourselves to the case when
\begin{gather*}
\liminf_{\varepsilon \to 0} \frac{a_{\varepsilon,u_k}(r_{\varepsilon,k})}{\sqrt{\log {d\choose k}}} >0,
\end{gather*}
which, together with (\ref{cond:inf_LB}), gives $a_{\varepsilon,u_k}(r_{\varepsilon,k})\asymp \sqrt{\log {d\choose k}}$ as $\varepsilon\to 0.$

First, consider  the set 
\begin{gather*}
\mathbb{H}_{k,d}^{\beta}=\left\{\etab_k=(\eta_{u_k})_{u_k\in \Ukd}:\eta_{u_k}\in\{0,1\}, \;c{d\choose k}^{1-\beta}\leq \sum_{u_k\in \Ukd}\eta_{u_k}\leq
C{d\choose k}^{1-\beta}\right\},
\end{gather*}
where $0<c<1<C<\infty$. For  $1\leq k\leq d$, the value $p_k={{d\choose k}^{-\beta}}$ is roughly 
the {probability} of occurrence of nonzero components in the vector $\etab_k\in \mathbb{H}_{k,d}^{\beta}$.
Since $\mathbb{H}_{k,d}^{\beta}\subset {H}_{k,d}^{\beta},$ it follows that 
the infimum of the normalized maximum Hamming risk 
$R^{\varepsilon}_{k,d}(\boldsymbol{\tilde{\eta}}_k)$ 
 over all selectors $\boldsymbol{\tilde{\eta}}_k$ in model (\ref{model2})
satisfies 
\begin{eqnarray*}
\inf_{\boldsymbol{\tilde{\eta}}_k} R^{\varepsilon}_{k,d}(\boldsymbol{\tilde{\eta}}_k)
&=&\inf_{\boldsymbol{\tilde{\eta}}_k} \sup_{\boldsymbol{\eta}_k \in {H}_{k,d}^{\beta}} \sup_{\thetab_k\in \mathring{\Theta}^{\sigma}_{k,d}(r_{\varepsilon,k})} {d \choose k}^{\beta-1} \E{\thetab_k,\etab_k} | \boldsymbol{\tilde{\eta}}_k-\boldsymbol{\eta}_k |
\\ & \geq& 
\inf_{\boldsymbol{\tilde{\eta}}_k}\sup_{\boldsymbol{\eta}_k \in \mathbb{H}_{k,d}^{\beta}} \sup_{\thetab_k\in \mathring{\Theta}^{\sigma}_{k,d}(r_{\varepsilon,k})} {d \choose k}^{\beta-1} \E{\thetab_k,\etab_k} | \boldsymbol{\tilde{\eta}}_k-\boldsymbol{\eta}_k |.
\end{eqnarray*}
Therefore, it is sufficient to show that 
\begin{gather*}
\liminf_{\varepsilon\to 0}\inf_{\boldsymbol{\tilde{\eta}}_k}\sup_{\boldsymbol{\eta}_k \in \mathbb{H}_{k,d}^{\beta}} \sup_{\thetab_k\in \mathring{\Theta}^{\sigma}_{k,d}(r_{\varepsilon,k})} {d \choose k}^{\beta-1} \E{\thetab_k,\etab_k} | \boldsymbol{\tilde{\eta}}_k-\boldsymbol{\eta}_k |>0.
\end{gather*}
This will be done by bounding the minimax risk from below by the Bayes risk with respect to some prior distributions for 
 $\etab_k$ and $\thetab_k$.

Specifically, let the prior distributions of $\etab_k$ and $\thetab_k$ be as follows, cf. the prior distributions in Section 7.3 of~\cite{DSA-2012}:
\begin{align*}
\pi_{\etab_k} &= \prod_{u_k \in \Ukd} \pi_{\eta_{u_k}}, \quad  \pi_{\eta_{u_k}}=(1-p_k) \delta_0 + p_k \delta_1, \\
\pi_{\thetab_k} &= \prod_{u_k \in \Ukd} \pi_{\thetab_{u_k}}, \quad  \pi_{\thetab_{u_k}} = \prod_{\lb \in \mathring{\mathbb{Z}}_{u_k}} \frac{\delta_{\theta_\lb^*} + \delta_{-\theta_\lb^*}}{2},
\end{align*}
where $\theta_\lb^*=\theta_\lb^*(r_{\varepsilon,k})$ is as in (\ref{def:a2_with_theta_star}) and
$\delta_x$ is the $\delta$-measure that puts a pointmass 1 at $x$.
We assume that $\etab_k$ and $\thetab_k$ are independent. Put
${\cal X}_{k,d}=\left\{\etab_k\in \{0,1\}^{d\choose k}\right\}$ and observe that ${\rm P}_{\pi_{\etab_k}}\left({\cal X}_{k,d}\right)=1.$
Next, define the event $B_{k,d}^{\beta}=\left\{\etab_k\in \mathbb{H}_{k,d}^{\beta}\right\}\subset {\cal X}_{k,d}$, and denote by
${\rm P}_{\hat{\pi}_{\etab_k}}$ the probability measure ${\rm P}_{\pi_{\etab_k}}$ conditioned on $B_{k,d}^{\beta}$, that is,
for any $A\subseteq{\cal X}_{k,d} $,  
$${\rm P}_{\hat{\pi}_{\etab_k}}(A)= 
\frac{{ {\rm P}_{\pi_{\etab_k}}\left(A\cap B_{k,d}^{\beta}\right) }}{{ {\rm P}_{\pi_{\etab_k}}\left( B_{k,d}^{\beta}\right) }}.$$
The measure ${\rm P}_{\hat{\pi}_{\etab_k}}$ is supported on $ B_{k,d}^{\beta}$. Therefore,  
using the notation $\overline{B_{k,d}^{\beta}}={\cal X}_{k,d}\setminus B_{k,d}^{\beta}$, we can write for all $\beta\in(0,1)$ and $\sigma>0$
\begin{gather*}
\inf_{\boldsymbol{\tilde{\eta}}_k}\sup_{\etab_k\in \mathbb{H}_{k,d}^\beta}
\sup_{\thetab_k\in \mathring{\Theta}^{\sigma}_{k,d}(r_{\varepsilon,k})}{d \choose k}^{\beta-1} \E{\thetab_k,\etab_k} | \boldsymbol{\tilde{\eta}}_k-\boldsymbol{\eta}_k |
\geq	
{d \choose k}^{\beta-1}\inf_{\boldsymbol{\tilde{\eta}}_k}  
\int_{B_{k,d}^{\beta}} {\rm E}_{\pi_{\thetab_k}}
\E{\thetab_k,\etab_k} | \boldsymbol{\tilde{\eta}}_k-\boldsymbol{\eta}_k |\d{\rm P}_{\hat{\pi}_{\etab_k}}(\etab_k) \nonumber \\
= \frac{{d \choose k}^{\beta-1} }{  {\rm P}_{\pi_{\etab_k}}( B_{k,d}^{\beta})}
\inf_{\boldsymbol{\tilde{\eta}}_k} 
\int_{B_{k,d}^{\beta}} {\rm E}_{\pi_{\thetab_k}}
\E{\thetab_k,\etab_k} | \boldsymbol{\tilde{\eta}}_k-\boldsymbol{\eta}_k |\, d{\rm P}_{{\pi}_{\etab_k}}(\etab_k)\\
\geq {d \choose k}^{\beta-1}
\inf_{\boldsymbol{\tilde{\eta}}_k} \Bigg(\int_{{\cal X}_{k,d}}{\rm E}_{\pi_{\thetab_k}}
\E{\thetab_k,\etab_k} | \boldsymbol{\tilde{\eta}}_k-\boldsymbol{\eta}_k |\, d{\rm P}_{{\pi}_{\etab_k}}(\etab_k)-\int_{\overline{B_{k,d}^{\beta}}} {\rm E}_{\pi_{\thetab_k}}
\E{\thetab_k,\etab_k} | \boldsymbol{\tilde{\eta}}_k-\boldsymbol{\eta}_k |\, d{\rm P}_{{\pi}_{\etab_k}}(\etab_k) \Bigg)\\
\geq {d \choose k}^{\beta-1} \!\!\!
\inf_{\boldsymbol{\tilde{\eta}}_k} {\rm E}_{\pi_{\etab_k}}{\rm E}_{\pi_{\thetab_k}}
\E{\thetab_k,\etab_k} | \boldsymbol{\tilde{\eta}}_k-\boldsymbol{\eta}_k | -
{d \choose k}^{\beta-1}\!\!\! \sup_{\boldsymbol{\tilde{\eta}}_k}
\int_{\overline{B_{k,d}^{\beta}}}
{\rm E}_{\pi_{\thetab_k}}\E{\thetab_k,\etab_k} | \boldsymbol{\tilde{\eta}}_k-\boldsymbol{\eta}_k |\, d{\rm P}_{{\pi}_{\etab_k}}(\etab_k)\\
=:\mathbb{L}_{\varepsilon,k}^{(1)}- \mathbb{L}_{\varepsilon,k}^{(2)},
\end{gather*}
where the terms $\mathbb{L}_{\varepsilon,k}^{(1)}$ and $\mathbb{L}_{\varepsilon,k}^{(2)}$ are nonnegative.
Thus, in order to complete the proof, it suffices to show that 
\begin{gather}
{\rm (a)} \;\liminf_{\varepsilon\to 0}\mathbb{L}_{\varepsilon,k}^{(1)}
\geq {\rm const}>0\qquad\mbox{and}\qquad
{\rm (b)} \;\limsup_{\varepsilon\to 0} \mathbb{L}_{\varepsilon,k}^{(2)}=0.
\label{LB1}
\end{gather}

Aiming to verify relation (a), recall  that  $\operatorname{E}_{\thetab_{u_k},\eta_{u_k}} $ denotes the expectation which corresponds
to the measure $\operatorname{P}_{\thetab_{u_k},\eta_{u_k}}$ 
 induced by the observation $\Xb_{\!u_k} = (X_\lb)_{\lb \in \mathring{\mathbb{Z}}_{u_k}}$.
Consider the data  $\{\Xb_{u_k}\}_{u_k\in \Ukd}$ 
and introduce the following continuous mixture of distributions:
\begin{equation*}
\operatorname{P}_{\pi,\eta_{u_k}} (d\Xb_{\!u_k}) = \E{\pi_{\thetab_{u_k}}} \operatorname{P}_{\thetab_{u_k},\eta_{u_k}} (dX_\lb),\quad 
\lb \in \mathring{\mathbb{Z}}_{u_k},\; u_k\in \Ukd,
\label{LB:mixture_distr_of_thetab_and_X}
\end{equation*}
that is,
\begin{equation*}
\E{\pi_{\thetab_{u_k}}} \operatorname{E}_{\thetab_{u_k},\eta_{u_k}} |\tilde{\eta}_{u_k}-\eta_{u_k}| = \E{\pi_{\thetab_{u_k}}} \int |\tilde{\eta}_{u_k}-\eta_{u_k} | d\operatorname{P}_{\thetab_{u_k},\eta_{u_k}} = \int |\tilde{\eta}_{u_k}-\eta_{u_k}| d\operatorname{P}_{\pi,\eta_{u_k}}.
\end{equation*}
Using the definition of $\pi_{\thetab_{u_k}}$ and the fact that $X_\lb \sim N (\eta_{u_k} \theta_\lb(u_k),\varepsilon^2)$, $\lb \in \mathring{\mathbb{Z}}_{u_k},$ $u_k\in \Ukd$,  we can alternatively express this mixture of distributions as
\begin{equation*}
\operatorname{P}_{\pi,\eta_{u_k}} = \prod_{\lb \in \mathring{\mathbb{Z}}_{u_k}} \left( \frac{N(\eta_{u_k} \theta_\lb^*,\varepsilon^2) + N(- \eta_{u_k} \theta_\lb^*,\varepsilon^2)}{2} \right), \quad u_k \in \Ukd.
\end{equation*}
In particular, $\operatorname{P}_{\pi,0}=\prod_{\lb \in \mathring{\mathbb{Z}}_{u_k}} N(0,\varepsilon^2).$
Let $\nu_{\lb}^* =\nu_{\lb}^* (r_{\varepsilon,k})$ be given by
\begin{equation*}
\nu_{\lb}^* := {\theta_\lb^*}/{\varepsilon},
\label{def:nu_star2}
\end{equation*}
that is, $a_{\varepsilon,u_k}^2(r_{\varepsilon,k})=(1/2) \sum_{\lb \in \mathring{\mathbb{Z}}_{u_k}}(\nu_{\lb}^*)^4$,
and define the independent random variables
\begin{equation*}
Y_\lb := {X_\lb}/{\varepsilon} = \eta_{u_k} \nu_\lb^* +  \xi_{\lb}
\sim N(\eta_{u_k} \nu_\lb^*,1), \quad \lb \in \mathring{\mathbb{Z}}_{u_k},\quad u_k \in \Ukd.
\label{LB:Yl_normal}
\end{equation*}
Then, denoting  ${\Yb}_{\!\!u_k} = (Y_\lb)_{\lb \in \mathring{\mathbb{Z}}_{u_k}}$, we can express the likelihood ratio in the form (see
the proof of Theorem 3.2 in \cite{ST-2023})
\begin{equation}
\Lambda_{\pi,u_k} := \frac{d \operatorname{P}_{\pi,1}}{d \operatorname{P}_{\pi,0}} (\Yb_{\!\!u_k}) = \prod_{\lb \in \mathring{\mathbb{Z}}_{u_k}} \exp \left( - \frac{(\nu_{\lb}^*)^2}{2} \right) \cosh \left( (\nu_{\lb}^*)^2 Y_\lb \right),
\label{LB:likelihood_ratio}
\end{equation}
where the quantities
$\nu_{\lb}^*$, $\lb \in \mathring{\mathbb{Z}}_{u_k}$, satisfy as $\varepsilon\to 0$
\begin{gather}\label{LB:nu_star_is_o1}
\nu_{\lb}^*
=o(1).
\end{gather}
Indeed, by  (\ref{cond:inf_LB}) and the ``continuity'' of $a_{\varepsilon,u_k}$, it holds that
$r_{\varepsilon,k}/r^*_{\varepsilon,k} <1$ for all small enough~$\varepsilon$. Next,
by (\ref{thetal}) and (\ref{rr}), as $\varepsilon\to 0$
\begin{gather*}
(\nu_{\lb}^*)^2 \asymp \varepsilon^{-2} r_{\varepsilon,k}^{2+k/\sigma}\quad\mbox{and}\quad r^*_{\varepsilon,k}\asymp\left\{\varepsilon\log^{1/4} {d\choose k}   \right\}^{{4\sigma}/{(4\sigma+k)}}.
\end{gather*}
Therefore, in view of conditions (\ref{th:assump_epsilon}) and (\ref{cond:inf_LB}), we obtain as $\varepsilon\to 0$
\begin{equation*}
(\nu_{\lb}^*)^2 \asymp \varepsilon^{-2} r_{\varepsilon,k}^{2+k/\sigma}  \leq \varepsilon^{-2} (r_{\varepsilon,k}^*)^{2+k/\sigma} \asymp \varepsilon^{{2k}/{(4 \sigma + k)}} \log^{(2 \sigma +k)/{(4 \sigma + k)}} {d \choose k}  = o(1),
\end{equation*}
and thus relation (\ref{LB:nu_star_is_o1}) is verified.

Returning to the term $\mathbb{L}_{\varepsilon,k}^{(1)}$, we may continue
\begin{eqnarray}
\mathbb{L}_{\varepsilon,k}^{(1)} &=& {d \choose k}^{\beta-1} 
\inf_{\boldsymbol{\tilde{\eta}}_k} {\rm E}_{\pi_{\etab_k}}{\rm E}_{\pi_{\thetab_k}}
\E{\thetab_k,\etab_k} | \boldsymbol{\tilde{\eta}}_k-\boldsymbol{\eta}_k | 
\geq {d \choose k}^{\beta-1} \sum_{u_k \in \Ukd} \inf_{{\tilde{\eta}}_{u_k}} \E{\pi_{\eta_{u_k}}} \operatorname{E}_{\pi,\eta_{u_k}} | \tilde{\eta}_{u_k}-\eta_{u_k} |\nonumber \\ &=& {d \choose k}^{\beta-1} \sum_{u_k \in \Ukd} \inf_{{\tilde{\eta}}_{u_k}} \left( (1-p_k) \E{\pi,0}(\tilde{\eta}_{u_k}) + p_k \E{\pi,1}(1-\tilde{\eta}_{u_k}) \right),
\label{LB:risk_uneq2}
\end{eqnarray}
where $\tilde{\eta}_{u_k}$ may be viewed as a (nonrandomized) test
in the problem of testing $H_0: \operatorname{P} = \operatorname{P}_{\pi,0}$ vs. $H_1: \operatorname{P} = \operatorname{P}_{\pi,1}$, and the quantity
\begin{equation*}
\inf_{\tilde{\eta}_{u_k}} \left( (1-p_k) \E{\pi,0}(\tilde{\eta}_{u_k}) + p_k \E{\pi,1}(1-\tilde{\eta}_{u_k}) \right)
\end{equation*}
coincides with the Bayes risk in this testing problem. The infimum over $\tilde{\eta}_{u_k}$ is attained for the Bayes test $\eta_B$ defined by (see, for example, Section 8.11 of \cite{DeGr})
\begin{equation*}
\eta_B (\Yb_{\!\!u_k}) = \ind{\Lambda_{\pi,u_k} \geq \frac{1-p_k}{p_k}},
\label{LB:opt_Bayes_test}
\end{equation*}
where $\Lambda_{\pi,u_k}$ is the likelihood ratio defined in (\ref{LB:likelihood_ratio}).
It now follows from (\ref{LB:risk_uneq2}) 
that for any $u_k\in\Ukd$ (from now on, we choose some $u_k$ and fix it)
and all small enough $\varepsilon$
\begin{align}
& \mathbb{L}_{\varepsilon,k}^{(1)}  \geq {d \choose k}^{\beta-1} {d \choose k} \bigg( (1-p_k) \E{\pi,0} \left( 	\eta_B (\Yb_{\!\!u_k}) \right) + p_k \E{\pi,1} \left( 1-	\eta_B (\Yb_{\!\!u_k}) \right) \bigg) \nonumber \\
&\quad = {d \choose k}^\beta  (1-p_k) \operatorname{P}_{\pi,0} \left( \Lambda_{\pi,u_k} \geq  \frac{1-p_k}{p_k} \right) + {d \choose k}^\beta p_k \operatorname{P}_{\pi,1} \left( \Lambda_{\pi,u_k} < \frac{1-p_k}{p_k} \right)\nonumber \\
&\quad \geq \frac12 {d \choose k}^\beta  \operatorname{P}_{\pi,0} \left( \Lambda_{\pi,u_k} \geq  \frac{1-p_k}{p_k} \right) + \operatorname{P}_{\pi,1} \left( \Lambda_{\pi,u_k} < \frac{1-p_k}{p_k} \right)
=: M_{\varepsilon,k}^{(1)} + M_{\varepsilon,k}^{(2)},
\label{LB:risk_uneq3}
\end{align}
where both terms $M_{\varepsilon,k}^{(1)}$ and $M_{\varepsilon,k}^{(2)}$ are nonnegative.
Hence,  we only need to show that at least one of these terms is positive for all small enough $\varepsilon$. Under the assumptions
that $r_{\varepsilon,k}>0$ is such that 
\begin{equation*}
0<\liminf_{\varepsilon \to 0} \frac{a_{\varepsilon,u_k}(r_{\varepsilon,k})}{\sqrt{\log {d\choose k}}} \leq\limsup_{\varepsilon \to 0} \frac{a_{\varepsilon,u_k}(r_{\varepsilon,k})}{\sqrt{\log {d\choose k}}} < \sqrt{2\beta},
\label{cond:inf2_lowpart}
\end{equation*}
which is assumed in the course of the proof (see condition (\ref{nomer1}) and the condition at the beginning of the proof),
the inequality $M_{\varepsilon,k}^{(2)}>0$ holds true for all small enough $\varepsilon$.
With relation (\ref{LB:nu_star_is_o1}) being valid, this latter inequality  was verified
in the proof of Theorem 3.2 in \cite{ST-2023} (see Case 1). Specifically, it is known that for all small enough $\varepsilon$ (see p. 2027 of \cite{ST-2023})
\begin{gather*}
M_{\varepsilon,k}^{(2)}\geq {1}/{4}.
\end{gather*}
From this and (\ref{LB:risk_uneq3}), for all small enough $\varepsilon$,
$$\mathbb{L}_{\varepsilon,k}^{(1)}\geq  M_{\varepsilon,k}^{(2)}\geq {1}/{4}.  $$
This shows that relation (a) in (\ref{LB1}) holds true, and it remains to prove relation (b).

To do this, we shall need an additional fact.
\medskip

\textbf{Fact 2. Bernstein's second inequality.} \textit{If for $k=1,\ldots,d$, 
$Y_{u_k}$, $u_k\in \Ukd,$ are i.i.d. random variables with mean 0, variance $\tau^2$, and such that $|Y_{u_k}|\leq L$ for some $L>0$, $u_k\in \Ukd$, then
for all $y>0$}
\begin{gather*}
{\rm P} \left(\left|\sum_{u_k\in \Ukd}Y_{u_k}\right|\geq y  \right)\leq 2\exp\left( -\frac{y^2}{2\left(\tau^2{d\choose k}+Ly/3\right)} \right).
\end{gather*}
(This version of Bernstein's inequality can be found on page 855 of \cite{SW}.)

\medskip
Applying Fact 2 with $L=1$ and letting $\zeta=\min(C-1, 1-c)$,  we have for all large enough $d$ and an appropriate constant $Q>0$
\begin{gather*}
{\rm P}_{\pi_{\etab_k}}\left(B_{k,d}^\beta\right)={\rm P}_{\pi_{\etab_k}}\left( \etab_k\in \mathbb{H}_{k,d}^{\beta} \right) 
\end{gather*}
\begin{gather*}
={\rm P}_{\pi_{\etab_k}}\left( (c-1){d\choose k}^{1-\beta}\leq 
\sum_{u_k\in \Ukd }(\eta_{u_k}-p_k)\leq (C-1){d\choose k}^{1-\beta} \right)\\ 
\geq 1-{\rm P}_{\pi_{\etab_k}}\left(\left|\sum_{u_k\in \Ukd }(\eta_{u_k}-p_k)\right|\geq \zeta {d\choose k}^{1-\beta}\right)\\
\geq  1-2\exp\left(-\frac{\left(\zeta {d\choose k}^{1-\beta}\right)^2  }{ 2p_k(1-p_k){d\choose k}+( 1/3)\zeta{d\choose k}^{1-\beta}  }  \right)
\geq 1-2\exp\left( -Q {d\choose k}^{1-\beta}   \right). 
\end{gather*}
From this,  for all $\beta\in(0,1)$, as $d\to \infty$
$${\rm P}_{\pi_{\etab_k}}\left(\overline{B_{k,d}^\beta}\right)=1-{\rm P}_{\pi_{\etab_k}}\left(B_{k,d}^\beta\right) 
\leq 2\exp\left( -Q {d\choose k}^{1-\beta} \right)=o\left(  {d\choose k}^{-1}\right).$$ 
Finally, noting that  $|\tilde{\etab}_k-\etab_k|=\sum_{u_k\in \Ukd }|\tilde \eta_{u_k}-\eta_{u_k }|\leq {d\choose k}$ for all selectors $\tilde{\etab}_k$ 
in model (\ref{model2}) and all vectors $\etab_k$ with prior distribution $\pi_{\etab_k}$, and that $d=d_\varepsilon\to \infty$ as $\varepsilon\to 0$, we obtain for all $\beta\in(0,1)$
\begin{eqnarray*}
 \mathbb{L}_{\varepsilon,k}^{(2)}&=&
{d \choose k}^{\beta-1}\!\!\! \sup_{\boldsymbol{\tilde{\eta}}_k}
\int_{\overline{B_{k,d}^{\beta}}}
{\rm E}_{\pi_{\thetab_k}}\E{\thetab_k,\etab_k} | \boldsymbol{\tilde{\eta}}_k-\boldsymbol{\eta}_k |\, d{\rm P}_{{\pi}_{\etab_k}}(\etab_k)\\
&\leq& {d\choose k}^{\beta-1} {d\choose k} {\rm P}_{\pi_{\etab_k}}\left(\overline{ B_{k,d}^\beta}\right)= o\left({d\choose k}^{\beta-1}  \right) =o(1),\quad \varepsilon\to 0,
\end{eqnarray*}
that is, the relation (b) in (\ref{LB1}) also holds true. The proof of Theorem \ref{th:sim_LB_fixedK} is complete.
\qedsymbol
\medskip

\textit{Proof of Theorem \ref{th4}.}
The proof  follows immediately from that of Theorem \ref{th:sim_UB_fixedK} by noting that
(i) $s$ is fixed; (ii) condition (\ref{eps}) implies condition (\ref{cond_epsilon}) for all $1\leq k\leq s$;
(iii) condition (\ref{nomer0}) implies condition (\ref{cond:inf}) for all $1\leq k\leq s$;
(iv)
$\log{d\choose k} =o\left(\varepsilon^{-2k/(2\sigma+k)}  \right)$ as $\varepsilon\to 0$ is equivalent to
$\varepsilon^2\left( \log {d\choose k} \right)^{1+2\sigma/k}=o(1)$ as $\varepsilon\to 0$, and
$\max_{1\leq k\leq s}\varepsilon^2\left( \log {d\choose k} \right)^{1+2\sigma/k}=\varepsilon^2\left(\log d\right)^{1+2\sigma}$ for $d\to \infty$ and $s$ being fixed or $s\to \infty$, $s=o(d)$.
\qedsymbol

\medskip
\textit{Proof of Theorem \ref{th5}.} The proof of this theorem follows from Theorem  \ref{th:sim_LB_fixedK}.
For every $k$ ($1\leq k\leq s$), we choose some $u_k\in \Ukd$ and fix it.
Let $k^{\prime}=k^{\prime}(\varepsilon)$ be a map from $(0,\infty)$ to $\{1,\ldots,s\}$ defined as follows:
\begin{equation*}
k^{\prime}=\arg\!\min_{\!\!\!\!\!\!\!\!\!\!1\leq k\leq s}\frac{a_{\varepsilon,u_k}(r_{\varepsilon,k})}{\sqrt{2\log {d\choose k}}}.
\end{equation*}
The infimum of the maximum normalized Hamming risk ${\mathcal{R}}^{\varepsilon}_{s,d}(\boldsymbol{\tilde{\eta}})$, introduced in (\ref{def:Ham_risk_s}), over all aggregate selectors $\boldsymbol{\tilde{\eta}}$ in model (\ref{model2_s}) satisfies
\begin{align}
\inf_{\boldsymbol{\tilde{\eta}}} {\mathcal{R}}^{\varepsilon}_{s,d}(\boldsymbol{\tilde{\eta}})=
&\inf_{\boldsymbol{\tilde{\eta}}} \sup_{\boldsymbol{\eta} \in  \mathcal{H}^{\beta}_{s,d}} \sup_{\thetab\in\Theta_{s,d}^{\sigma}(r_{\varepsilon})}  \E{\thetab,\etab} \left\{\sum_{k=1}^{s}{d\choose k}^{\beta-1} \sum_{u_k \in \Ukd} |\tilde{\eta}_{u_k}-{\eta}_{u_k} |\right\}\nonumber\\
\geq&
\inf_{\boldsymbol{\tilde{\eta}}_{k^{\prime}}} \sup_{\boldsymbol{\eta}_{k^{\prime}} \in {H}_{k^{\prime},d}^{\beta}} \sup_{\thetab_{k^{\prime}}\in \mathring{\Theta}^{\sigma}_{k^{\prime},d}(r_{\varepsilon,k^{\prime}})}
{d\choose k^{\prime}}^{\beta-1}\E{\thetab_{k^{\prime}},\etab_{k^{\prime}}} \left(\sum_{u_{k^{\prime}} \in \mathcal{U}_{k^{\prime},d}} |\tilde{\eta}_{u_{k^{\prime}}}-{\eta}_{u_{k^{\prime}}}  |\right). \label{eq1}
\end{align}
Noting that the condition $\varepsilon^2\left(\log d\right)^{1+2\sigma/k}=o(1)$ as $\varepsilon\to 0$ ensures that for all $1\leq k\leq s$ one has
$\log{d\choose k} =o\left(\varepsilon^{-2k/(2\sigma+k)}  \right)$ as $\varepsilon\to 0$, we obtain from (\ref{eq1}) and Theorem \ref{th:sim_LB_fixedK}
that
$$\liminf_{\varepsilon\to 0} \inf_{\boldsymbol{\tilde{\eta}}} {\mathcal{R}}^{\varepsilon}_{s,d}(\boldsymbol{\tilde{\eta}})  >0$$
provided
\begin{gather*}
\limsup_{\varepsilon \to 0} \frac{a_{\varepsilon,u_{k^{\prime}}}(r_{\varepsilon,k^{\prime}})}{\sqrt{2\log {d\choose k^{\prime}}}} < \sqrt{\beta},
\end{gather*}
which is true by the definition of $k^{\prime}$ and condition {\rm (\ref{nomer1})}. This completes the proof. \qedsymbol

\medskip
\textit{Proof of Theorem \ref{th6}.} 
It can be easily seen that, under the conditions $s=o(\log \varepsilon^{-1})$ and $\log \log d=o(s)$, one has as $\varepsilon\to 0$
\begin{equation*}\label{rek0}
k=o\left(\log \left( \varepsilon\log^{1/4} {d\choose k} \right)^{-1} \right),\quad 1\leq k\leq s.
\end{equation*}
Therefore, as follows from relations (93) and (95) in \cite{ST-2023}, for $1\leq k\leq s$ as $\varepsilon\to 0$, cf. (\ref{rr}),
\begin{equation*} 
r_{\varepsilon,k}^*\asymp r_{\varepsilon,k,m}^*\asymp \left( \varepsilon\log^{1/4} {d\choose k} \right)^{\,{4k}/{(4\sigma+k)}} k^{-{\sigma}/{2}},
\quad m=1,\ldots,M_k,
\end{equation*}
and, cf. (\ref{omegamax_asymp}), 
\begin{equation}\label{omega0}
\max_{\lb \in \mathring{\mathbb{Z}}_{u_k}} \omega_\lb(r^*_{\varepsilon,k,m})\asymp
\left( \varepsilon\log^{1/4} {d\choose k} \right)^{\,{2k}/{(4\sigma+k)}} \left({2\pi}/{e}\right)^{{k}/{4}} k^{5/4},	\quad m=1,\ldots,M_k.
\end{equation}
The statement of the theorem is obtained by acting as in the proof of Theorem \ref{th:sim_UB_fixedK}, while using~(\ref{omega0}) instead of (\ref{omegamax_asymp}) and observing that for $1\leq k\leq d$
\begin{equation} \label{eq:approx_log_dchoosek}
{d \choose k}  \geq \left( \frac{d}{k} \right)^k.
\end{equation}
Indeed, for each $k$, let index $m_0=m_{0,k}$ ($1\leq m_{0}\leq M_{k}-1$)  be such that
\begin{equation*} \label{interval_around_beta}
\beta_{k,m_{0}} \leq \beta <\beta_{k,m_{0}+1}.
\end{equation*} 
The maximum normalized Hamming risk ${\mathcal{R}}^{\varepsilon}_{s,d}(\boldsymbol{\hat{\eta}})$ of $\hat{\boldsymbol{\eta}}=(\hat{\boldsymbol{\eta}}_1,\ldots,\hat{\boldsymbol{\eta}}_s)$ can be estimated from above by using  (\ref{def_Q1_Q2}) as follows:
\begin{align}
{\cal R}^{\varepsilon}_{s,d}(\boldsymbol{\hat{\eta}}) 
&\leq \sum_{k=1}^s 
\sup_{\etab_k\in H^{\beta}_{k,d}}\sup_{\thetab_k\in \mathring{\Theta}^\sigma_{k,d}(r_{\varepsilon,k})}{d \choose k}^{\beta-1}{\rm E}_{\thetab_k,\etab_k}\Big(| \boldsymbol{\hat{\eta}}_k-\boldsymbol{\eta}_k | \, \big| \, \hat{m}_k \geq m_{0} \Big)
\operatorname{P}_{\thetab_k,\etab_k}\left(\hat{m}_k \geq m_{0} \right) \nonumber \\
& \;\; +
\sum_{k=1}^s 
\sup_{\etab_k\in H^{\beta}_{k,d}}\sup_{\thetab_k\in \mathring{\Theta}^\sigma_{k,d}(r_{\varepsilon,k})}{d \choose k}^{\beta-1}{\rm E}_{\thetab_k,\etab_k}\Big( |\boldsymbol{\hat{\eta}}_k-\boldsymbol{\eta}_k| \, \big| \, \hat{m}_k < m_{0} \Big)
\operatorname{P}_{\thetab_k,\etab_k}\left(\hat{m}_k < m_{0} \right)\nonumber \\
& = \sum_{k=1}^s Q_{\varepsilon,k}^{(1)} + \sum_{k=1}^s Q_{\varepsilon,k}^{(2)} =: \mathcal{Q}_{\varepsilon,s}^{(1)} + \mathcal{Q}_{\varepsilon,s}^{(2)},	 \label{def_Q1s_Q2s}
\end{align}  
where $Q_{\varepsilon,k}^{(1)}$ and $Q_{\varepsilon,k}^{(2)}$ are defined in (\ref{def_Q1_Q2}). 
For the term $\mathcal{Q}_{\varepsilon,s}^{(1)}$, in view of (\ref{eq:Q1_part1}), we get 

\begin{align}
\mathcal{Q}_{\varepsilon,s}^{(1)} & \leq \sum_{k=1}^s {d \choose k}^{\beta} \operatorname{P}_{\bf 0}\left(S_{u_k}(\beta_{k,m_0}) > \sqrt{(2 \beta_{k,m_0} + \epsilon) \log{d \choose k}}\right) \nonumber \\
&	 \;\;+ C \sum_{k=1}^s \sup_{\thetab_{u_k} \in \mathring{\Theta}_{c_{u_k}}(r_{\varepsilon,k}) }   \operatorname{P}_{\thetab_{u_k}}\left(S_{u_k}(\beta_{k,m_0}) \leq  \sqrt{(2 \beta_{k,m_0} + \epsilon) \log{d \choose k}}\right) +\sum_{k=1}^s  {\tau_{k,d}}^{-1}{{d \choose k}^{{\rho}_k}} \nonumber \\
&=:  \, \sum_{k=1}^s q_{\varepsilon,k}^{(1)}  + \sum_{k=1}^s q_{\varepsilon,k}^{(2)} +\sum_{k=1}^s  {\tau_{k,d}}^{-1}{{d \choose k}^{{\rho}_k}}.
\label{eq:Q1s_part1}
\end{align}

In order to show that $ \sum_{k=1}^{s} q_{\varepsilon,k}^{(1)} = o(1)$, we note that
condition (\ref{eq:exp_bound_condT}) is satisfied, and hence inequality (\ref{eq:exp_bound}) is applicable with $T_k = \sqrt{(2 \beta_{k,m_0} + \epsilon) \log{d \choose k}}$ for each $k=1,\ldots,s$. 
Indeed, by using (\ref{omega0}),  
we obtain
\begin{align*}
T_k \max_{\lb \in \mathring{\mathbb{Z}}_{u_k}} \omega_\lb (r^*_{\varepsilon,k,m}) &\asymp  \log^{1/2}{d \choose k}  \left(\varepsilon  \log^{1/4} {d\choose k} \right)^{{2k}/{(4 \sigma +k)}}  \left({2\pi}/{e}\right)^{k/4} k ^{5/4} = o(1), 
\end{align*}
where the last equality is due to the conditions $k \leq s=o(d), \; \log \log d = o(s)$ and $s = o(\log \varepsilon^{-1})$. 
Therefore, applying the same arguments as in (\ref{eq:bound_qa}), relations (\ref{eq:d_choose_k_to_rho_const}) and (\ref{eq:approx_log_dchoosek}), and the geometric series formula, we can write for all small enough $\varepsilon$
\begin{eqnarray} \label{eq:bound_q1s}
\sum_{k=1}^{s} q_{\varepsilon,k}^{(1)} &\leq& \sum_{k=1}^{s} {d \choose k}^{\beta} \exp \left\{- \left( \beta_{k,m_{0}} +{\epsilon}/{2}\right)  \log {d \choose k} (1 + o(1))\right\} \leq \sum_{k=1}^{s} {d \choose k}^{\beta-\beta_{k,m_{0}} - \epsilon/4}   \nonumber \\
&\leq& \sum_{k=1}^s {d\choose k}^{{\rho}_k}(d/k)^{-k\epsilon/4}
= {O} \left( \sum_{k=1}^{s}\left({s}/{d}\right)^{k\epsilon/4}\right) =  {O} \left( (s/d)^{\epsilon/4}\right) =o(1),
\end{eqnarray}
where the last equality is due to (\ref{eps1}).
Next, by acting as in the derivation of (\ref{qb_Chebychev}), in view of relations (\ref{Bmin}), (\ref{omega0}), and (\ref{eq:approx_log_dchoosek}), we obtain
\begin{gather}	
\sum_{k=1}^{s}{q_{\varepsilon,k}^{(2)}}  
\leq C  \sum_{k=1}^s \exp\left(- \frac{\beta_{k,m_{0}}}{2}\left(B_k^2-1\right)^2 \log {d\choose k}\right)+
\sum_{k=1}^s \frac{C_1\max_{\lb\in  \mathring{\mathbb{Z}}_{u_k}}\omega_\lb(r_{\varepsilon,k,m_{0}}^*)}{	\inf_{\thetab_{u_k} \in \mathring{\Theta}^{(2)}_{\cb_{u_k}}(r_{\varepsilon,k})} {\rm E}_{\thetab_{u_k}}(S_{u_k}(\beta_{k,m_{0}})) }\nonumber\\
\leq C \sum_{k=1}^s {d\choose k}^{-\beta_{k,m_{0}}(B_k^2-1)^2/2}+ C_2\sum_{k=1}^s \log^{-1/2}  {d\choose k}\left( \varepsilon \log^{1/4} {d\choose k} \right) ^{{2k}/{(4\sigma+k)}} \left({2\pi}/{e}\right)^{k/4} k^{5/4} \nonumber \\
\leq C \sum_{k=1}^s (d/k)^{-k\beta_{k,m_{0}}(B_k^2-1)^2/2}  +C_2 \sum_{k=1}^s  \varepsilon^{{2k}/{(4\sigma+k)}} \log^{1/2} {d\choose k} \left({2\pi}/{e}\right)^{k/4} k^{5/4} \nonumber \\
\leq C \sum_{k=1}^s (s/d)^{kc_0/2} +  C_2 \varepsilon^{2/(4\sigma +1)}  s^{9/4}\log^{1/2} {d\choose s}  \left({2\pi}/{e}\right)^{s/4}\nonumber\\ =
O\left( (s/d)^{c_0/2} \right) + o(1)=o(1),
\label{q2s_Chebychev}
\end{gather}  
where the constant $c_0>0$ is defined in (\ref{Bmin}), and the last two equalities hold due to  the conditions imposed on $s$ and $d$ in the statement of the theorem.

Furthermore, utilizing relations (\ref{eq:d_choose_k_to_rho_const}), (\ref{def:tau_s}), and (\ref{eq:approx_log_dchoosek}), we get 
\begin{align} \label{eq:sum_tauinv}
\sum_{k=1}^{s} \tau_{k,d}^{-1} {d \choose k}^{{\rho}_k} &= O(1)  \sum_{k=1}^{s}  {d \choose k}^{-\epsilon/8}  = O(1) \sum_{k=1}^{s}  ({k}/{d})^{k\epsilon/8} =   O\left(\sum_{k=1}^{s} \left( {s}/{d} \right)^{k\epsilon/8}\right) \nonumber \\
&= O\left(\left({s}/{d}\right)^{\epsilon/8}\right)= o(1),
\end{align}
where the last equality is due to condition (\ref{eps1}).
The substitution of (\ref{eq:bound_q1s})--(\ref{eq:sum_tauinv}) into (\ref{eq:Q1s_part1}) gives
\begin{equation} \label{eq:Q1s_o1}
\mathcal{Q}_{\varepsilon,s}^{(1)} = o(1),\quad \varepsilon \to 0.
\end{equation}

Now, by arguing, for every $1\leq k\leq s$, as in the proof of Theorem \ref{th:sim_UB_fixedK}, while keeping in mind the conditions on $s$ and $d$ in the statement of the theorem, conditions (\ref{vj}), (\ref{def:tau_s}), and (\ref{eps1}), and using relation (\ref{omega0}) instead of relation (\ref{omegamax_asymp}), we can show that
 the term
$\mathcal{Q}_{\varepsilon,s}^{(2)}=\sum_{k=1}^s Q_{\varepsilon,k}^{(2)}$ on the right-hand side of  (\ref{def_Q1s_Q2s}) satisfies, cf. the first two lines of (\ref{eq:Q2_v1}),
\begin{equation*} \label{eq:Q2s_part1}
\mathcal{Q}_{\varepsilon,s}^{(2)} \leq \sup_{\etab \in {\cal H}^{\beta}_{s,d}}\sup_{\thetab \in  \Theta_{s,d}^{\sigma} (r_{\varepsilon}) } \sum_{k=1}^{s} {d \choose k}^\beta \sum_{j=1}^{m_{0}-1 } \sum_{i=1}^{j} \Prob{\thetab_{k},\etab_{k}}{\sum_{u_k \in \Ukd} W_{u_k} > v_i (1+o(1))},\quad \varepsilon\to 0.
\end{equation*}
The application of Bernstein's first inequality (see Fact 1 in the proof of Theorem \ref{th:sim_UB_fixedK}) to the probability 
on the right-hand side of the above display
yields for all small enough $\varepsilon$,
cf. (\ref{eq:Q2_v1})--(\ref{Q2_part2}),
\begin{gather*}
\mathcal{Q}_{\varepsilon,s}^{(2)}  \leq  \sum_{k=1}^{s} {d \choose k}^{\beta} \sum_{j=1}^{m_0-1} \sum_{i=1}^{j} \exp \left( -({3v_i}/{8}) (1+o(1)) \right) \leq \sum_{k=1}^{s} {d \choose k}^{\beta} M_k^2 \exp \left( -(1/4) v_{m_{0}-1} \right) \nonumber \\
\leq \sum_{k=1}^{s} {d \choose k}^{\beta} M_k^2 \exp \left( -\frac{1}{4 \tau_{k,d}} {d \choose k}^{1-\beta_{k,m_{0}-1}} \right) \leq \sum_{k=1}^{s} \exp \left( -\frac{1}{4}  {d \choose k}^{1-\beta_{k,m_{0}-1}-\epsilon/2}(1+o(1)) \right) \nonumber \\
\leq \sum_{k=1}^s\exp\left(- \frac{1}{8} {d \choose k}^{1-\beta_{k,m_{0}-1}-\epsilon/2} \right).
\end{gather*}
From this, noting that  for all $1\leq k\leq s$ and all small enough $\varepsilon$, $1-\beta_{k,m_{0}-1}-\epsilon/2\geq 1-B-\epsilon/2
= 1-B(1+o(1))>0$, we may continue:
\begin{gather}
\mathcal{Q}_{\varepsilon,s}^{(2)} \leq
s\exp\left( -(1/8)  d^{1-B(1+o(1))}\right)=o(d)\exp\left( -(1/8)  d^{1-B(1+o(1))} \right)= o(1).
\label{eq:Q2s_part2}
\end{gather}

Finally, inserting relations (\ref{eq:Q1s_o1}) and (\ref{eq:Q2s_part2}) into (\ref{def_Q1s_Q2s}),  we get
\begin{align*}
{\cal R}^{\varepsilon}_{s,d}(\boldsymbol{\hat{\eta}}) \leq  \mathcal{Q}_{\varepsilon,s}^{(1)} + \mathcal{Q}_{\varepsilon,s}^{(2)} = o(1),	\quad \varepsilon\to 0,
\end{align*}
and the proof of Theorem \ref{th6} is complete.   \qedsymbol

\medskip
\textit{Proof of Theorem \ref{th7}.} The proof is similar to that of Theorem \ref{th5}, which in turn builds on the proof of Theorem \ref{th:sim_LB_fixedK}.
We only notice that relation (\ref{LB:nu_star_is_o1}), a key relation in the
proof of the lower bound on the normalized minimax risk in Theorem \ref{th:sim_LB_fixedK}, continues to hold.   \qedsymbol

\begin{appendix}
\section*{Appendix}\label{appn} %

\textit{Proof of Lemma 1.} The proof of this lemma goes, in part, along the lines of that of relation (56) in \cite{BSt-2017}. 
By the definition of $A_{u_k}(\beta)$ and taking into account the sparsity condition (\ref{sparsitycond}), 
for all $\etab_k=(\eta_{u_k})_{u_k\in {\cal U}_{k,d}} \in {H}^{\beta}_{k,d}$ and $\thetab_k=(\thetab_{u_k})_{u_k\in{\cal U}_{k,d}} \in \mathring\Theta_{k,d}^{\sigma}(r_{\varepsilon,k})$, 
all $1\leq i\leq j\leq m_0-1$,
and all small enough~$\varepsilon$, we obtain
\begin{align}
	&\sum_{u_k \in \Ukd} \left[ \operatorname{P}_{\thetab_{u_k},\eta_{u_k}}\left(\overline{A_{u_k}(\beta_{k,j+1})} \cap A_{u_k}(\beta_{k,i})\right) +  \operatorname{P}_{\thetab_{u_k},\eta_{u_k}}\left(A_{u_k}(\beta_{k,j+1}) \cap \overline{A_{u_k}(\beta_{k,i})} \right) \right]   \nonumber \\
	&= \sum_{u_k: \eta_{u_k} = 0} \left[  \operatorname{P}_{\thetab_{u_k},\eta_{u_k}}\left(\overline{A_{u_k}(\beta_{k,j+1})} \cap A_{u_k}({\beta_{k,i}})\right) +  \operatorname{P}_{\thetab_{u_k},\eta_{u_k}}\left(A_{u_k}(\beta_{j+1}) \cap \overline{A_{u_k}(\beta_{k,i})}\right) \right] \nonumber \\
	&\quad  + \sum_{u_k: \eta_{u_k} = 1} \left[ \operatorname{P}_{\thetab_{u_k},\eta_{u_k}}\left(\overline{A_{u_k}(\beta_{k,j+1})} \cap A_{u_k}(\beta_{k,i})\right) + \operatorname{P}_{\thetab_{u_k},\eta_{u_k}}\left(A_{u_k}(\beta_{k,j+1}) \cap \overline{A_{u_k}(\beta_{k,i})} \right) \right]  \nonumber
	\end{align}
	\begin{align}
	&\leq   {d\choose k} \left[ \operatorname{P}_{\bf 0}\left(\overline{A_{u_k}(\beta_{k,j+1})} \cap A_{u_k}(\beta_{k,i})\right) +
	\operatorname{P}_{\bf 0}\left({A_{u_k}(\beta_{k,j+1})} \cap \overline{A_{u_k}(\beta_{k,i})}\right) \right] + \nonumber \\
	& + C {d\choose k}^{1-\beta} \sup_{\thetab_{u_k} \in \mathring{\Theta}_{\cb_{u_k}}(r_{\varepsilon,k})} \bigg[ \operatorname{P}_{\thetab_{u_k}}\left(\overline{A_{u_k}(\beta_{k,j+1})} \cap A_{u_k}(\beta_{k,i})\right) 
	+ \operatorname{P}_{\thetab_{u_k}}\left(A_{u_k}(\beta_{k,j+1}) \cap \overline{A_{u_k}(\beta_{k,i})}\right) \bigg] \nonumber \\
	&\leq {d\choose k} \left[ \operatorname{P}_{\bf 0}\left(\overline{A_{u_k}(\beta_{k,j+1})} \right) + \operatorname{P}_{\bf 0}\left(\overline{A_{u_k}(\beta_{k,i})} \right) \right] + \nonumber \\
	& \quad +  C {d\choose k}^{1-\beta}\!\!\!\!\! \sup_{\thetab_{u_k} \in \mathring{\Theta}_{\cb_{u_k}}(r_{\varepsilon,k})}\!\! \left[ \operatorname{P}_{\thetab_{u_k}}\left(A_{u_k}(\beta_{k,i})\right) + \operatorname{P}_{\thetab_{u_k}}\left(A_{u_k}(\beta_{k,j+1})\right) \right] \nonumber \\
	\begin{split}
		&\leq   {d\choose k} \left[ \operatorname{P}_{\bf 0}\left(S_{u_k}(\beta_{k,j+1}) > \sqrt{(2\beta_{k,j+1} + \epsilon) \log {d\choose k}} \right) +
		\operatorname{P}_{\bf 0}\left( S_{u_k}(\beta_{k,i}) > \sqrt{(2 \beta_{k,i} + \epsilon) \log {d\choose k}}\right) \right]   \\
		&\quad  + C  {d\choose k}^{1-\beta} \sup_{\thetab_{u_k} \in \mathring{\Theta}_{\cb_{u_k}}(r_{\varepsilon,k})} \Bigg[ \operatorname{P}_{\thetab_{u_k}}\left(S_{u_k}(\beta_{k,j+1}) \leq \sqrt{(2\beta_{k,j+1} + \epsilon) \log {d\choose k}}\right)+   \\ 
		&\hspace{10mm} +\operatorname{P}_{\thetab_{u_k}}\left(S_{u_k}(\beta_{k,i}) \leq \sqrt{(2\beta_{k,i} + \epsilon) \log {d\choose k}}\right) \Bigg] =: J_{\varepsilon,k}^{(1)} (\beta_{k,j+1},\beta_{k,i}) + J_{\varepsilon,k}^{(2)}(\beta_{k,j+1},\beta_{k,i}).
		\label{def:J1_J2}
	\end{split}
\end{align}

Consider the first term on the right-hand side of (\ref{def:J1_J2}). Applying   (\ref{vj}), (\ref{def:tau}), and (\ref{eq:exp_bound}),
as~$\varepsilon\to 0$
\begin{gather*}
	{d\choose k} \operatorname{P}_{\bf 0}\left({S_{u_k}(\beta_{k,i}) > \sqrt{(2\beta_{k,i} + \epsilon) \log {d\choose k}} }\right)  \leq {d\choose k} \exp\left\{-\left( \beta_{k,i} + {\epsilon}/{2}\right)  \log {d\choose k}(1+o(1)) \right\}
	\nonumber \\ = {O} \left({d\choose k}^{1-\beta_{k,i}-{\epsilon}/{2}}\right) = O \left( v_i \tau_{k,d} {d\choose k}^{-{\epsilon}/{2}}\right) = o(v_i).
	\label{J1_part1}
\end{gather*}
Similarly, since $v_1 > v_2 > \ldots >v_{M_k}$ and $i < j+1$, utilizing inequality (\ref{eq:exp_bound}) gives as $\varepsilon \to 0$
\begin{multline*}
	{d\choose k} \operatorname{P}_{\bf 0}\left({S_{u_k}(\beta_{k,j+1}) > \sqrt{(2\beta_{k,j+1} + \epsilon) \log {d\choose k}} }\right)
	= O \left({d\choose k}^{1-\beta_{k,j+1}-{\epsilon}/{2}}\right)  = o(v_{j+1})= o(v_i).
	\label{J1_part2}
\end{multline*}
Combining the last two relations, we obtain for all $1\leq i\leq j\leq m_0-1$
\begin{equation} \label{J1_ovi}
	J_{\varepsilon,k}^{(1)} (\beta_{k,j+1},\beta_{k,i}) = o(v_i),\quad \varepsilon\to 0.
\end{equation}

Next, consider the term $J_{\varepsilon,k}^{(2)}(\beta_{k,j+1},\beta_{k,i})$ on the right-hand side of (\ref{def:J1_J2}). 
Fist, we note that for $1\leq i\leq j\leq m_0-1$
\begin{equation} \label{beta_rels}
	\beta_{k,j+1}  \leq \beta_{k,m_0} \leq  \beta \quad \mbox{and} \quad \beta_{k,i}  \leq \beta_{k,m_0-1} < \beta.
\end{equation}
By the definition of $r^*_{\varepsilon,k,m},$ $m=1,\ldots,M_k$,
as in (\ref{def:r_star_ekm}), and  (\ref{beta_rels}), 
we have for
$1\leq i\leq j\leq m_0-1$
\begin{align*}
	a_{\varepsilon,k}(r^*_{\varepsilon,k,m})&= \sqrt{2 \beta_{k,m} \log {d\choose k}} \leq \sqrt{2 \beta \log {d\choose k}},\quad m=i,j+1.
	\label{a_r_star_ij1_smaller_than _sqrt}
\end{align*}  
From this and the monotonicity of $a_{\varepsilon,k}(r_{\varepsilon,k})$, one can find constants $\Delta_{3,k}\in(0,1)$ and $\Delta_{4,k}\in(0,1)$ such that for all small enough $\varepsilon$, every $r_{\varepsilon,k}$ that satisfies (\ref{cond:inf}) also satisfies for 
all $1\leq i\leq j\leq m_0-1$
\begin{equation*} 
	r_{\varepsilon,k}  \geq r^*_{\varepsilon,k,i}(1+\Delta_{3,k}) ,  \quad \quad r_{\varepsilon,k}  \geq r^*_{\varepsilon,k,j+1}(1+\Delta_{4,k}).
\end{equation*}
Hence, acting as in the derivation of (\ref{bound_infES_knownBeta}), we obtain that for all small enough $\varepsilon$
\begin{equation} \label{infES_larger_than_sqrt}
	\inf_{\thetab_{u_k} \in \mathring{\Theta}_{\cb_{u_k}}(r_{\varepsilon,k})} {\rm E}_{\thetab_{u_k}}(S_{u_k}(\beta_{k,i})) \geq \sqrt{2 \beta_{k,i}\log {d\choose k}}(1+\Delta_{3,k})^2.
\end{equation}
 Then, applying similar arguments for estimating $\inf_{\thetab_{u_k} \in \mathring{\Theta}_{\cb_{u_k}}(r_{\varepsilon,k})} {\rm E}_{\thetab_{u_k}}(S_{u_k}(\beta_{k,j+1}))$,
we get for $1\leq i\leq j\leq m_0-1$, as $\varepsilon\to 0$ 
\begin{equation*}
	\sqrt{(2 \beta_{k,m}+\epsilon)\log {d\choose k}} -	\inf_{\thetab_{u_k} \in \mathring{\Theta}_{\cb_{u_k}}(r_{\varepsilon,k})} {\rm E}_{\thetab_{u_k}}(S_{u_k}(\beta_{k,m}))  \to -\infty, \quad m=i,j+1.
\end{equation*}

Now, consider
	the term $J_{\varepsilon,k}^{(2)}(\beta_{k,j+1},\beta_{k,i})$ on the right-hand side of (\ref{def:J1_J2}).
	We have
\begin{align}
	&C {d\choose k}^{1-\beta} \sup_{\thetab_{u_k} \in \mathring{\Theta}_{\cb_{u_k}}(r_{\varepsilon,k})} 	
	\operatorname{P}_{\thetab_{u_k}}\left(S_{u_k}(\beta_{k,i}) \leq \sqrt{(2\beta_{k,i} + \epsilon) \log {d\choose k}}\right) \nonumber  \\
	&\quad  \leq C {d\choose k}^{1-\beta} \sum_{p=1}^2 \sup_{\thetab_{u_k} \in \mathring{\Theta}^{(p)}_{\cb_{u_k}}(r_{\varepsilon,k})} 	
	\operatorname{P}_{\thetab_{u_k}}\left(S_{u_k}(\beta_{k,i}) \leq \sqrt{(2\beta_{k,i} + \epsilon) \log {d\choose k}}\right) \nonumber \\
	& \quad  \leq C {d\choose k}^{1-\beta} \!\!\!\sup_{\thetab_{u_k} \in \mathring{\Theta}^{(1)}_{\cb_{u_k}}(r_{\varepsilon,k})} 	\!\!
	\operatorname{P}_{\thetab_{u_k}}\Bigg(S_{u_k}(\beta_{k,i})  -  {\rm E}_{\thetab_{u_k}}(S_{u_k}(\beta_{k,i})) \leq \nonumber \\
	& \hspace{53mm} \leq \sqrt{(2\beta_{k,i} + \epsilon) \log {d\choose k}}  - \inf_{\thetab_{u_k} \in \mathring{\Theta}^{(1)}_{\cb_{u_k}}(r_{\varepsilon,k})}\!\!{\rm E}_{\thetab_{u_k}}(S_{u_k}(\beta_{k,i})) \Bigg)  \nonumber \\
	& \qquad + C {d\choose k}^{1-\beta} \sup_{\thetab_{u_k} \in \mathring{\Theta}^{(2)}_{\cb_{u_k}}(r_{\varepsilon,k})}
	\operatorname{P}_{\thetab_{u_k}}\Bigg(S_{u_k}(\beta_{k,i}) - {\rm E}_{\thetab_{u_k}}(S_{u_k}(\beta_{k,i})) \leq \nonumber \\
	& \hspace{53mm} \leq \sqrt{(2\beta_{k,i} + \epsilon) \log {d\choose k}}- {\rm E}_{\thetab_{u_k}}(S_{u_k}(\beta_{k,i})) \Bigg) \nonumber\\
	&\quad =: L_{\varepsilon,k,i}^{(1)} +  L_{\varepsilon,k,i}^{(2)}.
	\label{def:L1_L2_L3}
\end{align}

Let us show that each term on the right-hand side of (\ref{def:L1_L2_L3}) is $o(v_i)$ as $\varepsilon\to 0$. First, consider the term $L_{\varepsilon,k,i}^{(1)}$.
Using (\ref{infES_larger_than_sqrt}) and the relation $\sqrt{1+x} \sim 1 + {x}/{2}$ as $x \to 0$, we get as $\varepsilon\to 0$
\begin{align*}
	& \inf_{\thetab_{u_k} \in \mathring{\Theta}^{(1)}_{\cb_{u_k}}(r_{\varepsilon,k})} {\rm E}_{\thetab_{u_k}}(S_{u_k}(\beta_{k,i})) \geq \inf_{\thetab_{u_k} \in \mathring{\Theta}_{\cb_{u_k}}(r_{\varepsilon,k})}  {\rm E}_{\thetab_{u_k}}(S_{u_k}(\beta_{k,i})) \geq   \\ 
	&\geq \sqrt{2 \beta_{k,i}\log {d\choose k}}(1+\Delta_{3,k})^2 >\sqrt{(2 \beta_{k,i}+\epsilon)\log {d\choose k}} = \sqrt{2 \beta_{k,i}\log {d\choose k}}(1+o(1)),
\end{align*}  
and hence
\begin{eqnarray*} \label{def:T2_in_L1}
	\mathbb{T}_{k,i} &: =&  \inf_{\thetab_{u_k} \in \mathring{\Theta}^{(1)}_{\cb_{u_k}}(r_{\varepsilon,k})} {\rm E}_{\thetab_{u_k}}(S_{u_k}(\beta_{k,i}))  - \sqrt{(2 \beta_{k,i}+\epsilon)\log {d\choose k}}\nonumber  \\ &\geq &
	\sqrt{2 \beta_{k,i}\log {d\choose k}} \left[ (1+\Delta_{3,k})^2 -  \left(1+o(1)\right) \right].
\end{eqnarray*}  
From this, by applying inequality (\ref{eq:exp_bound2}) with $T_k=\mathbb{T}_{k,i}$ as above,
recalling 
	(\ref{vj}) and (\ref{def:tau}), noting that for all  $1\leq i\leq m_0-1$, we have $1-\beta < 1-\beta_{k,i}$ and $b\leq \beta_{k,i}\leq B$,
and setting  $\Delta_{5,k} = (1+\Delta_{3,k})^2 -1>0 $,
we obtain for all small enough $\varepsilon$
\begin{align}
	L_{\varepsilon,k,i}^{(1)} &\leq C{d \choose k}^{1-\beta} \exp \left( -\beta_{k,i}\log {d\choose k} \left[ (1+\Delta_{3,k})^2 -  1+o(1) \right]^2 (1+o(1)) \right) \nonumber \\
	& < C{d \choose k}^{1-\beta_{k,i}} \exp \left( - \beta_{k,i}\log {d\choose k} \left[ \Delta_{5,k}+o(1) \right]^2 (1+o(1)) \right) \nonumber \\
	&\leq C{d\choose k}^{1-\beta_{k,i}-\beta_{k,i} \Delta^2_{5,k}/2} = C v_i \tau_{k,d} {d\choose k}^{-\beta_{k,i} \Delta^2_{5,k}/2} = o(v_i).
	\label{L1_ovi}
\end{align}

Let us now bound from above the term $L_{\varepsilon,k,i}^{(2)} $ on the right-hand side of (\ref{def:L1_L2_L3}).
Using Chebyshev's inequality, the definition of $\mathring{\Theta}^{(2)}_{\cb_{u_k}}(r_{\varepsilon,k})$, and relations 
(\ref{omegamax_asymp}), (\ref{VarS_knownBeta}), (\ref{bound_infES_knownBeta}), (\ref{eq:log_dk_is_oE})
and (\ref{beta_rels}),  for all small enough $\varepsilon$ and some positive constants $C_3$ and $C_4$, the same for all $1\leq k\leq d$, 
we have for $1\leq i\leq m_0-1$
\begin{align}
	L_{\varepsilon,k,i}^{(2)} &= C{d\choose k}^{1-\beta}\!\!\!\!\!\! \!\!\!\sup_{\thetab_{u_k} \in \mathring{\Theta}^{(2)}_{\cb_{u_k}}(r_{\varepsilon,k})} \!\!
	\operatorname{P}_{\thetab_{u_k}}
	\Bigg(S_{u_k}(\beta_{k,i}) -  {\rm E}_{\thetab_{u_k}}(S_{u_k}(\beta_{k,i}))\leq  \nonumber  \\
	& \hspace{5cm} \leq  \sqrt{(2\beta_{k,i} + \epsilon) \log {d\choose k}}- {\rm E}_{\thetab_{u_k}}(S_{u_k}(\beta_{k,i}))  \Bigg) \nonumber \\
	& \leq C {d \choose k}^{1-\beta}  \sup_{\thetab_{u_k}\in \mathring{\Theta}^{(2)}_{\cb_{u_k}}(r_{\varepsilon,k})}
	\frac{{\rm var}_{\thetab_{u_k}} (S_{u_k}(\beta_{k,i}))}{\left( {\rm E}_{\thetab_{u_k}}(S_{u_k}(\beta_{k,i})) - \sqrt{(2\beta_{k,i} + \epsilon) \log {d\choose k}}\right)^2 }\nonumber\\
	&\leq C {d \choose k}^{1-\beta}  \sup_{\thetab_{u_k}\in \mathring{\Theta}^{(2)}_{\cb_{u_k}}(r_{\varepsilon,k})      }\frac{  1+ 
		4\max_{\lb\in \mathring{\mathbb{Z}}_{u_k}}\omega_\lb(r_{\varepsilon,k,i}^*){\rm E}_{\thetab_{u_k}}(S_{u_k}(\beta_{k,i}))  }{\left({\rm E}_{\thetab_{u_k}}(S_{u_k}(\beta_{k,i}))- \sqrt{(2\beta_{k,i} + \epsilon) \log {d\choose k}}\right)^2} \nonumber \\
	&\leq C_3 {d \choose k}^{1-\beta_{k,i}}  \sup_{\thetab_{u_k}\in \mathring{\Theta}^{(2)}_{\cb_{u_k}}(r_{\varepsilon,k})      }\frac{  \max_{\lb\in \mathring{\mathbb{Z}}_{u_k}}\omega_\lb(r_{\varepsilon,k,i}^*)}{{\rm E}_{\thetab_{u_k}}(S_{u_k}(\beta_{k,i}))  }\leq
	\frac{ C_3 {d \choose k}^{1-\beta_{k,i}} \max_{\lb\in \mathring{\mathbb{Z}}_{u_k}}\omega_\lb(r_{\varepsilon,k,i}^*)}{\inf_{\thetab_{u_k} \in \mathring{\Theta}_{\cb_{u_k}}(r_{\varepsilon,k})} {\rm E}_{\thetab_{u_k}}(S_{u_k}(\beta_{k,i}))  }\nonumber 
\end{align}
\begin{align}
	&\leq	C_4 {d \choose k}^{1-\beta_{k,i}} \log^{-{1}/{2}} {d \choose k }
	\left( \varepsilon \log^{{1}/{4}} {d\choose k} \right)^{{2k}/{(4 \sigma + k)}}  \nonumber\\
	&=  C_4 v_i \tau_{k,d} \log^{-{1}/{2}} {d \choose k }  \left( \varepsilon \log^{{1}/{4}} {d\choose k} \right)^{{2k}/{(4 \sigma + k)}}
	=o(v_i),\quad \varepsilon\to 0, \label{L2_part1} 
\end{align}
where the last equality follows from conditions (\ref{cond_epsilon}), (\ref{def:tau}), and (\ref{th:assump_epsilon}). 

By further analyzing  the term $J_{\varepsilon,k}^{(2)}(\beta_{k,j+1},\beta_{k,i})$ on the right-hand side of (\ref{def:J1_J2}),
we note that  $v_1 > v_2 > \ldots >v_{M_k}$ and $1-\beta\leq 1-\beta_{k,j+1}$ for all $1\leq j\leq m_0-1$. 
Therefore, applying the same arguments as above, for all  $1\leq i\leq j\leq m_0-1$, we arrive at, as $\varepsilon\to 0$,
\begin{equation} \label{J2_part1_ovi}
	C {d\choose k}^{1-\beta} \sup_{\thetab_{u_k} \in \mathring{\Theta}_{\cb_{u_k}}(r_{\varepsilon,k})}
	\operatorname{P}_{\thetab_{u_k}}
	\left(S_{u_k}(\beta_{k,j+1}) \leq \sqrt{(2\beta_{k,j+1} + \epsilon) \log {d\choose k}}\right) = o(v_{j+1}) = o(v_i).
\end{equation}
Combining (\ref{def:L1_L2_L3})--(\ref{J2_part1_ovi}) yields for $1\leq i\leq j\leq m_0-1$
\begin{equation*} \label{J2_ovi}
	J_{\varepsilon,k}^{(2)} (\beta_{k,j+1},\beta_{k,i}) = o(v_i), \quad \varepsilon \to 0.
\end{equation*}
From this, (\ref{def:J1_J2}), and (\ref{J1_ovi}), the statement of Lemma 1 
follows. \qedsymbol
\end{appendix}

\section*{Acknowledgments}
The authors would like to thank the anonymous referees, the Associate
Editor, and the Editor for their constructive comments  and suggestions that helped to improve the
quality of this article.

\section*{Funding}
The research of N. A. Stepanova was supported by an NSERC grant. The research of M.~Turcicova and X. Zhao was partially supported by NSERC grants. Additionally, the research of M. Turcicova was funded by the project ``Research of Excellence on Digital Technologies and Wellbeing CZ.02.01.01/00/22\_008/0004583" which is co-financed by the European Union.




\end{document}